\documentclass[11pt]{CSM_template}
\usepackage{amssymb,amsmath,amsthm,color,psfrag}
\usepackage{graphicx,epsfig}
\usepackage{algorithm,algorithmic,xspace}
\usepackage[sort,numbers]{natbib}
\usepackage{indentfirst}
\usepackage{subfigure}
\usepackage{multibib}
\newcites{jump}{\small References}
\newcites{sol}{\small References}
\newcites{add}{\small References}
\newcites{lip}{\small References}
\newcites{time}{\small References}

\usepackage{makeidx}
\usepackage{makesdx}
\makeindex
\makesymbolindex

\usepackage{vmargin}
\setpapersize{USletter}
\setmarginsrb{1in}{1in}{1in}{1in}%
           {9pt}{9pt}{9pt}{18pt}

\newcommand{\fmax}{f_{\operatorname{max}}}
\newcommand{\fmin}{f_{\operatorname{min}}}
\newcommand{\Imax}{I_{\operatorname{max}}}
\newcommand{\Imin}{I_{\operatorname{min}}}
\newcommand{\sign}{\operatorname{sign}}
\newcommand{\parts}[1]{\mathfrak{B}(#1)}
\newcommand{\setvaluedmap}{\mathcal{F}}
\newcommand{\domain}{\mathcal{D}}

\renewcommand{\SS}{\mathcal{S}}
\newcommand{\Uc}{\mathcal{U}}
\newcommand{\nearestneighbor}{\mathcal{N}}

\newcommand{\HHSP}{\mathcal{H}_{\text{SP}}}

\newcommand{\D}{\operatorname{dist}}
\newcommand{\sm}{\operatorname{sm}}
\newcommand{\normal}{\operatorname{n}}

\newcommand{\real}{{{\mathbb{R}}}}
\newcommand{\realpositive}{(0,\infty)}
\newcommand{\realnonnegative}{[0,\infty)}
\newcommand{\realnonpositive}{(-\infty,0]}
\newcommand{\rational}{{\mathbb{Q}}}

\newcommand{\co}{\operatorname{co}}
\newcommand{\union}{\cup}
\newcommand{\intersection}{\ensuremath{\operatorname{\cap}}}

\newcommand{\argmin}{\ensuremath{\operatorname{argmin}}}
\newcommand{\until}[1]{\{1,\dots,#1\}}

\newcommand{\boundary}[1]{\operatorname{bndry}(#1)}
\newcommand{\LN}{\operatorname{Ln}} 
\newcommand{\Lc}{{\mathcal{L}}}

\newcommand\oprocendsymbol{\hbox{$\blacksquare$}}
\newcommand\oprocend{\relax\ifmmode\else\unskip\hfill\fi\oprocendsymbol}

\newcounter{sidebar}
\newcounter{sidebarequation}
\newcounter{sidebarproposition}
\newcounter{example}
\newcommand{\cball}[2]{\overline{B}(#2,#1)}
\newcommand{\oball}[2]{B(#2,#1)}

\newcommand{\setLiederG}[2]{\widetilde{\Lc}_{#1} #2}

\newcommand{\setLiederlow}[2]{\underline{\Lc}_{#1} #2}
\newcommand{\setLiederup}[2]{\overline{\Lc}_{#1} #2}

\newcommand{\setdef}[2]{\{#1 \; : \; #2\}}

\newcommand{\epigraph}[1]{\operatorname{epi}(#1)}
\newcommand{\ov}{\overline}
\newcommand{\eps}{\varepsilon}

\newcommand{\diam}{\operatorname{diam}}

\renewcommand{\natural}{{\mathbb{N}}}

\newcommand{\map}[3]{#1: #2 \rightarrow #3}
\newtheorem{theorem}{\hspace*{1.8cm} Theorem}
\newtheorem{proposition}{\hspace*{1.8cm} Proposition}
\newtheorem{corollary}{\hspace*{1.8cm} Corollary}

\newcommand{\TwoNorm}[1]{\|#1\|_{2}}
\newcommand{\TwoNormBig}[1]{\Big\|#1\Big\|_{2}}

\renewcommand{\theenumi}{(\roman{enumi})}

\begin{document}

\title{Discontinuous Dynamical Systems\\
  {\large A tutorial on solutions, nonsmooth analysis, and stability}}


\author{Jorge Cort\'es}

\date{\today}

\maketitle

Discontinuous dynamical systems arise in a large number of
applications. In optimal control problems, open-loop bang-bang
controllers switch discontinuously between extreme values of the
inputs to generate minimum-time trajectories from the initial to the
final states~\cite{AAA-YS:04}. Thermostats implement on-off
controllers to regulate room temperature~\cite{SB:93}.  When the
temperature is above the desired value, the controller switches the
cooling system on. Once the temperature reaches a preset value, the
controller switches the cooling system off.  The controller is
therefore a discontinuous function of the room temperature.  In
nonsmooth mechanics, the motion of rigid bodies is subject to velocity
jumps and force discontinuities as a result of friction and
impact~\cite{BB:98,FP-CG:96}.  In the robotic manipulation of objects
by means of mechanical contact~\cite{GASP-MFMM-VK:04}, discontinuities
occur naturally from interaction with the environment.

Discontinuities are also intentionally designed to achieve regulation
and stabilization.  Sliding mode control~\cite{VIU:92,CE-SKS:98} uses
discontinuous feedback controllers for stabilization.  The design
procedure for sliding mode control begins by identifying a surface in
the state space with the property that the dynamics of the system
restricted to this surface are easily stabilizable. Feedback
controllers are then synthesized on each side of the surface to steer
the solutions of the system toward the surface. The resulting
closed-loop system is discontinuous on the surface.  In
robotics~\cite{RCA:98}, it is of interest to induce emergent behavior
in a swarm of robots by prescribing interaction rules among individual
agents.  Simple laws such as ``move away from the nearest other robot
or environmental boundary'' give rise to discontinuous dynamical
systems. For example, consider a robot placed in a corner of a room.
On opposite sides of the bisector line between the two walls forming
the corner, the ``move away'' law translates into different velocity
vectors for the robot. The dynamical system is thus discontinuous
along the bisector line.  In optimization~\cite{FHC:83,AB-EK:06},
continuous-time algorithms that perform generalized gradient descent,
which is a discontinuous function of the state, are used when the
objective function is not smooth.  In adaptive
control~\cite{MMP-PAI:93}, switching algorithms are employed to select
the most appropriate controller from a given finite family in order to
enhance robustness, ensure boundedness of the estimates, and prevent
the system from stepping into undesired regions of the state space.

Many control systems cannot be stabilized by continuous
state-dependent feedback.  As a consequence, it is necessary to
consider either time-dependent or discontinuous feedback.  As an
illustration, consider the one-dimensional system $\dot x = X(x,u)$
with an equilibrium at the origin for $u=0$. When $x>0$, we look for
$u$ such that $X(x,u)<0$ (``go left'').  Likewise, when $x<0$, we look
for $u$ such that $X(x,u)>0$ (``go right''). Together, these
conditions can be stated as $x X(x,u)<0$.  When the control $u$
appears nonlinearly in $X$, it may be the case that no continuous
state-dependent feedback controller $x \mapsto u(x)$ exists such that
$x X(x,u(x))<0$ for all $x \in \real$. The following example
illustrates the above discussion.

\refstepcounter{example}%
\subsection*{Example~\theexample: System requiring discontinuous
  stabilization}\label{ex:sontag}

Consider the one-dimensional system~\cite{EDS:98}
\begin{align}\label{eq:sontag}
  \dot x = x [(u-1)^2 - (x-1)] [(u+1)^2 + (x-2) ] .
\end{align}
The shaded areas in Figure~\ref{fig:sontag} represent the regions in
the space $(x,u)$ where $x X(x,u) < 0$. From the plot, it can be seen
that there exists no continuous function $x \mapsto u(x)$ defined on
$\real$ whose graph belongs to the union of the shaded areas.  Even
control systems whose inputs appear linearly (see
Example~\ref{ex:Brockett}) may be subject to obstructions that
preclude the existence of continuous state-dependent stabilizing
feedback~\cite{RWB:83a,EDS:98,AB-LR:05}.
\hfill $\blacksquare$

Numerous fundamental questions arise when dealing with discontinuous
dynamical systems.  The most basic question is the notion of a
solution.  For a discontinuous vector field, the existence of a
continuously differentiable solution, that is, a continuously
differentiable curve whose derivative follows the direction of the
vector field, is not guaranteed.  The following examples illustrate
the difficulties that arise in defining solutions of discontinuous
dynamical systems.

\refstepcounter{example}%
\subsection*{Example~\theexample: Brick on a frictional
  ramp}\label{ex:brick}

Consider a brick sliding on a ramp~\cite{DES:00}. As the brick moves,
it experiences a friction force in the opposite direction (see
Figure~\ref{fig:brick}(a)).  During sliding, the Coulomb friction
model states that the magnitude of the friction force is independent
of the magnitude of the velocity and is equal to the normal contact
force times the coefficient of friction.  The application of this
model to the sliding brick yields
\begin{align}\label{eq:brick}
  \dot v(t) = g (\sin \theta) - \nu g ( \cos \theta ) \sign (v(t)) ,
\end{align}
where $v$ is the velocity of the brick, $g$ is the acceleration due to
gravity, $\theta >0$ is the inclination of the ramp, $\nu$ is the
coefficient of friction, and $\sign (0) = 0$.  The right-hand side
of~\eqref{eq:brick} is a discontinuous function of $v$ because of the
presence of the sign function.  Figure~\ref{fig:brick}(b) shows the
phase plot of this system for several values of $\nu$.

Depending on the magnitude of the friction force, experiments show
that the brick stops and stays stopped. In other words, the brick
attains $v=0$ in finite time, and maintains $v=0$.  However, there is
no continuously differentiable solution of~\eqref{eq:brick} that
exhibits this type of behavior. To see this, note that $v=0$ and $\dot
v=0$ in~\eqref{eq:brick} imply $\sin \theta = 0$, which contradicts
$\theta>0$.  In order to explain this type of physical evolution, we
need to understand the discontinuity in~\eqref{eq:brick} and expand
our notion of solution beyond continuously differentiable solutions.
\hfill $\blacksquare$

\refstepcounter{example}%
\subsection*{Example~\theexample: Nonsmooth harmonic
  oscillator}\label{ex:oscillator}

Consider a unit mass subject to a discontinuous spring force.  The
spring does not exert any force when the mass is at the reference
position $x=0$.  When the mass is displaced to the right, that is,
$x>0$, the spring exerts a constant negative force that pulls it back
to the reference position. When the mass is displaced to the left,
that is, $x<0$, the spring exerts a constant positive force that pulls
it back to the reference position.  According to Newton's second law,
the system evolution is described by~\cite{FC:99}
\begin{align}\label{eq:nonsmooth-harmonic-oscillator}
  \ddot x + \sign(x) = 0 .
\end{align}
By defining the state variables $x_1 = x$ and $x_2 = \dot
x$,~\eqref{eq:nonsmooth-harmonic-oscillator} can be rewritten as
\begin{align}
  \dot x_1 (t) & =  x_2 (t) ,
  \label{eq:nonsmooth-harmonic-oscillator-a}
  \\
  \dot x_2 (t) & = -\sign (x_1 (t)) ,
  \label{eq:nonsmooth-harmonic-oscillator-b}
\end{align}
which is a nonsmooth version of the classical harmonic oscillator.
The phase portrait of this system is plotted in
Figure~\ref{fig:oscillator}(a).

It can be seen that $(0,0)$ is the unique equilibrium
of~\eqref{eq:nonsmooth-harmonic-oscillator-a}-\eqref{eq:nonsmooth-harmonic-oscillator-b}.
By discretizing the equations of motion, we find that the trajectories
of the system approach, as the time step is made smaller and smaller,
the set of curves in Figure~\ref{fig:oscillator}(b), which are the
level sets of the function $(x_1,x_2) \mapsto |x_1| +
\frac{x_2^2}{2}$.  These level sets are analogous to the level sets of
the function $(x_1,x_2) \mapsto x_1^2 + x_2^2$, which are the
trajectories of the classical harmonic oscillator $\ddot x + x = 0
$~\cite{HG:80}.  However, the trajectories in
Figure~\ref{fig:oscillator}(b) are not continuously differentiable
along the $x_2$-coordinate axis since the limiting velocity vectors
from the right of the axis and from the left of the axis do not
coincide.
\hfill $\blacksquare$

\refstepcounter{example}%
\subsection*{Example~\theexample: Move-away-from-nearest-neighbor
  interaction law}\label{ex:away}

\begin{figure}[htbp]
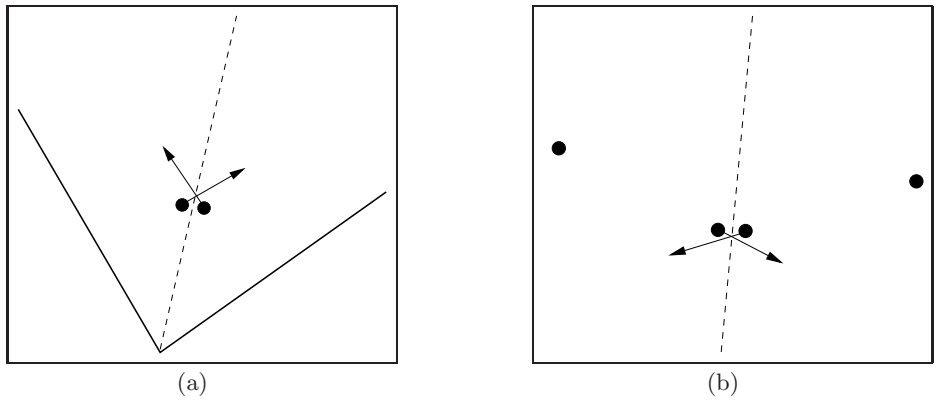

  \centering \subfigure[]{
    \fbox{\resizebox{.3\linewidth}{!}{\input{figure4a.tex}}}}
  \qquad\qquad \subfigure[]{ \fbox{
      \resizebox{.3\linewidth}{!}{\input{figure4b.tex}}}}
  \caption{Move-away-from-nearest-neighbor interaction law.  A node in
    (a) computes different directions of motion if placed slightly to
    the right or to the left of the bisector line defined by two
    polygonal boundaries.  A node in (b) computes different directions
    of motion if placed slightly to the right or to the left of the
    half plane defined by two other nodes. In both cases, the vector
    field is discontinuous.}
  \label{fig:away}
\end{figure}

Consider $n$ nodes $p_1,\dots,p_n$ evolving in a square $Q$ according
to the interaction rule ``move diametrically away from the nearest
neighbor.''  Note that this rule is not defined when two nodes are
located at the same point, that is, when the configuration is an
element of $\SS \triangleq \setdef{(p_1,\dots,p_n) \in Q^n}{p_i =p_j
  \; \text{for some} \; i \neq j}$. We thus consider only
configurations that belong to $Q^n \! \setminus \! \SS$.  Next, we
define the map that assigns to each node its nearest neighbor, where
the nearest neighbor may be an element of the boundary $\boundary{Q}$
of the square.  Note that the nearest neighbor of a node might not be
unique, that is, more than one node can be located at the same
(nearest) distance.  Hence, we define the \emph{nearest-neighbor map}
$\map{\nearestneighbor=(\nearestneighbor_1,\dots,\nearestneighbor_n)}{Q^n
  \! \setminus \! \SS}{Q^n}$ by arbitrarily selecting, for each $i \in
\until{n}$, an element,
\begin{align*}
  \nearestneighbor_i(p_1,\dots,p_n) \in \argmin
  \setdef{\TwoNorm{p_i-q}}{q \in \boundary{Q} \cup \{ p_1,\dots,p_n \}
    \! \setminus \! \{p_i\}} ,
\end{align*}
where $\argmin$ stands for the minimizing value of $q$, and
$\TwoNorm{\cdot}$ denotes the Euclidean norm. By definition,
$\nearestneighbor_i(p_1,\dots,p_n) \neq p_i$.
\sindex{NearestNeighborMap}{{$\nearestneighbor$}}{Nearest-neighbor map}%
For $i \in \until{n}$, we thus define the
move-away-from-nearest-neighbor interaction law
\begin{align}\label{eq:away}
  \dot p_i = \frac{p_i(t) -
    \nearestneighbor_i(p_1(t),\dots,p_n(t))}{\TwoNorm{p_i(t) -
      \nearestneighbor_i(p_1(t),\dots,p_n(t))}} .
\end{align}
Changes in the value of the nearest-neighbor map $\nearestneighbor$
induce discontinuities in the dynamical system.  For instance,
consider a node sufficiently close to a vertex of $\boundary{Q}$ so
that the closest neighbor to the node is an element of the boundary.
Depending on how the node is positioned with respect to the bisector
line passing through the vertex, the node computes different
directions of motion, see Figure~\ref{fig:away} for an illustration.

To analyze the dynamical system~\eqref{eq:away}, we need to understand
how the discontinuities affect its evolution.  Since each node moves
away from its nearest neighbors, it is reasonable to expect that the
nodes never run into each other, however, rigorous verification of
this property requires a proof.  We would also like to characterize
the asymptotic behavior of the trajectories of the
system~\eqref{eq:away}.  In order to study these questions, we need to
extend our notion of solution.  \hfill $\blacksquare$

\subsection*{Beyond continuously differentiable solutions}

Examples~\ref{ex:brick}-\ref{ex:away} above can be described by a
dynamical system of the form
\begin{align}\label{eq:mother-equation}
  \dot{x}(t) = X(x(t)) , \quad x(t_0) = x_0 ,
\end{align}
where $x \in \real^d$, $d$ is a positive integer, and
$\map{X}{\real^d}{\real^d}$ is not necessarily continuous.  We refer
to a continuously differentiable solution $t \mapsto x(t) $
of~\eqref{eq:mother-equation} as \emph{classical}.  Clearly, if $X$ is
continuous, then every solution is classical.  Without loss of
generality, we take $t_0 = 0$.  We consider only solutions that run
forward in time.

Examples~\ref{ex:brick}-\ref{ex:away} illustrate the limitations of
classical solutions, and confront us with the need to identify a
suitable notion for solutions of~\eqref{eq:mother-equation}.
Unfortunately, there is not a unique answer to this question.
Depending on the problem and objective at hand, different notions are
appropriate.  In this article, we restrict our attention to solutions
that are absolutely continuous.  Although not treated in detail in
this article, it is also possible to consider solutions that admit
discontinuities, and hence are not absolutely continuous. These
solutions are discussed in  ``Solutions with Jumps.''

The function $\map{\gamma}{[a,b]}{\real}$ is \emph{absolutely
  continuous}\index{function!absolutely continuous} if, for all $\eps
\in \realpositive$, there exists $\delta \in \realpositive$ such that,
for each finite collection $ \{(a_1,b_1), \ldots, (a_n,b_n)\}$ of disjoint
open intervals contained in $[a,b]$ with $ \sum_{i=1}^n (b_i-a_i)<
\delta$, it follows that
\begin{align*}
  \sum_{i=1}^n | \gamma(b_i)-\gamma(a_i) |< \eps.
\end{align*}
Equivalently~\cite{EDS:98}, $\gamma$ is absolutely continuous if there
exists a Lebesgue integrable function $\map{\kappa}{[a,b]}{\real}$
such that
\begin{align*}
  \gamma (t) = \gamma(a) + \int_a^t \kappa (s) ds, \quad t \in
  [a,b].
\end{align*}
Every absolutely continuous function is continuous. However, the
converse is not true, since the function $\map{\gamma}{[-1,1]}{\real}$
defined by $\gamma (t) = t \sin \big( \frac{1}{t} \big)$ for $t \neq
0$ and $\gamma (0) = 0$ is continuous, but not absolutely continuous.
Moreover, every continuously differentiable function is absolutely
continuous, but the converse is not true. For instance, the function
$\map{\gamma}{[-1,1]}{\real}$ defined by $\gamma (t) = |t|$ is
absolutely continuous but not continuously differentiable at $0$.  As
this example suggests, every absolutely continuous function is
differentiable almost everywhere.  Finally, every locally Lipschitz
function (see ``Locally Lipschitz Functions'') is absolutely
continuous, but the converse is not true.  For instance, the function
$\map{\gamma}{[0,1]}{\real}$ defined by $\gamma (t) = \sqrt{t}$ is
absolutely continuous but not locally Lipschitz at $0$.

Caratheodory solutions~\cite{AFF:88} are a generalization of classical
solutions.  Roughly speaking, Caratheodory solutions are absolutely
continuous curves that satisfy the integral version of the
differential equation~\eqref{eq:mother-equation}, that is,
\begin{align}\label{eq:integrable-mother-equation}
  x(t) = x(t_0) + \int_{t_0}^t X(x(s)) ds , \quad t>t_0 ,
\end{align}
where the integral is the Lebesgue integral.  By using the integral
form~\eqref{eq:integrable-mother-equation}, Caratheodory solutions
relax the classical requirement that the solution must follow the
direction of the vector field at all times, that is, the differential
equation~\eqref{eq:mother-equation} need not be satisfied on a set of
measure zero.  As shown in this article, Caratheodory solutions exist
for Example~\ref{ex:oscillator} but do not exist for
examples~\ref{ex:brick} and~\ref{ex:away}.

Alternatively, Filippov solutions~\cite{AFF:88} replace the
differential equation~\eqref{eq:mother-equation} by a differential
inclusion of the form
\begin{align}\label{eq:diff-inclusion-intro}
  \dot x(t) \in \setvaluedmap (x(t)),
\end{align}
where $\map{\setvaluedmap}{\real^d}{\parts{\real^d}}$ and
$\parts{\real^d}$ denotes the collection of all subsets of $\real^d$.
Filippov solutions are absolutely continuous curves. At a given
state~$x$, instead of focusing on the value of the vector field at
$x$, the idea behind Filippov solutions is to introduce a set of
directions that are determined by the values of the vector field $X$
in a neighborhood of $x$.  Differential
inclusions~\cite{GVS:01,JPA-AC:94} thus involve \emph{set-valued}
maps.  Just as a standard map or function takes a point in its domain
to \emph{a single} point in another space, a set-valued map takes a
point in its domain to \emph{a set of points} in another space.  A
differential inclusion thus specifies that the state derivative
belongs to a set of directions, rather than being a specific
direction.  This flexibility is crucial for providing conditions on
the discontinuous vector field under which Filippov solutions exist.
This solution notion plays a key role in many of the applications
mentioned above, including sliding mode control and mechanics with
Coulomb-like friction.

Unfortunately, the obstructions to continuous stabilization
illustrated in Example~\ref{ex:sontag} also
hold~\cite{EPR:94bis,JMC-LR:94} for Filippov solutions.
Sample-and-hold solutions~\cite{NNK-AIS:88}, which are also absolutely
continuous curves, turn out to be the appropriate notion for
circumventing these
obstructions~\cite{FHC-YSL-RJS-PRW:98,EDS:99,FHC:04}.  ``Additional
Solution Notions for Discontinuous Systems'' describes additional
solution notions for discontinuous systems.  In this article, we focus
on Caratheodory, Filippov, and sample-and-hold solutions.

\subsection*{Existence, uniqueness, and stability of solutions}

In addition to the notion of solution, we consider existence and
uniqueness of solutions as well as stability.  For ordinary
differential equations, it is well known that continuity of the vector
field does not guarantee uniqueness. Not surprisingly, no matter what
notion of solution is chosen for a discontinuous vector field,
nonuniqueness can occur.  We thus provide sufficient conditions for
uniqueness. We also present results specifically tailored to piecewise
continuous vector fields and differential inclusions.

The lack of uniqueness of solutions must be considered when we try to
establish properties such as local stability.  This issue is reflected
in the use of the adjectives weak and strong.  The word ``weak'' is
used when a property is satisfied by at least one solution starting
from each initial condition. On the other hand, ``strong'' is used
when a property is satisfied by all solutions starting from each
initial condition.  Therefore, for example, ``weakly stable
equilibrium point'' means that at least one solution starting close to
the equilibrium point remains close to it, whereas ``strongly stable
equilibrium point'' means that all solutions starting close to the
equilibrium point remain close to it.  For detailed definitions,
see~\cite{AFF:88,AB-LR:05}.

We present weak and strong stability results for discontinuous
dynamical systems and differential inclusions.  As suggested by
Example~\ref{ex:oscillator}, smooth Lyapunov functions do not suffice
to analyze the stability of discontinuous systems.  This fact leads
naturally to the subject of nonsmooth analysis.  In particular, we pay
special attention to the generalized gradient of a locally Lipschitz
function~\cite{FHC:83} and the proximal subdifferential of a lower
semicontinuous function~\cite{FHC-YSL-RJS-PRW:98}.  Building on these
notions, weak and strong monotonicity properties of candidate Lyapunov
functions can be established along the solutions of discontinuous
dynamical systems.  These results are used to provide generalizations
of Lyapunov stability theorems and the invariance principle, which
help us study the stability of solutions. To illustrate the
applicability of these results, we discuss in detail a class of
nonsmooth gradient flows.

There are two ways to apply the stability results presented here to
control systems. The first way is to choose a specific input function
and consider the resulting dynamical system.  The second way is to
associate with the control system the set-valued map that assigns each
state to the set of all vectors generated by the allowable inputs, and
consider the resulting differential inclusion.  Rather than focusing
on a particular input, this viewpoint allows us to consider the full
range of trajectories of the control system.  To analyze the stability
of the control system under this approach, we can use the nonsmooth
tools developed for differential inclusions.  We explore this idea in
detail.

Given the large body of work on discontinuous systems, our aim is to
provide a clear exposition of a few useful and central results.
``Additional Topics on Discontinuous Systems and Differential
Inclusions'' briefly discusses issues that are not considered in the
main exposition.

\subsection*{Organization of this article}

We start by reviewing basic results on the existence and uniqueness of
classical, that is, continuously differentiable, solutions of ordinary
differential equations. We also present several examples in which the
vector field fails to satisfy the standard smoothness properties.  We
then introduce various notions of solution for discontinuous systems,
discuss existence and uniqueness results, and present useful tools for
analysis. In preparation for the statement of stability results, we
introduce the generalized gradient and proximal sub\-differential from
nonsmooth analysis, and present various tools for their explicit
computation. Then, we develop analysis results to characterize the
stability and asymptotic convergence properties of the solutions of
discontinuous dynamical systems.  We illustrate these nonsmooth
stability results by means of examples, paying special attention to
gradient systems.  Throughout the discussion, we interchangeably use
``differential equation,'' ``dynamical system,'' and ``vector field.''
For reference, ``Index of Symbols'' summarizes the notation used
throughout.

To simplify the presentation, we have chosen to restrict our attention
to time-invariant vector fields, although most of the development can
be adapted to the time-varying setting. We briefly discuss
time-varying systems in ``Caratheodory Conditions for Time-varying
Vector Fields'' and ``Caratheodory Solutions of Differential
Inclusions.''  Likewise, for simplicity, we mostly consider vector
fields defined over the whole Euclidean space, although the exposition
can be carried out in more general settings such as open and connected
subsets of the Euclidean space.

\section*{Existence and Uniqueness for Ordinary Differential
  Equations}

In this section, we review basic results on existence and uniqueness
of classical solutions for ordinary differential equations. We also
present examples that do not satisfy the hypotheses of these results
but nevertheless exhibit existence and uniqueness of classical
solutions, as well as examples that do not possess such desirable
properties.

\subsection*{Existence of classical
  solutions}

Consider the differential equation
\begin{align}
  \label{eq:ODE-auto}
  \dot x (t) = X(x(t)) ,
\end{align}
where $\map{X}{\real^d}{\real^d}$ is a vector field. The point $x_e
\in \real^d$ is an \emph{equilibrium} of~\eqref{eq:ODE-auto} if $0 =
X(x_e)$.  A \emph{classical solution
  of~\eqref{eq:ODE-auto}}\index{solution!classical} on $[0,t_1]$ is a
continuously differentiable map $\map{x}{[0,t_1]}{\real^d}$ that
satisfies~\eqref{eq:ODE-auto}.  Note that, without loss of generality,
we consider only solutions that start at time $0$.  Usually, we refer
to $t \mapsto x(t)$ as a classical solution with initial condition $x
(0)$.  We sometimes write the initial condition as $x_0$ instead of
$x(0)$.  The solution $t \mapsto x(t)$ is
\emph{maximal}\index{solution!maximal} if it cannot be extended
forward in time, that is, if $t \mapsto x(t)$ is not the result of the
truncation of another solution with a larger interval of definition.
Note that the interval of definition of a maximal solution is either
of the form $ [0,T)$, where $T>0$, or $[0,\infty)$.

Continuity of the vector field suffices to guarantee the existence of
classical solutions, as stated by Peano's theorem~\cite{EAC-NL:55}.

\begin{proposition}\label{prop:ODE-existence}
  Let $\map{X}{\real^d}{\real^d}$ be continuous.  Then, for all $x_0
  \in \real^d$, there exists a classical solution
  of~\eqref{eq:ODE-auto} with initial condition $x(0) = x_0$.
\end{proposition}

The following example shows that, if the vector field is
discontinuous, then classical solutions of~\eqref{eq:ODE-auto} might
not exist.

\refstepcounter{example}%
\subsubsection*{Example~\theexample: Discontinuous vector field with
  nonexistence of classical
  solutions}\label{ex:discontinuous-non-existence}

Consider the vector field $\map{X}{\real}{\real}$ defined by
\begin{align}\label{eq:X1}
  X(x) =
  \begin{cases}
    - 1, & x>0 ,\\
    1, & x\le0 ,
  \end{cases}
\end{align}
which is discontinuous at $0$ (see
Figure~\ref{fig:counter-examples}(a)). 
Suppose that there exists a continuously differentiable function
$\map{x}{[0,t_1]}{\real}$ such that $\dot x (t) = X (x (t))$ and $x
(0) = 0$.  Then $\dot x (0) = X( x (0)) = X(0) = 1$, which implies
that, for all positive $t$ sufficiently small, $x(t) >0$ and hence
$\dot x (t) = X( x (t)) = -1$, which contradicts the fact that $t
\mapsto \dot x(t) $ is continuous.  Hence, no classical solution
starting from $0$ exists.
\hfill $\blacksquare$

In contrast to Example~\ref{ex:discontinuous-non-existence}, the
following example shows that the lack of continuity of the vector
field does not preclude the existence of classical solutions.

\refstepcounter{example}%
\subsubsection*{Example~\theexample: Discontinuous vector field with
  existence of classical
  solutions}\label{ex:discontinuous-existence}

Consider the vector field $\map{X}{\real}{\real}$,
\begin{align}\label{eq:X2}
  X(x) = - \sign (x) =
  \begin{cases}
    - 1, & x>0 ,\\
    0, & x = 0,\\
    1, & x<0 ,
  \end{cases}
\end{align}
which is discontinuous at $0$ (see
Figure~\ref{fig:counter-examples}(b)).  If $x(0)>0$, then the maximal
solution $\map{x}{[0, x(0))}{\real}$ is $x(t) = x(0) - t$, whereas, if
$x(0)<0$, then the maximal solution $\map{x}{[0, -x(0))}{\real}$ is $
x(t) = x(0) + t$.  Finally, if $x(0)=0$, then the maximal solution
$\map{x}{\realnonnegative}{\real}$ is $x(t) = 0$.  Hence the
associated dynamical system $\dot x(t) = X(x(t))$ has a classical
solution starting from every initial condition.  Although the vector
fields~\eqref{eq:X1} and~\eqref{eq:X2} are identical except for the
value at $0$, the existence of classical solutions is surprisingly
different.
\hfill $\blacksquare$

\subsection*{Uniqueness of classical
  solutions}

\emph{Uniqueness of classical solutions}\index{solution!uniqueness of}
of the differential equation~\eqref{eq:ODE-auto} means that every pair
of solutions with the same initial condition coincide on the
intersection of their intervals of existence.  In other words, if
$\map{x_1}{[0,t_1]}{\real^d}$ and $\map{x_2}{[0,t_2]}{\real^d}$ are
classical solutions of~\eqref{eq:ODE-auto} with $x_1 (0) = x_2(0)$,
then $x_1 (t) = x_2 (t) $ for all $t \in [0,t_1] \cap [0,t_2] =
[0,\min\{ t_1,t_2 \}]$.  Equivalently, we say that there exists a
unique maximal solution starting from each initial condition.

Uniqueness of classical solutions is guaranteed under a wide variety
of conditions. The book~\cite{RPA-VL:93}, for instance, is devoted to
collecting various uniqueness criteria.  Here, we focus on a
uniqueness criterion based on one-sided Lipschitzness.  The vector
field $\map{X}{\real^d}{\real^d}$ is \emph{one-sided Lipschitz} on $U
\subset \real^d$ if there exists $L >0$ such that, for all $y,y' \in
U$,
\begin{align}\label{eq:unique-ODE-sol}
  [X(y) - X(y')]^T (y-y') \le L \, \TwoNorm{y-y'}^2 .
\end{align}
This property is one-sided because it imposes a requirement on $X$
only when the angle between the two vectors on the left-hand side
of~\eqref{eq:unique-ODE-sol} is between $0$ and $180$ degrees.  The
following result, given in~\cite{RPA-VL:93}, uses this property to
provide a sufficient condition for uniqueness.

\begin{proposition}\label{prop:ODE-uniqueness}
  Let $\map{X}{\real^d}{\real^d}$ be continuous.  Assume that, for all
  $x \in \real^d$, there exists $\eps >0$ such that $X$ is one-sided
  Lipschitz on $\oball{\eps}{x}$.  Then, for all $x_0 \in \real^d$,
  there exists a unique classical solution of~\eqref{eq:ODE-auto} with
  initial condition $x(0) = x_0$.
\end{proposition}

Every vector field that is locally Lipschitz at $x$ (see ``Locally
Lipschitz Functions'') satisfies the one-sided Lipschitz condition on
a neighborhood of $x$, but the converse is not true. For example, the
vector field $\map{X}{\real}{\real}$ defined by $X(0) = 0$ and $X(x) =
x \log(|x|)$ for $x \neq 0$ is one-sided Lipschitz on a neighborhood
of $0$, but is not locally Lipschitz at $0$.  Furthermore, a one-sided
Lipschitz vector field can be discontinuous.  For example, the
discontinuous vector field~\eqref{eq:X2} is one-sided Lipschitz on a
neighborhood of~$0$.

The following result is an immediate consequence of
Proposition~\ref{prop:ODE-uniqueness}.
\begin{corollary}\label{cor:locally-Lipschitz-uniqueness}
  Let $\map{X}{\real^d}{\real^d}$ be locally Lipschitz. Then, for all
  $x_0 \in \real^d$, there exists a unique classical solution
  of~\eqref{eq:ODE-auto} with initial condition $x(0) = x_0$.
\end{corollary}

Although local Lipschitzness is typically invoked as in
Corollary~\ref{cor:locally-Lipschitz-uniqueness} to guarantee
uniqueness, Proposition~\ref{prop:ODE-uniqueness} shows that
uniqueness is guaranteed under weaker hypotheses.  The following
example shows that, if the hypotheses of
Proposition~\ref{prop:ODE-uniqueness} are not satisfied, then
solutions might not be unique.

\refstepcounter{example}%
\subsubsection*{Example~\theexample: Continuous, not one-sided
  Lipschitz vector field with nonunique classical
  solutions}\label{ex:continuous-non-uniqueness}

Consider the vector field $\map{X}{\real}{\real}$ defined by
\begin{align}\label{eq:X3}
  X(x) = \sqrt{|x|} .
\end{align}
This vector field is continuous everywhere, locally Lipschitz on
$\real \! \setminus \!  \{0 \}$ (see Figure~\ref{fig:counter-examples}(c)),
but is neither locally Lipschitz at $0$ nor one-sided Lipschitz on any
neighborhood of $0$.  The associated dynamical system $\dot x (t) =
X(x(t))$ has infinitely many maximal solutions starting from $0$,
namely, for all $a>0$, $\map{x_a}{\realnonnegative}{\real}$, where
\begin{align*}
  x_a(t) =
  \begin{cases}
    0, & 0 \le t \le a ,\\
    {(t-a)^2}/{4}, & \phantom{0 \le} t \ge a ,
  \end{cases}
\end{align*}
and $\map{x}{\realnonnegative}{\real}$, where $x(t) = 0$.
\hfill
$\blacksquare$

However, the following example shows that a differential equation
can possess a unique classical solution even when the hypotheses of
Proposition~\ref{prop:ODE-uniqueness} are not satisfied.

\refstepcounter{example}%
\subsubsection*{Example~\theexample: Continuous, not one-sided
  Lipschitz vector field with unique classical
  solutions}\label{ex:continuous-uniqueness}

Consider the vector field $\map{X}{\real}{\real}$ defined by
\begin{align}\label{eq:X4}
  X(x) = 
  \begin{cases}
    - x \log x, & x > 0,\\
    0, & x = 0,\\
    x \log (-x), & x < 0.
  \end{cases}
\end{align}
This vector field is continuous everywhere, and locally Lipschitz on
$\real \! \setminus \! \{0 \}$ (see Figure~\ref{fig:counter-examples}(d)).
However, $X$ is not locally Lipschitz at $0$ nor one-sided Lipschitz
on any neighborhood of $0$.  Nevertheless, the associated dynamical
system $\dot x (t) = X(x(t))$ has a unique solution starting from each
initial condition.  If $x(0)>0$, then the maximal solution
$\map{x}{\realnonnegative}{\real}$ is $x(t) = \exp (\log x(0)
\exp(-t))$, whereas, if $x(0)<0$, then the maximal solution
$\map{x}{\realnonnegative}{\real}$ is $x(t) = -\exp (\log (-x(0))
\exp(t))$.  Finally, if $x(0)=0$, then the maximal solution
$\map{x}{\realnonnegative}{\real}$ is $x(t) = 0$.
\hfill
$\blacksquare$

Note that Proposition~\ref{prop:ODE-uniqueness} assumes that the
vector field is continuous.  However, the discontinuous vector
field~\eqref{eq:X2} in Example~\ref{ex:discontinuous-existence} is
one-sided Lipschitz in a neighborhood of $0$ and, indeed, has a unique
classical solution starting from each initial condition.  This
observation suggests that discontinuous systems are not necessarily
more complicated or less ``well-behaved'' than continuous systems. The
continuous system~\eqref{eq:X3} in
Example~\ref{ex:continuous-non-uniqueness} does not have a unique
classical solution starting from each initial condition, whereas the
discontinuous system~\eqref{eq:X2} in
Example~\ref{ex:discontinuous-existence} does.  A natural question to
ask is under what conditions does a discontinuous vector field have a
unique solution starting from each initial condition. Of course, the
answer to this question relies on the notion of solution itself.  We
explore these questions in the next section.
  
\section*{Notions of Solution for Discontinuous Dynamical
  Systems}
  
The above discussion shows that the classical notion of solution is
too restrictive when considering a discontinuous vector field.  We
thus explore alternative notions of solution to reconcile this
mismatch.  To address the discontinuities of the differential
equation~\eqref{eq:ODE-auto}, we first relax the requirement that
solutions follow the direction specified by the vector field at all
times.  The precise mathematical notion corresponding to this idea is
that of Caratheodory solutions, which we introduce next.

\subsection*{Caratheodory solutions}

A \emph{Caratheodory solution
  of~\eqref{eq:ODE-auto}}\index{solution!Caratheodory} defined on
$[0,t_1] \subset \real$ is an absolutely continuous map
$\map{x}{[0,t_1]}{\real^d}$ that satisfies~\eqref{eq:ODE-auto} for
almost all $t \in [0,t_1]$ (in the sense of Lebesgue measure).  In
other words, a Caratheodory solution follows the direction specified
by the vector field except for a set of time instants that has measure
zero. Equivalently, Caratheodory solutions are absolutely continuous
functions that solve the integral version of~\eqref{eq:ODE-auto}, that
is,
\begin{align}\label{eq:dfn-Caratheodory}
  x (t) = x (0) + \int_{0}^t X (x (s)) ds.
\end{align}
Of course, every classical solution is also a Caratheodory solution.

\refstepcounter{example}%
\subsubsection*{Example~\theexample:  System with Caratheodory
  solutions and no classical solutions}\label{ex:Caratheodory}

The vector field $\map{X}{\real}{\real}$ defined by
\begin{align*}
  X(x) = 
  \begin{cases}
    1, & x>0 ,\\
    \frac{1}{2}, & x = 0, \\
    -1, & x<0 ,
  \end{cases}
\end{align*}
is discontinuous at $0$.  The associated dynamical system $\dot x(t) =
X(x(t))$ does not have a classical solution starting from $0$.
However, this system has two Caratheodory solutions starting from $0$,
namely, $\map{x_1}{\realnonnegative}{\real}$, where $x_1 (t) = t$, and
$\map{x_2}{\realnonnegative}{\real}$, where $x_2 (t) = -t$. Note that
both $x_1$ and $x_2$ violate the differential equation only at $t=0$,
that is, $\dot{x}_1 (0) \neq X (x_1(0))$ and $\dot{x}_2 (0) \neq X
(x_2(0))$.
\hfill $\blacksquare$

\subsubsection*{Example~\ref{ex:oscillator} revisited: Existence of
  Caratheodory solutions for the nonsmooth harmonic oscillator}
The nonsmooth harmonic oscillator in Example~\ref{ex:oscillator} does
not possess a classical solution starting from any initial condition
on the $x_2$-axis. However, the closed level sets depicted in
Figure~\ref{fig:oscillator}(b), when traversed clockwise, are
Caratheodory solutions.
\hfill $\blacksquare$

Unfortunately, the good news is quickly over since it is easy to find
examples of discontinuous dynamical systems that do not admit
Caratheodory solutions.  For example, the physical motions observed in
Example~\ref{ex:brick}, where the brick slides for a while and then
remains stopped, are not Caratheodory solutions.  Furthermore, the
discontinuous vector field~\eqref{eq:X1} does not admit a Caratheodory
solution starting from $0$.  Finally, the
move-away-from-nearest-neighbor interaction law in
Example~\ref{ex:away} is yet another example where Caratheodory
solutions do not exist, as we show next.

\subsubsection*{Example~\ref{ex:away} revisited: Nonexistence of
  Caratheodory solutions for the move-away-from-nearest-neighbor
  interaction law}

Consider one agent moving in the square $[-1,1]^2 \subset \real^2$
under the move-away-from-nearest-neighbor interaction law described in
Example~\ref{ex:away}.  Since no other agent is present in the square,
the agent moves away from the nearest polygonal boundary,
according to the vector field
\begin{align}\label{eq:no-Caratheodory-sol}
  X(x_1,x_2) = 
  \begin{cases}
    (-1,0), & -x_1 < x_2 \le x_1,\\
    (0,1), & x_2 < x_1 \le -x_2,\\
    (1,0), & x_1 \le x_2 < -x_1,\\
    (0,-1), & -x_2 \le x_1 < x_2.
  \end{cases}
\end{align}
Since the move-away-from-nearest-neighbor interaction law takes
multiple values on the diagonals $\setdef{(a, \pm a) \in [-1,1]^2}{a
  \in [-1,1]}$ of the square, we choose one of these values in the
definition~\eqref{eq:no-Caratheodory-sol} of~$X$.  The phase portrait
in Figure~\ref{fig:ex1}(a) shows that the vector field $X$ is
discontinuous on the diagonals.

The dynamical system $\dot x (t) = X (x(t))$ has no Caratheodory
solution if and only if the initial condition belongs to the
diagonals.  This fact can be justified as follows.  Consider the four
open regions of the square separated by the diagonals, see
Figure~\ref{fig:ex1}(a). On the one hand, if the initial condition
belongs to one of these regions, then the dynamical system has a
classical solution; depending on which region the initial condition
belongs to, the agent moves either vertically or horizontally toward
the diagonals.  On the other hand, on the diagonals of the square, $X$
pushes trajectories outward, whereas, outside the diagonals, $X$
pushes trajectories inward. Therefore, if the initial condition
belongs to the diagonals, then the only candidate trajectory for a
Caratheodory solution is a trajectory that moves along the diagonals.
From the definition~\eqref{eq:dfn-Caratheodory} of Caratheodory
solution, it is clear that a trajectory that moves along the diagonals
of the square is not a Caratheodory solution.  We show later that this
trajectory is instead a Filippov solution.
\hfill
$\blacksquare$

\bigskip
\noindent \textbf{\textit{Sufficient conditions for the existence of
    Caratheodory solutions}} \medskip

Conditions under which Caratheodory solutions exist are discussed in
``Caratheodory Conditions for Time-varying Vector Fields.'' For
time-invariant vector fields, the Caratheodory conditions given by
Proposition~\ref{prop:Caratheodory-existence} specialize to continuity
of the vector field. This requirement provides no improvement over
Proposition~\ref{prop:ODE-existence}, which guarantees the existence
of classical solutions under continuity.

Therefore, it is of interest to determine conditions for the existence
of Caratheodory solutions specifically tailored to time-invariant
vector fields. For example, directionally continuous vector fields are
considered in~\cite{AP:76}.  A vector field
$\map{X}{\real^d}{\real^d}$ is \emph{directionally continuous}
\index{vector field!directionally continuous}%
if there exists $\delta \in \realpositive$ such that, for every $x \in
\real^d$ with $X(x) \neq 0$ and every sequence $\{x_n\}_{n \in
  \natural} \subset \real^d$ with $x_n \to x$ and
\begin{align}\label{eq:directionally-continuous}
  \TwoNormBig{\frac{x_n-x}{\TwoNorm{x_n -x}} -
    \frac{X(x)}{\TwoNorm{X(x)}}} < \delta , \quad n \in \natural ,
\end{align}
it follows that $X(x_n) \to X(x)$.  If the vector field $X$ is
directionally continuous, then, for all $x_0 \in \real^d$, there
exists a Caratheodory solution of~\eqref{eq:ODE-auto} with initial
condition $x(0) = x_0$.

Patchy vector fields~\cite{FA-AB:99} are another class of
time-invariant, discontinuous vector fields that have Caratheodory
solutions.  Additional conditions for the existence and uniqueness of
Caratheodory conditions can be found
in~\cite{AB:98,AB-WS:98,JII-AJS:00}.  Beyond differential equations,
Caratheodory solutions can also be defined for differential
inclusions, as explained in ``Set-valued Maps'' and ``Caratheodory
Solutions of Differential Inclusions.''

\subsection*{Filippov solutions}

The above discussion shows that the relaxation of the value of the
vector field on a set of times of measure zero in the definition of
Caratheodory solution is not always sufficient to guarantee that such
solutions exist.  Due to the discontinuity of the vector field, its
value can exhibit significant variations arbitrarily close to a given
point, and this mismatch might make it impossible to construct a
Caratheodory solution.

What if, instead of focusing on the value of the vector field at
individual points, we consider how the vector field looks like
\emph{around} each point?  The idea of looking at a neighborhood of
each point is at the core of the notion of Filippov
solution~\cite{AFF:88}. Closely related notions are those of
Krasovskii solution~\cite{NNK:63} and Sentis solution~\cite{RS:78}.

The mathematical framework for formalizing this neighborhood idea uses
set-valued maps.  The idea is to associate a set-valued map to
$\map{X}{\real^d}{\real^d}$ by looking at the neighboring values of
$X$ around each point. Specifically, for $x \in \real^d$, the vector
field $X$ is evaluated at the points belonging to $\oball{\delta}{x}$,
which is the open ball centered at $x$ with radius $\delta>0$.  We
examine the effect of $\delta$ approaching $0$ by performing this
evaluation for smaller and smaller $\delta$. For additional
flexibility, we exclude an arbitrary set of measure zero in
$\oball{\delta}{x}$ when evaluating $X$, so that the outcome is the
same for two vector fields that differ on a set of measure zero.

Mathematically, the above procedure can be summarized as follows.  Let
$\parts{\real^d}$ denote the collection of subsets of $\real^d$.  For
$\map{X}{\real^d}{\real^d}$, define the \emph{Filippov set-valued
  map}\index{set-valued map!Filippov}
$\map{F[X]}{\real^d}{\parts{\real^d}}$ by
\begin{align}\label{eq:Filippov-set-valued-map}
  F[X](x) & \triangleq \bigcap_{\delta >0} \bigcap_{\mu (S) = 0}
  \overline{\co} \{ X(\oball{\delta}{x} \! \setminus \! S) \} , \quad
  x \in \real^d .
\end{align}
In~\eqref{eq:Filippov-set-valued-map}, $\overline{\co}$ denotes convex
closure, and $\mu$ denotes Lebesgue measure. Because of the way the
Filippov set-valued map is defined, the value of $F[X]$ at a point $x$
is independent of the value of the vector field $X$ at $x$.
\sindex{ConvexClosure}{{$\overline{\co}(S)$}}{Convex closure of a set
  $S \subseteq \real^d$}%
\sindex{FilippovSetValuedMap}{{$F[X]$}}{Filippov set-valued map
  associated with a vector field $X:\real^d \rightarrow \real^d$}%

\subsubsection*{Examples~\ref{ex:discontinuous-non-existence}
  and~\ref{ex:discontinuous-existence} revisited: Filippov set-valued
  map of the sign function}

Let us compute the Filippov set-valued map for the vector
fields~\eqref{eq:X1} and~\eqref{eq:X2}.  Since both vector fields
differ only at $0$, which is a set of measure zero, their associated
Filippov set-valued maps are identical and equal to
$\map{F[X]}{\real}{\parts{\real}}$, where
\begin{align}\label{eq:Filippov-set-sign}
  F[X] (x) = 
  \begin{cases}
    -1 , & x>0 ,\\
    [-1,1],  & x=0 ,\\
    1 , & x<0 .
  \end{cases}
\end{align}
Note that this Filippov set-valued map is multiple-valued only at the
point of discontinuity of the vector field. This observation is valid
for all vector fields.
\hfill
$\blacksquare$

We are now ready to handle the discontinuities of the vector field $X$
by using the Filippov set-valued map of $X$.  We replace the
differential equation $\dot x (t) = X(x(t))$ by the differential
inclusion
\begin{align}\label{eq:Filippov-inclusion}
  \dot x (t) \in F[X] (x(t)) .
\end{align}
A \emph{Filippov solution
  of~\eqref{eq:ODE-auto}}\index{solution!Filippov} on $[0,t_1] \subset
\real$ is an absolutely continuous map $\map{x}{[0,t_1]}{\real^d}$
that satisfies~\eqref{eq:Filippov-inclusion} for almost all $t \in
[0,t_1]$.  Equivalently, a Filippov solution of~\eqref{eq:ODE-auto} is
a Caratheodory solution of the differential
inclusion~\eqref{eq:Filippov-inclusion}, see ``Caratheodory Solutions
of Differential Inclusions.''

Because of the way the Filippov set-valued map is defined, a vector
field that differs from $X$ on a set of measure zero has the same
Filippov set-valued map, and hence the same set of solutions.  The
next result establishes mild conditions under which Filippov solutions
exist~\cite{AFF:88,JPA-AC:94}.

\begin{proposition}\label{prop:existence-Filippov} 
  Let $\map{X}{\real^d}{\real^d}$ be measurable and locally
  essentially bounded, that is, bounded on a bounded neighborhood of
  every point, excluding sets of measure zero.  Then, for all $x_0 \in
  \real^d$, there exists a Filippov solution of~\eqref{eq:ODE-auto}
  with initial condition $x(0) = x_0$.
\end{proposition}

In Proposition~\ref{prop:existence-Filippov}, the hypotheses on the
vector field imply that the associated Filippov set-valued map
satisfies all of the hypotheses of
Proposition~\ref{prop:existence-solution} (see ``Caratheodory
Solutions of Differential Inclusions''), which in turn guarantees the
existence of Filippov solutions.

\subsubsection*{Examples~\ref{ex:discontinuous-non-existence}
  and~\ref{ex:discontinuous-existence} revisited: Existence of
  Filippov solutions for the sign function}

The application of Proposition~\ref{prop:existence-Filippov} to the
bounded vector fields in examples~\ref{ex:discontinuous-non-existence}
and~\ref{ex:discontinuous-existence} guarantees that a Filippov
solution of~\eqref{eq:ODE-auto} exists for both examples starting from
each initial condition.  Furthermore, since the vector
fields~\eqref{eq:X1} and~\eqref{eq:X2} have the same Filippov
set-valued map~\eqref{eq:Filippov-set-sign},
examples~\ref{ex:discontinuous-non-existence}
and~\ref{ex:discontinuous-existence} have the same maximal Filippov
solutions. If $x(0)>0$, then the maximal solution
$\map{x}{\realnonnegative}{\real}$ is
\begin{align*}
  x(t) =
  \begin{cases}
    x(0) - t, & t \le x(0) , \\
    0, & t \ge x(0) ,
  \end{cases}
\end{align*}
whereas, if $x(0)<0$, then the maximal solution
$\map{x}{\realnonnegative}{\real}$ is
\begin{align*}
  x(t) =
  \begin{cases}
    x(0) + t, & t \le -x(0) , \\
    0, & t \ge -x(0) .
  \end{cases}
\end{align*}
Finally, if $x(0)=0$, then the maximal solution
$\map{x}{\realnonnegative}{\real}$ is $x(t) = 0$.
\hfill
$\blacksquare$

Similar computations can be made for the
move-away-from-nearest-neighbor interaction law in
Example~\ref{ex:away} to show that Filippov solutions exist starting
from each initial condition, as we show next.

\subsubsection*{Example~\ref{ex:away} revisited: Filippov solutions
  for the move-away-from-nearest-neighbor interaction
  law}

Consider again the discontinuous vector field for one agent moving in
the square $[-1,1]^2 \subset \real^2$ under the
move-away-from-nearest-neighbor interaction law described in
Example~\ref{ex:away}.  The corresponding set-valued map
$\map{F[X]}{[-1,1]^2}{\parts{\real^2}}$ is given by
\begin{align}\label{eq:Filippov-set-away}
  F[X](x_1,x_2) =
  \begin{cases}
    \setdef{(y_1,y_2) \in \real^2}{|y_1+y_2| \le 1, |y_1-y_2| \le
      1}, &  (x_1,x_2)=(0,0),\\
    \{ (-1,0) \}, & -x_1 < x_2 < x_1,\\
    \setdef{(y_1,y_2) \in \real^2}{y_1+y_2=-1,y_1 \in[-1,0]},
    &  0 < x_2 = x_1,\\
    \{ (0,1) \}, & x_2 < x_1 < -x_2,\\
    \setdef{(y_1,y_2) \in \real^2}{y_1-y_2=-1,y_1 \in[-1,0]},
    &  0 < -x_1 = x_2,\\
    \{ (1,0) \}, & x_1 < x_2 < -x_1,\\
    \setdef{(y_1,y_2) \in \real^2}{y_1+y_2=1,y_1 \in[0,1]},
    &  x_2 = x_1 < 0,\\
    \{ (0,-1) \}, & -x_2 < x_1 < x_2,\\
    \setdef{(y_1,y_2) \in \real^2}{y_1-y_2=1,y_1 \in[0,1]}, & 0 < x_1
    = - x_2.
  \end{cases}
\end{align}
Since $X$ is bounded, it follows from
Proposition~\ref{prop:existence-Filippov} that a Filippov solution
exists starting from each initial condition, see
Figure~\ref{fig:ex1}(b).  In particular, each solution starting from a
point on a diagonal is a straight line flowing along the diagonal
itself and reaching $(0,0)$.  For example, the maximal solution
$\map{x}{\realnonnegative}{\real^2}$ starting from $(a,a) \in \real^2$
is given by
\begin{align}\label{eq:sliding-away}
  t \mapsto x (t) =
  \begin{cases}
    (a- \frac{1}{2}  \sign(a) t, a - \frac{1}{2}  \sign(a) t), & t \le
    |a| ,\\ 
    (0,0), & t \ge |a| .
  \end{cases}
\end{align}
Note that the solution \emph{slides} along the diagonals, following a
convex combination of the limiting values of $X$ around the diagonals,
rather than the direction specified by $X$ itself.  We study this type
of behavior in more detail in the section ``Piecewise continuous
vector fields and sliding motions.''
\hfill
$\blacksquare$

\bigskip
\noindent \textbf{\textit{Relationship between Caratheodory and
    Filippov solutions}} \medskip

In general, Caratheodory and Filippov solutions are not related.  A
vector field for which both notions of solution exist but Filippov
solutions are not Caratheodory solutions is given in~\cite{FA-AB:99}.
The following is an example of the opposite case.

\refstepcounter{example}%
\subsubsection*{Example~\theexample: Vector field with a Caratheodory
  solution that is not a Filippov
  solution}\label{ex:Caratheodory-not-Filippov}

Consider the vector field $\map{X}{\real}{\real}$ given by
\begin{align}
  X (x) =
  \begin{cases}
    1 , & x\neq  0 ,\\
    0, & x=0 .
  \end{cases}
\end{align}
The dynamical system~\eqref{eq:ODE-auto} has two Caratheodory
solutions starting from $0$, namely,
$\map{x_1}{\realnonnegative}{\real}$, where $x_1(t) = 0$, and
$\map{x_2}{\realnonnegative}{\real}$, where $x_2(t) = t$. However, the
associated Filippov set-valued map $\map{F[X]}{\real}{\parts{\real}}$
is $F[X](x) = \{1\}$, and hence $t \mapsto x_2(t)$ is the unique
Filippov solution starting from $0$.
\hfill
$\blacksquare$

On a related note, Caratheodory solutions are always Krasovskii
solutions (see ``Additional Solution Notions for Discontinuous
Systems'') but the converse is not true~\cite{FC:99}.

\bigskip
\noindent \textbf{\textit{Computing the Filippov set-valued map}}
\medskip

Computing the Filippov set-valued map can be a daunting task.  A
calculus is developed in~\cite{BP-SSS:87} to simplify this
calculation. We summarize some useful facts below.  Note that the
Filippov set-valued map can also be constructed for maps of the form
$\map{X}{\real^d}{\real^m}$, where $d$ and $m$ are not necessarily
equal.

\textbf{Consistency.} If $\map{X}{\real^d}{\real^m}$ is continuous at
$x \in \real^d$, then
\begin{align}\label{eq:Filippov-consistency}
  F[X] (x) = \{X(x)\} .
\end{align}

\textbf{Sum rule.} If $\map{X_1,X_2}{\real^d}{\real^m}$ is locally
bounded at $x \in \real^d$, then
\begin{align}\label{eq:Filippov-sum}
  F[X_1+X_2] (x) \subseteq F[X_1] (x) + F[X_2] (x) .
\end{align}
Moreover, if either $X_1$ or $X_2$ is continuous at $x$, then equality
holds.

\textbf{Product rule.} If $\map{X_1}{\real^d}{\real^m}$ and
$\map{X_2}{\real^d}{\real^n}$ are locally bounded at $x \in \real^d$,
then
\begin{align}\label{eq:Filippov-product}
  F [(X_1,X_2)^T] (x) \subseteq F[X_1] (x) \times F[X_2] (x) .
\end{align}
Moreover, if either $X_1$ or $X_2$ is continuous at $x$, then equality
holds.

\textbf{Chain rule.} If $\map{Y}{\real^d}{\real^n}$ is continuously
differentiable at $x \in \real^d$ with Jacobian rank $n$, and
$\map{X}{\real^n}{\real^m}$ is locally bounded at $Y(x) \in \real^n$,
then
\begin{align}\label{eq:Filippov-chain}
  F[X \circ Y] (x) = F[X] (Y(x)) .
\end{align}

\textbf{Matrix transformation rule.} If $\map{X}{\real^d}{\real^{m}}$
is locally bounded at $x \in \real^d$ and $\map{Z}{\real^d}{\real^{d
    \times m}}$ is continuous at $x \in \real^d$, then
\begin{align}\label{eq:Filippov-matrix}
  F[Z \, X](x) = Z(x) \, F[X] (x) .
\end{align}

\bigskip
\noindent \textbf{\textit{Piecewise continuous vector fields and sliding
  motions}} 
\medskip

In this section we consider vector fields that are continuous
everywhere except on a surface of the state space.  Consider, for
instance, two continuous dynamical systems, one per side of the
surface, glued together to give rise to a discontinuous dynamical
system.  Here, we analyze the properties of the Filippov solutions of
this type of systems.

The vector field $\map{X}{\real^d}{\real^d}$ is \emph{piecewise
  continuous}
\index{vector field!piecewise continuous}%
if there exists a finite collection of disjoint, open, and connected
sets $\domain_1,\dots,\domain_m \subset \real^d$ whose closures cover
$\real^d$, that is, $\real^d = \union_{k=1}^m \overline{\domain_k}$,
such that, for all $k = 1, \dots, m$, the vector field $X$ is
continuous on $\domain_k$.  We further assume that the restriction of
$X$ to $\domain_k$ admits a continuous extension to the closure
$\ov{\domain_k}$, which we denote by $X_{|\ov{\domain_k}}$.  Every
point of discontinuity of $X$ must therefore belong to the union of
the boundaries of the sets $\domain_1,\dots,\domain_m$. Let us denote
by $S_X \subseteq \boundary{\domain_1} \union \dots \union
\boundary{\domain_m}$ the set of points where $X$ is discontinuous.
Note that $S_X$ has measure zero.

\sindex{DiscontinuitySet}{{$S_X$}}{Set of points where the vector
  field $\map{X}{\real^d}{\real^d}$ is discontinuous}%

The Filippov set-valued map of a piecewise continuous vector field $X$
is given by the expression
\begin{align}
  F[X] (x) = \overline{\co} \setdef{\lim_{i \rightarrow \infty}
    X(x_i)}{x_i \rightarrow x \, , \; x_i \not \in S_X} .
\end{align}
This set-valued map can be computed as follows. At points of
continuity of $X$, that is, for $x \not \in S_X$, the consistency
property~\eqref{eq:Filippov-consistency} implies that $F[X](x) = \{
X(x) \}$. At points of discontinuity of $X$, that is, for $x \in S_X$,
$F[X](x)$ is a convex polyhedron in $\real^d$ of the form
\begin{align}\label{eq:Filippov-convex-polyhedron}
  F[X](x) = \co \setdef{X_{\ov{\domain_k}} (x)}{x \in
    \boundary{\domain_k}} .
\end{align}
 
As an illustration, let us revisit
examples~\ref{ex:brick}-\ref{ex:away}.

\subsubsection*{Examples~\ref{ex:brick}-\ref{ex:away} revisited:
  Computation of the Filippov set-valued map and Filippov
  solutions}

The vector field of the sliding brick in Example~\ref{ex:brick} is
piecewise continuous with $\domain_1 = \setdef{v \in \real}{v<0}$ and
$\domain_2 = \setdef{v \in \real}{v>0}$.  Note that the restriction of
the vector field to $\domain_1$ can be continuously extended to
$\ov{\domain_1}$ by setting $ X_{|\ov{\domain_1}} (0) = g (\sin \theta
+ \nu \cos \theta)$.  Likewise, the restriction of the vector field to
$\domain_2$ can be continuously extended to $\ov{\domain_2}$ by
setting $X_{|\ov{\domain_2}} (0) = g (\sin \theta - \nu \cos \theta)$.
Therefore, the associated Filippov set-valued map
$\map{F[X]}{\real}{\parts{\real}}$ is given by
\begin{align*}
  F[X] (v) =
  \begin{cases}
    \{ g (\sin \theta - \nu \cos \theta) \}, & v>0, \\
    \setdef{g (\sin \theta - d\, \nu \cos \theta)}{d \in [-1,1]}, &
    v = 0,\\
    \{ g (\sin \theta + \nu \cos \theta) \}, & v<0,
  \end{cases}
\end{align*}
which is singleton-valued for all $v \not \in S_X = \{0\}$, and a
closed segment at $v=0$.  If the friction coefficient $\nu$ is large
enough, that is, satisfies $\nu > \tan \theta$, then $F[x](v)<0$ if
$v>0$ and $F[x](v)>0$ if $v<0$. Therefore, the Filippov solutions that
start with an initial positive velocity $v$ eventually reach $v=0$ and
remain at $0$. This fact precisely corresponds to the observed
physical motions for Example~\ref{ex:brick}, where the brick slides
for a while and then remains stopped. This example shows that Filippov
solutions have physical significance.

The vector field $\map{X}{\real^2}{\real^2}$ for the nonsmooth
harmonic oscillator in Example~\ref{ex:oscillator} is continuous on
each of the half planes $\{ \domain_{1}, \domain_{2}\}$, with
\begin{alignat*}{2}
  \domain_{1} &= \setdef{(x_1,x_2) \in \real^2}{x_1 < 0}
  , \\
  \domain_{2} &= \setdef{(x_1,x_2) \in \real^2}{x_1 > 0} ,
\end{alignat*}
and discontinuous on $S_X = \setdef{(0,x_2)}{x_2 \in \real}$.
Therefore, $X$ is piecewise continuous.  Its Filippov set-valued map
$\map{F[X]}{\real^2}{\parts{\real^2}}$ is given by
\begin{align*}
  F[X] (x_1,x_2) =
  \begin{cases}
    \{ (x_2,-\sign (x_1)) \}, & x_1 \neq 0  , \\
    \{ x_2 \} \times [-1,1], & x_1 = 0 .
  \end{cases}
\end{align*}
Therefore, the closed level sets depicted in
Figure~\ref{fig:oscillator}(b), when traversed clockwise, are Filippov
solutions.

The discontinuous vector field $\map{X}{[-1,1]^2}{\real^2}$ for one
agent moving in the square $[-1,1]^2 \subset \real^2$ under the
move-away-from-nearest-neighbor interaction law described in
Example~\ref{ex:away} is piecewise continuous, with
\begin{align*}
  \domain_1 & = \setdef{(x_1,x_2) \in [-1,1]^2}{-x_1 < x_2 < x_1} ,\\
  \domain_{2} &= \setdef{(x_1,x_2) \in [-1,1]^2}{x_2 < x_1
    < -x_2} ,\\
  \domain_3 & = \setdef{(x_1,x_2) \in [-1,1]^2}{x_1 < x_2 < -x_1} ,\\
  \domain_{4} &= \setdef{(x_1,x_2) \in [-1,1]^2}{-x_2 < x_1 <
    x_2} .
\end{align*}
Its Filippov set-valued map, described
in~\eqref{eq:Filippov-set-away}, maps points outside the diagonals
$S_X = \setdef{(a,\pm a) \in [-1,1]^2}{a \in [-1,1]}$ to singletons,
maps points in $S_X \! \setminus \! \{ (0,0)\}$ to closed segments, and maps
$(0,0)$ to a square polygon. Several Filippov solutions starting from
various initial conditions are plotted in Figure~\ref{fig:ex1}(b).
\hfill
$\blacksquare$

\begin{figure}[htbp]
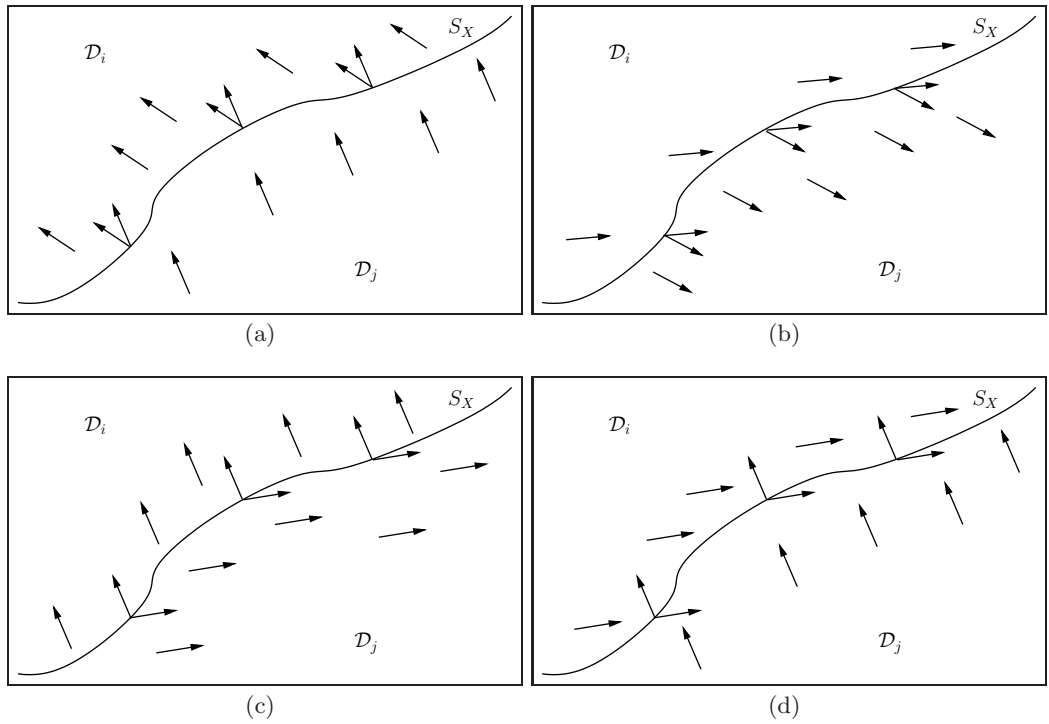

  \centering
  \subfigure[]{\fbox{\resizebox{0.4\linewidth}{!}{\input{figure8a.tex}}}}
  \subfigure[]{\fbox{\resizebox{0.4\linewidth}{!}{\input{figure8b.tex}}}}
  \\
  \subfigure[]{\fbox{\resizebox{0.4\linewidth}{!}{\input{figure8c.tex}}}}
  \subfigure[]{\fbox{\resizebox{0.4\linewidth}{!}{\input{figure8d.tex}}}}
  \caption{Piecewise continuous vector fields. These dynamical systems
    are continuous on $\domain_1$ and $\domain_2$, and discontinuous
    on $S_X$.  In (a) and (b), Filippov solutions, known as
    transversally crossing trajectories, cross $S_X$.  In (c), there
    are two Filippov solutions, known as repulsive trajectories,
    starting from each point in $S_X$.  Finally, in (d), the Filippov
    solutions that reach $S_X$, known as attractive trajectories,
    continue sliding along $S_X$.}\label{fig:sliding}
\end{figure}

Let us now discuss what happens on the points of discontinuity of the
vector field $\map{X}{\real^d}{\real^d}$.  Suppose that $x \in S_X$
belongs to just two boundary sets, that is, $x \in
\boundary{\domain_i} \intersection \boundary{\domain_j}$, for some
distinct $i,j \in \until{m}$, but $x \not \in \boundary{\domain_k}$,
for $k \in \until{m} \! \setminus \! \{i, j\}$. In this case,
\begin{align*}
  F[X] (x) = \co \{X_{|\ov{\domain_i}} (x) , X_{|\ov{\domain_j}} (x)
  \}.
\end{align*}
We consider three possibilities.  First, if all of the vectors
belonging to $F[X] (x)$ point into $\domain_i$, then a Filippov
solution that reaches $S_X$ at $x$ continues its motion in $\domain_i$
(see Figure~\ref{fig:sliding}(a)).  Likewise, if all of the vectors
belonging to $F[X] (x)$ point into $\domain_j$, then a Filippov
solution that reaches $S_X$ at $x$ continues its motion in $\domain_j$
(see Figure~\ref{fig:sliding}(b)).  Finally, if a vector belonging to
$F[X] (x)$ is tangent to $S_X$, then either all Filippov solutions
that start at $x$ leave $S_X$ immediately (see
Figure~\ref{fig:sliding}(c)), or there exist Filippov solutions that
reach the set $S_X$ at $x$, and \emph{remain} in $S_X$ afterward (see
Figure~\ref{fig:sliding}(d)).

The last kind of trajectories are called \emph{sliding motions}, since
they slide along the boundaries of the sets where the vector field is
continuous.  This type of behavior is illustrated for
Example~\ref{ex:away} in~\eqref{eq:sliding-away}.  Sliding motions can
also occur along points belonging to the intersection of more than two
sets in $\ov{\domain_1},\dots,\ov{\domain_m}$.  The theory of sliding
mode control builds on the existence of this type of trajectories to
design stabilizing feedback controllers. These controllers induce
sliding surfaces with the desired properties so that the closed-loop
system is stable~\cite{VIU:92,CE-SKS:98}.

The solutions of piecewise continuous vector fields appear frequently
in state-dependent switching dynamical
systems~\cite{DL:03,MDB-CJB-ARC-PK:07}.  Consider, for instance, the
case of two dynamical systems with the same unstable equilibrium.  By
identifying an appropriate switching surface between the two systems,
it is often possible to synthesize a discontinuous dynamical system
for which the equilibrium is stable.

\bigskip
\noindent \textbf{\textit{Uniqueness of Filippov
    solutions}}
\medskip

A discontinuous dynamical system does not necessarily have a unique
Filippov solution starting from each initial condition.  The situation
depicted in Figure~\ref{fig:sliding}(c) is a qualitative example in
which multiple Filippov solutions exist for the same initial
condition.  The following is another example of lack of uniqueness.

\refstepcounter{example}%
\subsubsection*{Example~\theexample: Vector field with nonunique
  Filippov solutions}\label{ex:nonuniqueness-Filippov}

Consider the vector field $\map{X}{\real}{\real}$ defined by $X(x) =
\sign (x)$.  For all $x_0 \in \real \! \setminus \! \{0\}$, the
system~\eqref{eq:ODE-auto} has a unique Filippov solution starting
from $x_0$. However, the system~\eqref{eq:ODE-auto} has three maximal
solutions $\map{x_1,x_2,x_3}{\realnonnegative}{\real}$ starting from
$x_0=0$ given by $x_1(t) = -t$, $ x_2(t) = 0$, and $ x_3(t) = t$.
\hfill $\blacksquare$

We now provide two complementary uniqueness results for Filippov
solutions.  The first result~\cite{AFF:88} considers the Filippov
set-valued map associated with a discontinuous vector field, and
identifies conditions under which
Proposition~\ref{prop:unique-sol-inclusion} in ``Caratheodory
Solutions of Differential Inclusions'' can be applied to the resulting
differential inclusion. In order to state this result, we need to
introduce the following definition.  The vector field
$\map{X}{\real^d}{\real^d}$ is \emph{essentially one-sided Lipschitz}
on $U \subset \real^d$ if there exists $L >0$ such that, for almost
all $y,y' \in U$,
\begin{align}\label{eq:unique-Filipov-sol-I}
  [X(y) - X(y')]^T (y-y') \le L \, \TwoNorm{y-y'}^2 .
\end{align}
The first uniqueness result for Filippov solutions is stated next.

\begin{proposition}\label{prop:uniqueness-Filippov-I}
  Let $\map{X}{\real^d}{\real^d}$ be measurable and locally
  essentially bounded.  Assume that, for all $x \in \real^d$, there
  exists $\eps >0$ such that $X$ is essentially one-sided Lipschitz on
  $\oball{\eps}{x}$.  Then, for all $x_0 \in \real^d$, there exists a
  unique Filippov solution of~\eqref{eq:ODE-auto} with initial
  condition $x(0) = x_0$.
\end{proposition}

Note the parallelism of this result with
Proposition~\ref{prop:ODE-uniqueness} for ordinary differential
equations with a continuous vector field $X$.  Let us apply
Proposition~\ref{prop:uniqueness-Filippov-I} to an example.

\refstepcounter{example}%
\subsubsection*{Example~\theexample: Vector field with unique Filippov
  solutions}\label{ex:uniqueness-Filippov}

Let $\rational$ denote the set of rational numbers, and define the
vector field $\map{X}{\real}{\real}$ by
\begin{align*}
  X(x) =
  \begin{cases}
    1, & x \in \rational ,\\
    -1, & x \not \in \rational ,
  \end{cases}
\end{align*}
which is discontinuous everywhere on $\real$.  Since $\rational$ has
measure zero in $\real$, the value of the vector field at rational
points plays no role in the computation of $F[X]$. Hence, the
associated Filippov set-valued map $\map{F[X]}{\real}{\parts{\real}}$
is $F[X](x) = \{-1 \}$.  Since~\eqref{eq:unique-Filipov-sol-I} holds
for all $y,y' \not \in \rational$, there exists a unique solution
starting from each initial condition, more precisely,
$\map{x}{\realnonnegative}{\real}$, where $x (t) = x(0) - t$.
\hfill
$\blacksquare$

The hypotheses of Proposition~\ref{prop:uniqueness-Filippov-I} are
somewhat restrictive.  This assertion is justified by the observation
that, for $d>1$, piecewise continuous vector fields on $\real^d$ are
not essentially one-sided Lipschitz.  We justify this assertion in
``Uniqueness of Filippov Solutions of Piecewise Continuous Vector
Fields.''  Figure~\ref{fig:Lipschitz-type-counter-example} shows an
example of a piecewise continuous vector field with unique solutions
starting from each initial condition.  However, this uniqueness cannot
be guaranteed by means of
Proposition~\ref{prop:uniqueness-Filippov-I}.

The following result~\cite{AFF:88} identifies sufficient conditions
for uniqueness specifically tailored for piecewise continuous vector
fields.

\begin{proposition}\label{prop:uniqueness-Filippov-II}
  Let $\map{X}{\real^d}{\real^d}$ be a piecewise continuous vector
  field, with $\real^d = \domain_1 \union \domain_2$. Let $S_X =
  \boundary{\domain_1} = \boundary{\domain_2}$ be the set of points at
  which $X$ is discontinuous, and assume that $S_X$ is a
  $C^2$-manifold.  Furthermore, assume that, for $i \in \{1, 2\}$,
  $X_{|\ov{\domain_i}}$ is continuously differentiable on $\domain_i$
  and $X_{|\ov{\domain_1}} - X_{|\ov{\domain_2}}$ is continuously
  differentiable on $S_X$. If, for each $x \in S_X$, either
  $X_{|\ov{\domain_1}}(x)$ points into $\domain_2$ or
  $X_{|\ov{\domain_2}}(x)$ points into $\domain_1$, then there exists
  a unique Filippov solution of~\eqref{eq:ODE-auto} starting from each
  initial condition.
\end{proposition}

Note that the continuous differentiability hypothesis on $X$ already
guarantees uniqueness of solutions on each of the sets $\domain_1$ and
$\domain_2$. Roughly speaking, the additional assumptions on $X$ along
$S_X$ in Proposition~\ref{prop:uniqueness-Filippov-II} guarantee that
uniqueness on $\real^d$ is not disrupted by the discontinuities.
Under the stated assumptions, when reaching $S_X$, Filippov solutions
can cross it or slide along it. The situation depicted in
Figure~\ref{fig:sliding}(c) is thus ruled out.

\subsubsection*{Examples~\ref{ex:brick}-\ref{ex:away} revisited:
  Unique Filippov solutions for the sliding brick, the nonsmooth
  harmonic oscillator, and the move-away-from-nearest-neighbor
  interaction law}

As an application of Proposition~\ref{prop:uniqueness-Filippov-II},
let us reconsider Example~\ref{ex:brick}. At $v=0$, the vector
$X_{|\ov{\domain_1}}(0)$ points into $\domain_{2}$, see
Figure~\ref{fig:brick}(b).
Proposition~\ref{prop:uniqueness-Filippov-II} then ensures that there
exists a unique Filippov solution starting from each initial
condition.
  
For Example~\ref{ex:oscillator}, the vector $X_{|\ov{\domain_1}}$
points into $\domain_{2}$ at every point in $S_X \cap \{x_2>0 \}$,
while the vector $X_{|\ov{\domain_2}}$ points into $\domain_{1}$ at
every point in $S_X \cap \{x_2<0 \}$, see
Figure~\ref{fig:oscillator}(a).  Moreover, there is only one solution
(the equilibrium solution) starting from $(0,0)$.  Therefore, using
Proposition~\ref{prop:uniqueness-Filippov-II}, we conclude that this
system has a unique Filippov solution starting from each initial
condition.

For Example~\ref{ex:away}, it is convenient to define $\domain_5 =
\domain_1$.  Then, at $(x_1,x_2) \in \boundary{\domain_i} \cap
\boundary{\domain_{i+1}} \! \setminus \! \{(0,0) \}$, with $i \in
\until{4}$, the vector $X_{|\ov{\domain_i}}(x_1,x_2)$ points into
$\domain_{i+1}$, and the vector $X_{|\ov{\domain_{i+1}}}(x_1,x_2)$
points into $\domain_{i}$, see Figure~\ref{fig:ex1}(a).  Moreover,
there is only one solution (the equilibrium solution) starting from
$(0,0)$.  Therefore, using
Proposition~\ref{prop:uniqueness-Filippov-II}, we conclude that
Example~\ref{ex:away} has a unique Filippov solution starting from
each initial condition.
\hfill
$\blacksquare$

Proposition~\ref{prop:uniqueness-Filippov-II} can also be applied to
piecewise continuous vector fields with an arbitrary number of
partitioning domains, provided that the set where the vector field is
discontinuous is a disjoint union of surfaces resulting from pairwise
intersections of the boundaries of pairs of domains.  Alternative
versions of this result can also be stated for time-varying piecewise
continuous vector fields, as well as for situations in which more than
two domains intersect at a point of discontinuity~\cite[Theorem 4 at
page 115]{AFF:88}.

The literature contains additional results guaranteeing uniqueness of
Filippov solutions tailored to specific classes of dynamical systems.
For instance,~\cite{AYP-WPMHH-HN:03} studies uniqueness for relay
linear systems,~\cite{MMP-PAI:93} investigates uniqueness for adaptive
control systems, while~\cite{RMC-AM:03} establishes uniqueness for a
class of discontinuous differential equations whose vector field
depends on the solution of a scalar conservation law.

\subsection*{Solutions of control systems with discontinuous
  inputs}

Let $\map{X}{\real^d \times \Uc}{\real^d}$, where $\Uc \subseteq
\real^m$ is the set of allowable control-function values, and consider
the control equation on $\real^d$ given by
\begin{align}\label{eq:control-equation}
  \dot{x}(t) = X(x(t),u (t)) .
\end{align}
At first sight, a natural way to identify a notion of solution
for~\eqref{eq:control-equation} is to select a control input, either
an open-loop $\map{u}{\realnonnegative}{\Uc}$, a closed-loop
$\map{u}{\real^d}{\Uc}$, or a combination $\map{u}{\realnonnegative
  \times \real^d}{\Uc}$, and then consider the resulting differential
equation.  When the selected control input $u$ is a discontinuous
function of $x\in \real^d$, then we can consider the solution notions
of Caratheodory or Filippov. At least two alternatives are considered
in the literature. We discuss them next.

\bigskip
\noindent \textbf{\textit{Solutions by means of differential inclusions}}
\medskip

A first alternative to defining a solution notion consists of
associating a differential inclusion with the control
equation~\eqref{eq:control-equation}. In this approach, the set-valued
map $\map{G[X]}{\real^d}{\parts{\real^d}}$ is defined by
\begin{align}\label{eq:control-set-valued-map}
  G[X](x) \triangleq \setdef{X(x,u)}{u \in \Uc} .
\end{align}
\sindex{ControlSetValuedMap}{{$G[X]$}}{Set-valued map associated with
  a control system $\map{X}{\real^d \times \Uc}{\real^d}$}%
In other words, the set $G[X](x)$ captures all of the directions in
$\real^d$ that can be generated at $x$ with controls belonging to
$\Uc$.  Consider now the differential inclusion
\begin{align}\label{eq:control-inclusion}
  \dot{x}(t) \in G[X](x(t)) .
\end{align}
A \emph{solution of~\eqref{eq:control-equation}} on $[0,t_1] \subset
\real$ is defined to be a Caratheodory solution of the differential
inclusion~\eqref{eq:control-inclusion}, that is, an absolutely
continuous map $\map{x}{[0,t_1]}{\real^d}$ such that $\dot{x}(t) \in
G[X] (x(t))$ for almost all $t \in [0,t_1]$.

If we choose an open-loop input $\map{u}{\realnonnegative}{\Uc}$
in~\eqref{eq:control-equation}, then a Caratheodory solution of the
resulting dynamical system is also a Caratheodory solution of the
differential inclusion~\eqref{eq:control-inclusion}.  Alternatively,
it can be shown~\cite{AFF:88} that, if $X$ is continuous and $\Uc$ is
compact, then the converse is also true.  The differential
inclusion~\eqref{eq:control-inclusion} has the advantage of not
focusing attention on a particular control input, but rather allows us
to comprehensively study and understand the properties of the control
system.

\bigskip
\noindent \textbf{\textit{Sample-and-hold solutions}} 
\medskip

A second alternative to defining a solution notion for the control
equation~\eqref{eq:control-equation} uses the notion of
sample-and-hold solution~\cite{FHC-YSL-EDS-AIS:97}. As discussed in
the section ``Stabilization of control systems,'' this notion plays a
key role in the stabilization of asymptotically controllable systems.

\index{partition}%
\index{partition!diameter}%
\sindex{Partition}{{$\pi$}}{Partition of a closed interval}%
\sindex{DiameterPartition}{{$\diam (\pi)$}}{Diameter of the partition $\pi$}%

A \emph{partition} of the interval $[t_0,t_1]$ is an increasing
sequence $\pi = \{s_i\}_{i=0}^N$ with $s_0 = t_0$ and $s_N = t_1$.
The partition need not be finite.  The notion of a partition of
$[t_0,\infty)$ is defined similarly.  The diameter of $\pi$ is $ \diam
(\pi) \triangleq \sup \setdef{s_{i}-s_{i-1}}{i \in \until{N}}$.  Given
a control input $\map{u}{\realnonnegative \times \real^d}{\Uc}$, an
initial condition $x_0$, and a partition $\pi$ of $[0,t_1]$, a
\emph{$\pi$-solution
  of~\eqref{eq:control-equation}}\index{solution!sample-and-hold}
defined on $[0,t_1] \subset \real$ is the map
$\map{x}{[0,t_1]}{\real^d}$, with $x(0) = x_0$, recursively defined by
requiring the curve $ t \in [s_{i-1},s_{i}] \mapsto x(t)$, for $ i \in
\until{N-1}$, to be a Caratheodory solution of the differential
equation
\begin{align}
  \dot x (t) = X(x(t),u(s_{i-1},x(s_{i-1}))) .
\end{align}
$\pi$-solutions are also referred to as \emph{sample-and-hold}
solutions because the control is held fixed throughout each interval
of the partition at the value according to the state at the beginning
of the interval.  Figure~\ref{fig:sample-and-hold} shows an example of
a $\pi$-solution. The existence of $\pi$-solutions is
guaranteed~\cite{FHC-YSL-RJS-PRW:98} if, for all $u \in \Uc \subseteq
\real^m$, the map $x \mapsto X(x,u)$ is continuous.

\section*{Nonsmooth Analysis}

We now consider candidate nonsmooth Lyapunov functions for
discontinuous differential equations.  The level of generality
provided by nonsmooth analysis is not always necessary.  The stability
properties of some discontinuous dynamical systems and differential
inclusions can be analyzed with smooth functions, as the following
example shows.

\subsubsection*{Examples~\ref{ex:discontinuous-non-existence}
  and~\ref{ex:discontinuous-existence} revisited: Stability of the
  origin for the sign function}

We have already established that the vector fields~\eqref{eq:X1}
and~\eqref{eq:X2} in examples~\ref{ex:discontinuous-non-existence}
and~\ref{ex:discontinuous-existence}, respectively, have unique
Filippov solutions starting from each initial condition. Now, consider
the smooth Lyapunov function $\map{f}{\real}{\real}$, where $f(x) =
x^2/2$. Now, for all $x \neq 0$, we have
\begin{align*}
  \nabla f (x) \cdot X(x) = - |x| < 0.
\end{align*}
Since, according to~\eqref{eq:Filippov-set-sign}, $F[X](x) = \{X(x)
\}$ on $\real \! \setminus \! \{0\}$, we deduce that the function $f$ is
decreasing along every Filippov solution of~\eqref{eq:X1}
and~\eqref{eq:X2} that starts on $\real \! \setminus \! \{0\}$. Therefore,
we conclude that the equilibrium $x=0$ is globally asymptotically
stable.
\hfill $\blacksquare$

However, nonsmooth Lyapunov functions may be needed if dealing with
discontinuous dynamics, as the following example shows.

\subsubsection*{Example~\ref{ex:oscillator} revisited: The nonsmooth
  harmonic oscillator does not admit a smooth Lyapunov function}

Consider the vector field for the nonsmooth harmonic oscillator in
Example~\ref{ex:oscillator}.  We reason as in~\cite{FC:99}.  As stated
in the section ``Piecewise continuous vector fields and sliding
motions,'' all of the Filippov solutions
of~\eqref{eq:nonsmooth-harmonic-oscillator} are periodic, and
correspond to the curves plotted in Figure~\ref{fig:oscillator}(b).
Therefore, if a smooth Lyapunov function $\map{f}{\real^2}{\real}$
exists, then it must be constant on each Filippov solution.  Since the
level sets of $f$ must be one-dimensional, it follows that each curve
must be a level set.  This property contradicts the fact that the
function is smooth, since the level sets of a smooth function are also
smooth, and the Filippov solutions plotted in
Figure~\ref{fig:oscillator}(b) are not smooth along the
$x_2$-coordinate axis.
\hfill $\blacksquare$

The above example illustrates the need to consider nonsmooth analysis.
Similar examples can be found in~\cite[Section~2.2.2]{AB-LR:05}.  It
is also worth noting that the presentation in ``Stability analysis by
means of the generalized gradient of a nonsmooth Lyapunov function''
boils down to classical stability analysis when the candidate Lyapunov
function is smooth.

In this section we discuss two tools from nonsmooth analysis, namely,
the generalized gradient and the proximal
subdifferential~\cite{FHC:83,FHC-YSL-RJS-PRW:98}. As with the notions
of solution of discontinuous differential equations, multiple
generalized derivative notions are available in the literature when a
function fails to be differentiable. These notions include, in
addition to the two considered in this section, the generalized (super
or sub) differential, the (upper or lower, right or left) Dini
derivative, and the contingent
derivative~\cite{FHC-YSL-RJS-PRW:98,RTR-RJBW:98,JPA-HF:90,JPA-AC:94}.
Here, we focus on the notions of the generalized gradient and proximal
subdifferential because of their role in providing stability tools for
discontinuous differential equations.
\subsection*{The generalized gradient of a locally Lipschitz
  function}

Rademacher's theorem~\cite{FHC:83} states that every locally Lipschitz
function is differentiable almost everywhere in the sense of Lebesgue
measure.  When considering a locally Lipschitz function as a candidate
Lyapunov function, this statement may raise the question of whether to
disregard those points where the gradient does not exist.
Conceivably, the solutions of the dynamical system stay for almost all
time away from ``bad'' points where no gradient of the function
exists.  However, this assumption is not always valid.  As we show in
Example~\ref{ex:sm} below, the solutions of a dynamical system may
insist on staying on the ``bad'' points forever.  In that case, having
some gradient-like information is helpful.
\sindex{OmegaSet}{{$\Omega_f$}}{Set of points where the
  locally Lipschitz function $\map{f}{\real^d}{\real}$ fails to be
  differentiable}%

Let $\map{f}{\real^d}{\real}$ be a locally Lipschitz function, and let
$\Omega_f \subset \real^d$ denote the set of points where $f$ fails to
be differentiable. The \emph{generalized gradient}
\index{generalized gradient}
$\map{\partial f}{\real^d}{\parts{\real^d}}$ of $f$ is defined by
\begin{align}\label{eq:generalized-gradient}
  \partial f (x) \triangleq \co \setdef{\lim_{i \rightarrow \infty}
    \nabla f (x_i)}{x_i \rightarrow x \, , \; x_i \not \in S \cup
    \Omega_f} ,
\end{align}
\sindex{ConvexHull}{{$\co(S)$}}{Convex hull of a set $S \subseteq
  \real^d$}%
\sindex{Gradient}{$\nabla f$}{Gradient of the differentiable function
  $\map{f}{\real^d}{\real}$}%
\sindex{GeneralizedGradient}{$\partial f$}{Generalized gradient of the
  locally Lipschitz function $\map{f}{\real^d}{\real}$}%
where $\co$ denotes convex hull.  In this definition, $S \subset
\real^d$ is a set of measure zero that can be arbitrarily chosen to
simplify the computation. The resulting set $\partial f (x)$ is
independent of the choice of $S$.  From the definition, the
generalized gradient of $f$ at $x$ consists of all convex combinations
of all of the possible limits of the gradient at neighboring points
where $f$ is differentiable.  Equivalent definitions of the
generalized gradient are given in~\cite{FHC:83}.

Some useful properties of the generalized gradient are summarized in
the following result~\cite{FHC-YSL-RJS-PRW:98}.

\begin{proposition}\label{prop:generalized-gradient-properties}
  If $\map{f}{\real^d}{\real}$ is locally Lipschitz at $x\in \real^d$,
  then the following statements hold:
  \begin{enumerate}
  \item $\partial f (x)$ is nonempty, compact, and convex.
  \item The set-valued map $\map{\partial
      f}{\real^d}{\parts{\real^d}}$, $x\mapsto \partial f(x)$, is
    upper semicontinuous and locally bounded at $x$.
  \item If $f$ is continuously differentiable at $x$, then $\partial f
    (x) = \{ \nabla f (x) \}$.
  \end{enumerate}
\end{proposition}

Let us compute the generalized gradient for an illustrative example.

\refstepcounter{example}%
\subsubsection*{Example~\theexample: Generalized gradient of the
  absolute value function}\label{ex:abs}

Consider the locally Lipschitz function $\map{f}{\real}{\real}$, where
$f(x) = |x|$, which is continuously differentiable everywhere except
at $0$.  Since $\nabla f (x) = 1$ for $x>0$ and $\nabla f (x) = -1$
for $x<0$, from~\eqref{eq:generalized-gradient}, we deduce
\begin{align*}
  \partial f (0) = \co \{1,-1 \} = [-1,1] ,
\end{align*}
while, by
Proposition~\ref{prop:generalized-gradient-properties}\textit{(iii)},
we deduce $\partial f(x) = \{1\}$ for $x>0$ and $\partial f(x) =
\{-1\}$ for $x<0$.
\hfill
$\blacksquare$

The next example shows that
Proposition~\ref{prop:generalized-gradient-properties}\textit{(iii)}
may not be true if $f$ is not continuously differentiable.

\refstepcounter{example}%
\subsubsection*{Example~\theexample: Generalized gradient of a
  differentiable but not continuously differentiable function
}\label{ex:not-continuously-differentiable}

Consider the locally Lipschitz function $\map{f}{\real}{\real}$, where
$f(x) = x^2 \sin (\frac{1}{x})$, which is continuously differentiable
everywhere except at $0$, where it is only differentiable. It can be
shown that $\nabla f (0) = 0$, while $\nabla f (x) = 2 x \sin
(\frac{1}{x}) - \cos(\frac{1}{x})$ for $x \neq 0$.
From~\eqref{eq:generalized-gradient}, we deduce that $\partial f (0) =
[-1,1]$, which is different from $\{ \nabla f (0) \} = \{0\}$. This
example illustrates that
Proposition~\ref{prop:generalized-gradient-properties}\textit{(iii)}
may not be valid if $f$ is not continuously differentiable at the
point $x$.
\hfill
$\blacksquare$

\bigskip
\noindent \textbf{\textit{Computing the generalized gradient}}
\medskip

Computation of the generalized gradient of a locally Lipschitz
function is often a difficult task. In addition to the brute force
approach of working from the
definition~\eqref{eq:generalized-gradient}, various results are
available to facilitate this computation.  Many results that are valid
for ordinary derivatives have a counterpart in this setting.  We
summarize some of these results here, and refer the reader
to~\cite{FHC:83,FHC-YSL-RJS-PRW:98} for a complete exposition.  In the
statements of these results, the notion of a regular function plays a
prominent role, see ``Regular Functions'' for a precise definition.

\textbf{Dilation rule.} If $\map{f}{\real^d}{\real}$ is locally
Lipschitz at $x \in \real^d$ and $s\in \real$, then the function $ s
f$ is locally Lipschitz at $x$, and
\begin{align}\label{eq:gradient-dilation}
  \partial (s f)(x) = s \; \partial f(x).
\end{align}

\textbf{Sum rule.} If $\map{f_1,f_2}{\real^d}{\real}$ are locally
Lipschitz at $x \in \real^d$ and $s_1,s_2 \in \real$, then the
function $s_1 f_1 + s_2 f_2$ is locally Lipschitz at $x$, and
\begin{align}\label{eq:gradient-sum}
  \partial \big( s_1 f_1 + s_2 f_2\big)(x) \subseteq s_1
  \partial f_1(x) + s_2 \partial f_2(x) ,
\end{align}
where the sum of two sets $A_1$ and $A_2$ is defined as $A_1 + A_2 =
\setdef{a_1 + a_2}{a_1 \in A_1, \; a_2 \in A_2}$.  Moreover, if $f_1$
and $f_2$ are regular at $x$, and $s_1,s_2 \in \realnonnegative$, then
equality holds and $s_1 f_1 + s_2 f_2$ is regular at $x$.

\textbf{Product rule.} If $\map{f_1,f_2}{\real^d}{\real}$ are locally
Lipschitz at $x \in \real^d$, then the function $f_1 f_2$ is locally
Lipschitz at $x$, and
\begin{equation}\label{eq:gradient-product}
  \partial \big( f_1 f_2\big)(x) \subseteq f_2 (x) 
  \partial f_1(x) + f_1(x) \partial f_2(x).
\end{equation}
Moreover, if $f_1$ and $f_2$ are regular at $x$, and $f_1 (x), f_2(x)
\ge 0$, then equality holds and $f_1 f_2$ is regular at $x$.

\textbf{Quotient rule.} If $\map{f_1,f_2}{\real^d}{\real}$ are locally
Lipschitz at $x \in \real^d$, and $f_2(x) \neq 0$, then the function
$f_1 /f_2$ is locally Lipschitz at $x$, and
\begin{equation}\label{eq:gradient-quotient}
  \partial \Big( \frac{f_1}{f_2} \Big)(x) \subseteq \frac{1}{f_2^2
    (x)} \Big(f_2 (x) 
    \partial f_1(x) - f_1(x) \partial f_2(x) \Big) .
\end{equation}
Moreover, if $f_1$ and $-f_2$ are regular at $x$, and $f_1 (x) \ge 0$
and $f_2(x) > 0$, then equality holds and $f_1 /f_2$ is regular at
$x$.

\textbf{Chain rule.} If each component of $\map{h}{\real^d}{\real^m}$
is locally Lipschitz at $x \in \real^d$ and $\map{g}{\real^m}{\real}$
is locally Lipschitz at $h(x) \in \real^m$, then the function $ g
\circ h$ is locally Lipschitz at $x$, and
\begin{equation}\label{eq:gradient-chain}
  \partial \big( g \circ h \big)(x) \subseteq \co \Big \{ \sum_{k=1}^m
  \alpha_k \zeta_k \; : \;     (\alpha_1,\dots,\alpha_m) \in
  \partial g(h(x)) , \; (\zeta_1,\dots,\zeta_m)  \in \partial h_1
  (x)  \times \dots \times \partial h_m (x) \Big\}.
\end{equation}
Moreover, if $g$ is regular at $h(x)$, each component of $h$ is
regular at $x$, and $\partial g (h(x)) \subset \realnonnegative^m$,
then equality holds and $ g \circ h$ is regular at $x$.

The following useful result from~\cite[Proposition~2.3.12]{FHC:83}
concerns the generalized gradient of the maximum and the minimum of a
finite set of functions.
\begin{proposition}\label{prop:gradients-max-min}
  For $k \in \until{m}$, let $f_k:\real^d \rightarrow \real$ be
  locally Lipschitz at $x \in \real^d$, and define the functions
  $\map{\fmax,\fmin}{\real^d}{\real}$ by
  \begin{align*}
    \fmax(y) & \triangleq \max \setdef{f_k(y)}{k \in \until{m}}, \\
    \fmin(y) & \triangleq \min \setdef{f_k(y)}{k \in \until{m}} .
  \end{align*}
  Then, the following statements hold:
  \begin{enumerate}
  \item $\fmax$ and $\fmin$ are locally Lipschitz at $x$.
  \item Let $\Imax(x)$ denote the set of indices $k$ for which
    $f_k(x)=\fmax(x)$. Then
    \begin{align}\label{eq:gradient-max}
      \partial \fmax (x) \subseteq \co \bigcup \setdef{\partial f_i
        (x)}{i \in \Imax(x)} .
    \end{align}
    Furthermore, if $f_i$ is regular at $x$ for all $i \in \Imax(x)$,
    then equality holds and $\fmax$ is regular at~$x$,
  \item Let $\Imin(x)$ denote the set of indices $k$ for which
    $f_k(x)=\fmin(x)$. Then
    \begin{align}\label{eq:gradient-min}
      \partial \fmin (x) \subseteq \co \bigcup \setdef{\partial f_i
        (x)}{i \in \Imin(x)} .
    \end{align}
    Furthermore, if $-f_i$ is regular at $x$ for all $i \in \Imin(x)$,
    then equality holds and $-\fmin$ is regular at~$x$.
  \end{enumerate}
\end{proposition}

It follows from Proposition~\ref{prop:gradients-max-min} that the
maximum of a finite set of continuously differentiable functions is a
locally Lipschitz and regular function whose generalized gradient at
each point $x$ is easily computable as the convex hull of the
gradients of the functions that attain the maximum at $x$.

\subsubsection*{Example~\ref{ex:abs} revisited: Generalized gradient
  of the absolute value function}

The absolute value function $f(x) = |x|$ can be rewritten as $f(x) =
\max \{ x , -x \}$.  Both $x \mapsto x$ and $x \mapsto -x$ are
continuously differentiable, and hence locally Lipschitz and regular.
Therefore, according to
Proposition~\ref{prop:gradients-max-min}\textit{(i)}
and~\textit{(ii)}, $f$ is locally Lipschitz and regular, and its
generalized gradient is
\begin{align}\label{eq:max-gradient}
  \partial f (x) =
  \begin{cases}
    \{1\}, & x>0,\\
    [-1,1], & x=0,\\
    \{-1\}, & x<0 .
  \end{cases}
\end{align}
This result is obtained in Example~\ref{ex:abs} by direct computation.
\hfill $\blacksquare$

\refstepcounter{example}%
\subsubsection*{Example~\theexample: Generalized gradient of the minus
  absolute value function}\label{ex:minus-abs}

The minimum of a finite set of regular functions is not always
regular.  A simple example is given by $g(x) = \min \{ x, -x\} =
-|x|$, which is not regular at $0$, as we show in ``Regular
Functions.'' However, according to
Proposition~\ref{prop:gradients-max-min}\textit{(iii)}, this fact does
not mean that its generalized gradient cannot be computed. Indeed,
\begin{align}\label{eq:min-gradient}
  \partial g (x) =
  \begin{cases}
    \{-1\}, & x>0,\\
    [-1,1], & x=0,\\
    \{1\}, & x<0 .
  \end{cases}
\end{align}
This result can also be obtained by combining~\eqref{eq:max-gradient}
with the application of the dilation rule~\eqref{eq:gradient-dilation}
to the function $f(x) = |x|$ with the parameter $s=-1$.
\hfill $\blacksquare$

\bigskip
\noindent \textbf{\textit{Critical points and directions of
  descent}}
\medskip

A \emph{critical}\index{critical point} point of
$\map{f}{\real^d}{\real}$ is a point $x \in\real^d$ such that $ 0 \in
\partial f(x)$.  According to this definition, the maximizers and
minimizers of a locally Lipschitz function are critical points. As an
example, $x=0$ is the unique minimizer of $f(x) = |x|$, and, indeed,
we see that $0 \in
\partial f (0)$ in~\eqref{eq:max-gradient}.

If a function $f$ is continuously differentiable, then the gradient
$\nabla f$ provides the direction of maximum ascent of $f$
(respectively, $-\nabla f$ provides the direction of maximum descent
of $f$).  When we consider a locally Lipschitz function, however, the
question of choosing the directions of descent among all of the
available directions in the generalized gradient arises.  Without loss
of generality, we restrict our discussion to directions of descent,
since a direction of descent of $-f$ corresponds to a direction of
ascent of $f$, while $f$ is locally Lipschitz if and only if $-f$ is
locally Lipschitz.

Let $\map{\LN}{\parts{\real^d}}{\parts{\real^d}}$ be the set-valued
map that associates to each subset $S$ of $\real^d$ the set of
least-norm elements of its closure $\overline{S}$.  If the set $S$ is
convex and closed, then the set $\LN (S)$ is a singleton, which
consists of the orthogonal projection of $0$ onto $S$.  For a locally
Lipschitz function $f$, consider the \emph{generalized gradient vector
  field}\index{vector field!generalized gradient} $ \map{\LN (\partial
  f)}{\real^d}{\real^d}$ defined by
\begin{align*}
  x \mapsto \LN (\partial f) (x) \triangleq \LN (\partial f (x)) .
\end{align*}
It turns out that $-\LN (\partial f) (x)$ is a direction of descent of
$f$ at $x \in \real^d$. More precisely~\cite{FHC:83},
\sindex{LeastNorm}{{$\LN (S)$}}{Least-norm elements in the closure of
  the set $S \subseteq \real^d$}%
if $0 \not \in \partial f (x)$, then there exists $T>0$ such that
\begin{equation}\label{eq:not-critical}
  f(x-t  \LN (\partial f)(x)) \le f(x) - \frac{t}{2} \TwoNorm{\LN
    (\partial f)(x)}^2 < f(x) , \quad 0  < t <  T ,
\end{equation}
that is, by taking a small step in the direction $-\LN (\partial f)
(x)$, the function $f$ is guaranteed to decrease by an amount that
scales linearly with the stepsize.

\refstepcounter{example}%
\subsubsection*{Example~\theexample:
  Minimum-distance-to-polygonal-boundary function}\label{ex:sm}

Let $Q \subset \real^2$ be a convex polygon. Consider the minimum
distance function $\map{\sm_Q}{Q}{\real}$ from a point within the
polygon to its boundary defined by
\sindex{SmallestDistanceToBoundary}{$\sm_Q$}{Minimum distance
  function from a point in a convex polygon $Q \subset \real^d$ to
  the boundary of $Q$}
\begin{align*}
  \sm_Q (p) \triangleq \min \setdef{\TwoNorm{p-q}}{q \in \boundary{Q}}
  .
\end{align*}
Note that the value of $\sm_Q$ at $p$ is the radius of the largest
disk with center $p$ contained in the polygon.  Moreover, the function
$\sm_Q$ is locally Lipschitz on $Q$. To show this, rewrite $\sm_Q$ as
\begin{align*}
  \sm_Q (p) \triangleq \min \setdef{\D (p,e)}{e \; \text{is an edge of
      $Q$}} ,
\end{align*}
where $\D (p,e)$ denotes the Euclidean distance from the point $p$
to the edge $e$.
\sindex{Distance}{$\D(p,S)$}{Euclidean distance from the point $p
  \in \real^d$ to the set $S \subseteq \real^d$}
Indeed, the function $\sm_Q$ is concave on $Q$.

Let us consider the generalized gradient vector field corresponding to
$\sm_Q$, where the definition of $\sm_Q$ is extended outside $Q$ by
setting $\sm_Q(p) = - \min \setdef{\TwoNorm{p-q}}{q \in \boundary{Q}}$
for $p \not \in Q$.  Applying
Proposition~\ref{prop:gradients-max-min}\textit{(iii)}, we deduce that
$-\sm_Q$ is regular on $Q$ and its generalized gradient at $ p \in Q$
is
\begin{align*}
  \partial \sm_Q (p) &= \co \setdef{\normal_e}{e \; \text{edge of $Q$
      such that $ \sm_Q (p)= \D(p,e)$} } ,
\end{align*}
where $\normal_e$ denotes the unit normal to the edge $e$ pointing
toward the interior of $Q$.
\sindex{Normal}{$\normal_e$}{Unit normal to the edge $e$ of a polygon
  $Q$ pointing toward the interior of $Q$}
Therefore, at points $p$ in $Q$ for which there exists a unique edge
$e$ of $Q$ nearest to $p$, the function $\sm_Q$ is differentiable, and
its generalized gradient vector field is given by $\LN (\sm_Q)(p) =
\normal_e$.  Note that this vector field corresponds to the
move-away-from-nearest-neighbor interaction law for one agent moving
in the polygon introduced in Example~\ref{ex:away}! 

At points $p$ of $Q$ where various edges $\{e_1,\dots,e_m\}$ are at
the same minimum distance to $p$, the function $\sm_Q$ is not
differentiable, and its generalized gradient vector field is given by
the least-norm element in $\co \{\normal_{e_1},\dots,\normal_{e_m}\}$.
If $p$ is not a critical point, $0$ does not belong to $\co
\{\normal_{e_1},\dots,\normal_{e_m}\}$, and the least-norm element
points in the direction of the bisector line between two of the edges
in $\{e_1,\dots,e_m\}$.  Figure~\ref{fig:sm} shows a plot of the
generalized gradient vector field of $\sm_Q$ on the square $Q =
[-1,1]^2$. Note the similarity with the plot in
Figure~\ref{fig:ex1}(a).  Indeed, the critical points of $\sm_Q$ are
characterized~\cite{JC-FB:02m} by the statement that $0 \in
\partial \sm_Q(p)$ if and only if $p$ belongs to the incenter set
of~$Q$.

The incenter set of $Q$ is composed of the centers of the
largest-radius disks contained in~$Q$. In general, the incenter set is
a segment (consider the example of a rectangle).  However, if $0 \in
\text{interior} ( \partial \sm_Q(p))$, then the incenter set of $Q$ is
the singleton $\{p\}$.
\hfill
$\blacksquare$

\bigskip
\noindent \textbf{\textit{Nonsmooth gradient
    flows}}
\medskip

Given a locally Lipschitz function $\map{f}{\real^d}{\real}$, the
nonsmooth analog of the classical gradient descent flow of a
differentiable function is defined by
\begin{equation}\label{eq:natural-gradient}
  \dot{x}(t) = - \LN (\partial f)(x(t)) . 
\end{equation}
According to~\eqref{eq:not-critical}, unless the flow is already at a
critical point, $-\LN (\partial f)(x)$ is a direction of descent of
$f$ at $x$. Since this nonsmooth gradient vector field is
discontinuous, we have to specify a notion of solution
for~\eqref{eq:natural-gradient}. In this case, we use the notion of
Filippov solution due largely to the remarkable fact~\cite{BP-SSS:87}
that the Filippov set-valued map associated with the nonsmooth
gradient flow of $f$ given by~\eqref{eq:natural-gradient} is precisely
the generalized gradient of the function, as the next result states.

\begin{proposition}[Filippov set-valued map of nonsmooth
  gradient]\label{prop:Filippov-set-valued-nonsmooth-gradient}
  If $\map{f}{\real^d}{\real}$ is locally Lipschitz, then the Filippov
  set-valued map $\map{F[\LN (\partial f)]}{\real^d}{\parts{\real^d}}$
  of the nonsmooth gradient of $f$ is equal to the generalized
  gradient $\map{\partial f}{\real^d}{\parts{\real^d}}$ of $f$, that
  is, for $x \in \real^d$,
  \begin{align}\label{eq:magic-Filippov}
    F[\LN (\partial f)](x) = \partial f (x) .
  \end{align}
\end{proposition}

As a consequence of
Proposition~\ref{prop:Filippov-set-valued-nonsmooth-gradient}, the
discontinuous system~\eqref{eq:natural-gradient} is equivalent to the
differential inclusion
\begin{align}
  \dot{x}(t) \in - \partial f (x(t)) .
\end{align}
Solutions of~\eqref{eq:natural-gradient} are sometimes called ``slow
motions'' of the differential inclusion because of the use of the
least-norm operator, which selects the vector in the generalized
gradient of the function with the smallest magnitude.  How can we
analyze the asymptotic behavior of the solutions of this system? If
the function~$f$ is differentiable, then the invariance principle
allows us to deduce that, for functions with bounded level sets, the
solutions of the gradient flow converge to the set of critical points.
The key tool behind this result is to establish that the function
decreases along solutions of the system.  This behavior is formally
expressed through the notion of the Lie derivative.  In the section
``Nonsmooth Stability Analysis'' we discuss generalizations of the
notion of Lie derivative to the nonsmooth case.  These notions allow
us to study the asymptotic convergence properties of the trajectories
of nonsmooth gradient flows.

\subsection*{The proximal subdifferential of a lower semicontinuous
  function}

A complementary set of nonsmooth analysis tools for dealing with
Lyapunov functions arises from the concept of proximal
subdifferential. This concept has the advantage of being defined for
the class of lower semicontinuous functions, which is larger than the
class of locally Lipschitz functions.  The generalized gradient
provides us with directional descent information, that is, directions
along which the function decreases.  The price we pay for using the
proximal subdifferential is that explicit descent directions are not
generally known to us. Nevertheless, the proximal subdifferential
allows us to reason about the monotonic properties of the function,
which, as we show in the section ``Stability analysis by means of the
proximal subdifferential of a nonsmooth Lyapunov function,'' turn out
to be sufficient to provide stability results.  The proximal
subdifferential is particularly appropriate when dealing with convex
functions.  Here, we briefly touch on the topic of convex analysis.
For further information, see~\cite{RTR:70,JMB-ASL:00,JBHU-CL:04}.

A function $\map{f}{\real^d}{\real}$ is \emph{lower
  semicontinuous}\index{function!lower
  semicontinuous}~\cite{FHC-YSL-RJS-PRW:98} at $x \in \real^d$ if, for
all $\eps \in \realpositive$, there exists $\delta \in \realpositive$
such that, for $ y \in B(x,\delta)$, $f(y) \ge f(x) - \eps$.  The
\emph{epigraph} of $f$ is the set of points lying on or above its
graph, that is, $\epigraph{f} = \setdef{(x,\mu) \in \real^d \times
  \real}{f(x) \le \mu} \subset \real^{d+1}$. The function~$f$ is lower
semicontinuous if and only if its epigraph is closed.  The function
$\map{f}{\real^d}{\real}$ is \emph{upper
  semicontinuous}\index{function!upper semicontinuous} at $x \in
\real^d$ if $-f$ is lower semicontinuous at $x$. Note that $f$ is
continuous at $x$ if and only if $f$ is both upper semicontinuous and
lower semicontinuous at $x$.

For a lower semicontinuous function $\map{f}{\real^d}{\real}$, the
vector $\zeta \in \real^d$ is a \emph{proximal subgradient of $f$ at
  $x \in \real^d$}
\index{proximal subgradient}
\sindex{Proximal subdifferential}{$\partial_P f$}{Proximal
  subdifferential  of the  lower semicontinuous function
  $\map{f}{\real^d}{\real}$}%
if there exist $\sigma, \delta \in \realpositive$ such that, for all
$y \in \oball{\delta}{x}$,
\begin{align}\label{eq:prox-gradient}
  f(y) \ge f(x) + \zeta (y-x) - \sigma^2 \TwoNorm{y - x}^2.
\end{align}
The set $\partial_P f (x)$ of all proximal subgradients of $f$ at $x$
is the \emph{proximal subdifferential of $f$ at $x$}.
\index{proximal subdifferential}
The proximal subdifferential at $x$, which may be empty, is convex but
not necessarily open, closed, or bounded.

\begin{figure}[htbp]
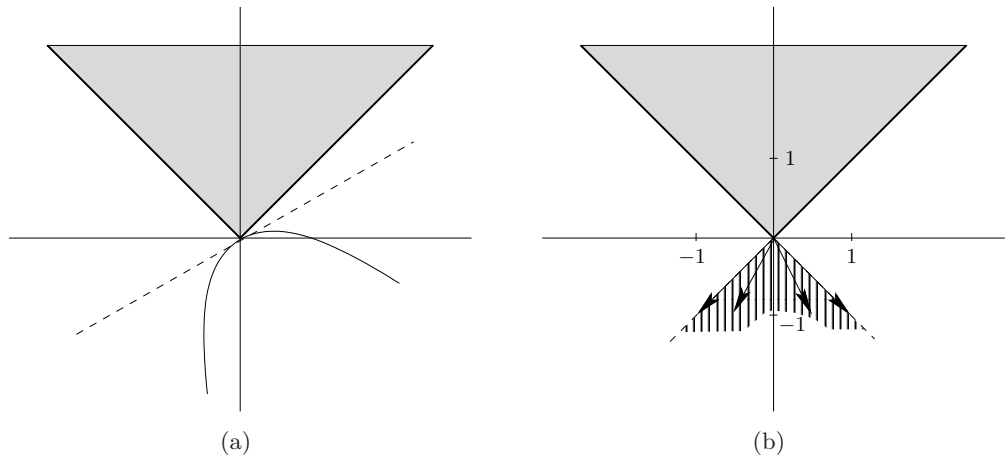

  \centering%
  \subfigure[]{\resizebox{.375\linewidth}{!}{\input{figure10a.tex}}} \qquad
  \subfigure[]{\resizebox{.375\linewidth}{!}{\input{figure10b.tex}}}
  \caption{Geometric interpretations of the proximal subdifferential
    of the function $x \mapsto |x|$ at $x=0$ computed in
    Example~\ref{ex:abs-prox}. The epigraph of the function is shaded,
    while the proximal normal cone to the epigraph at $0$ is striped.
    In (a), according to~\eqref{eq:prox-gradient}, each proximal
    subgradient corresponds to a direction tangent to a parabola that
    fits under the epigraph of the function. In (b), according
    to~\eqref{eq:proximal-differential-normal-cone}, each proximal
    subgradient can be uniquely associated with an element of the
    proximal normal cone to the epigraph of the function.}
  \label{fig:prox-subdifferential}
\end{figure}

Geometrically, the definition of a proximal subgradient can be
interpreted as follows.  Equation~\eqref{eq:prox-gradient} is
equivalent to saying that, in a neighborhood of $x$, the function $y
\mapsto f(y)$ majorizes the quadratic function $y \mapsto f(x) +
\zeta (y-x) - \sigma^2 \TwoNorm{y - x}^2$. In other words, there
exists a parabola that locally fits under the epigraph of $f$ at
$(x,f(x))$. This geometric interpretation is useful for explicitly
computing the proximal subdifferential, see
Figure~\ref{fig:prox-subdifferential}(a) for an illustration.  An
additional geometric characterization of the proximal subdifferential
of $f$ can be given in terms of normal cones and the epigraph of $f$.
Given a closed set $S \subseteq \real^d$ and $x \in S$, the
\emph{proximal normal cone} $N_S(x)$ to $S$ is the set of all $y \in
\real^d$ such that $x$ is the nearest point in $S$ to $x + \lambda y$,
for $\lambda >0$ sufficiently small. Then, $\zeta \in
\partial_P f (x)$ if and only if
\begin{align}\label{eq:proximal-differential-normal-cone}
  (\zeta, -1) \in N_{\epigraph{f}}(x,f(x)) .
\end{align}
Figure~\ref{fig:prox-subdifferential}(b) illustrates this geometric
interpretation.

\refstepcounter{example}%
\subsubsection*{Example~\theexample: Proximal subdifferentials of the
  absolute value function and its negative}\label{ex:abs-prox}
  
Consider the locally Lipschitz functions $\map{f,g}{\real}{\real}$,
$f(x) = |x|$ and $g(x) = - |x|$.  Using the geometric interpretation
of~\eqref{eq:prox-gradient}, it can be seen that
\begin{align}
  \partial_P f(x) &=
  \begin{cases}
    \{ 1 \}, & x <0,\\
    [-1,1], &  x=0,\\
    \{ -1 \}, & x >0 ,
  \end{cases}
  \label{eq:comp-proximal-subdifferential-a}
  \\
  \partial_P g(x) &=
  \begin{cases}
    \{ -1 \}, & x <0,\\
    \emptyset, &  x=0,\\
    \{ 1 \}, & x >0 .
  \end{cases}
  \label{eq:comp-proximal-subdifferential-b}
\end{align}
Note
that~\eqref{eq:comp-proximal-subdifferential-a}-\eqref{eq:comp-proximal-subdifferential-b}
is an example of the fact that $\partial_P ( -f)$ and $-
\partial_P f$ are not necessarily equal.  The generalized gradient of
$g$ in~\eqref{eq:min-gradient} is different from the proximal
subdifferential of $g$ in~\eqref{eq:comp-proximal-subdifferential-b}.
\hfill $\blacksquare$

Unlike the case of the generalized gradient, the proximal
subdifferential may not coincide with $\nabla f(x)$ when $f$ is
continuously differentiable. The function $\map{f}{\real}{\real}$,
where $f(x)= -|x|^{3/2}$, is continuously differentiable with $\nabla
f (0) = 0$, but $\partial_P f (0) = \emptyset$.  In fact, a
continuously differentiable function on $\real$ with an empty proximal
subdifferential almost everywhere is provided
in~\cite{FHC-YSL-PRW:95}.  However, the density theorem
(cf.~\cite[Theorem 3.1]{FHC-YSL-RJS-PRW:98}) states that the proximal
subdifferential of a lower semicontinuous function is nonempty on a
dense set of its domain of definition, although a dense set can have
zero Lebesgue measure.

On the other hand, the proximal subdifferential can be more useful
than the generalized gradient in some situations, as the following
example shows.

\refstepcounter{example}%
\subsubsection*{Example~\theexample: Proximal subdifferential of the
  square root of the absolute
  value}\label{ex:prox-subdifferential-root-absolute}

Consider the function $\map{f}{\real}{\real}$, where $f(x) =
\sqrt{|x|}$, which is continuous at $0$, but not locally Lipschitz at
$0$, which is the global minimizer of $f$.  The generalized gradient
does not exist at $0$, and hence we cannot characterize this point as
a critical point.  However, the function $f$ is lower semicontinuous
with proximal subdifferential,
\begin{align}\label{eq:prox-subdifferential-root-absolute}
  \partial_P f (x) = 
  \begin{cases}
    \big \{ \frac{1}{2\sqrt{x}} \big \}, & x>0
    ,\\
    \real , & x=0 ,\\
    \big \{ -\frac{1}{2 \sqrt{-x}} \big \}, & x<0 .
  \end{cases}
\end{align}
As we discuss later, it follows
from~\eqref{eq:prox-subdifferential-root-absolute} that $0$ is the
unique global minimizer of $f$.
\hfill
$\blacksquare$

If $\map{f}{\real^d}{\real}$ is locally Lipschitz at $x \in \real^d$,
then the proximal subdifferential of $f$ at $x$ is bounded.  In
general, the relationship between the generalized gradient and the
proximal subdifferential of a function $f$ that is locally Lipschitz at $x \in
\real^d$ is expressed by
\begin{align*}
  \partial f (x) = \co \setdef{\lim_{n \to \infty} \zeta_n \in
    \real^d}{\zeta_n \in \partial_P f(x_n) \; \text{and} \; \lim_{n
      \to \infty} x_n = x} .
\end{align*}

\bigskip
\noindent \textbf{\textit{Computing the proximal subdifferential}}
\medskip

As with the generalized gradient, computation of the proximal
subdifferential of a lower semicontinuous function is often far from
straightforward. Here we provide some useful results following the
exposition in~\cite{FHC-YSL-RJS-PRW:98}.

\textbf{Dilation rule.} If $\map{f}{\real^d}{\real}$ is lower
semicontinuous at $x \in \real^d$ and $s\in \realpositive$, then the
function $ s f$ is lower semicontinuous at $x$, and
\begin{align}\label{eq:differential-dilation}
  \partial_P (s f)(x) = s \; \partial_P f(x).
\end{align}

\textbf{Sum rule.} If $\map{f_1,f_2}{\real^d}{\real}$ are lower
semicontinuous at $x \in \real^d$, then the function $ f_1 + f_2$ is
lower semicontinuous at $x$, and
\begin{align}\label{eq:differential-sum}
  \partial_P f_1(x) + \partial_P f_2(x) \subseteq \partial_P \big( f_1
  + f_2\big)(x) .
\end{align}
Moreover, if either $f_1$ or $f_2$ are twice continuously
differentiable, then equality holds.

\textbf{Fuzzy sum rule.} If $\map{f_1,f_2}{\real^d}{\real}$ is lower
semicontinuous at $x \in \real^d$, $\zeta \in \partial_P \big( f_1 +
f_2\big)(x)$, and $\eps >0$, then there exist $x_1,x_2 \in
\oball{\eps}{x}$ with $|f_i(x_i) - f_i(x)| < \eps$, $i \in \{ 1, 2
\}$, such that
\begin{align}\label{eq:fuzzy-sum}
  \zeta \in \partial_P f_1(x_1) + \partial_P f_2(x_2) + \eps
  \oball{1}{0}.
\end{align}

\textbf{Chain rule.} Assume that either $\map{h}{\real^d}{\real^m}$
linear and $\map{g}{\real^m}{\real}$ lower semicontinuous at $h(x) \in
\real^m$, or $\map{h}{\real^d}{\real^m}$ locally Lipschitz at $x \in
\real^d$ and $\map{g}{\real^m}{\real}$ locally Lipschitz at $h(x) \in
\real^m$. If $\zeta \in \partial_P (g \circ h) (x)$ and $\eps >0$,
then there exist $\tilde x \in \real^d$, $\tilde y \in \real^m$, and
$\gamma \in \partial_P g(\tilde y)$ with $ \max \{\TwoNorm{\tilde x -
  x} , \TwoNorm{\tilde y - h(x)} \}<\eps$ such that $\TwoNorm{h(\tilde
  x) - h(x)} < \eps$ and
\begin{align}\label{eq:chain-rule-proximal}
  \zeta \in \partial_P (\gamma^T h) (\tilde x) + \eps \oball{1}{0} ,
\end{align}
where $\map{\gamma^T h}{\real^d}{\real}$ is defined by $ (\gamma^T h)
(x) = \gamma^T h(x)$.

The statement of the chain rule~\eqref{eq:chain-rule-proximal} shows a
characteristic feature of the proximal subdifferential.  That is,
rather than at an specific point of interest, the proximal
subdifferential can only be expressed with objects evaluated at
neighboring points.

Computation of the proximal subdifferential of a twice continuously
differentiable function is particularly simple.  If
$\map{f}{\real^d}{\real}$ is twice continuously differentiable on the
open set $U \subseteq \real^d$, then, for all $x \in U$,
\begin{align}
  \partial_P f (x) = \{ \nabla f (x) \} .
\end{align}
This simplicity also works for continuously differentiable convex
functions, as the following
result~\cite{FHC-YSL-RJS-PRW:98,RTR:70,GVS:01} states.  Note that
every convex function on $\real^d$ is locally Lipschitz, and hence
continuous.

\begin{proposition}\label{prop:proximal-convex}
  If $\map{f}{\real^d}{\real}$ is convex, then the following
  statements hold:
  \begin{enumerate}
  \item For $x \in \real^d$, $\zeta \in \partial_P f(x)$ if and only
    if, for all $y \in \real^d$, $ f(y) \ge f(x) + \zeta (y-x)$.
    
  \item The set-valued map $\map{\partial_P
      f}{\real^d}{\parts{\real^d}}$, $x\mapsto \partial_P f (x)$,
    takes nonempty, compact, and convex values, and is upper
    semicontinuous and locally bounded.
    
  \item If, in addition, $f$ is continuously differentiable, then, for
    all $x \in \real^d$, $\partial_P f (x) = \{ \nabla f (x) \}$.
  \end{enumerate}
\end{proposition}

\subsubsection*{Example~\ref{ex:abs-prox} revisited: Proximal
  subdifferential of the absolute value function}

The computation~\eqref{eq:comp-proximal-subdifferential-a} of the
proximal subdifferential of the function $f (x) = |x|$ can be
performed using Proposition~\ref{prop:proximal-convex} by observing
that $f$ is convex.  By
Proposition~\ref{prop:proximal-convex}\textit{(i)}, for $x=0$, $\zeta
\in
\partial_P f(0)$ if and only if, for all $y \in \real$,
\begin{align}\label{eq:prox-subdifferential-aux}
  |y| \ge |0| + \zeta (y-0) = \zeta y.
\end{align}
Hence, for $y \ge 0$,~\eqref{eq:prox-subdifferential-aux} implies
$\zeta \in (- \infty,1]$, while, for $y \le
0$,~\eqref{eq:prox-subdifferential-aux} implies that $\zeta \in
[-1,+\infty)$.  Therefore, $\zeta \in (- \infty,1] \cap [-1,+\infty)
$, and $\partial_P f (0) = [-1,1]$. By
Proposition~\ref{prop:proximal-convex}\textit{(iii)}, for $x > 0$,
$\partial_P f (x) = \{1\}$, and, for $x < 0$, $\partial_P f (x) =
\{-1\}$.
\hfill $\blacksquare$

Regarding critical points, if $x$ is a local minimizer of the lower
semicontinuous function $\map{f}{\real^d}{\real}$, then $0 \in
\partial_P f(x)$. Conversely, if $f$ is convex, and $0 \in \partial_P
f(x)$, then $x$ is a global minimizer of $f$.  For the study of maximizers,
instead of lower semicontinuous function, convex function, and
proximal subdifferential, the relevant notions are upper
semicontinuous function, concave function, and proximal
superdifferential, respectively~\cite{FHC-YSL-RJS-PRW:98}.

\bigskip
\noindent \textbf{\textit{Gradient differential inclusion of a convex
    function}} \medskip

It is not always possible to associate a nonsmooth gradient flow with
a lower semicontinuous function because the proximal subdifferential
might be empty almost everywhere.  However, following
Proposition~\ref{prop:proximal-convex}\textit{(ii)}, we can associate
a nonsmooth gradient flow with each convex function, as we briefly
discuss next~\cite{JPA-IE:84}.

Let $\map{f}{\real^d}{\real}$ be convex, and consider the gradient
differential inclusion
\begin{align}\label{eq:gradient-inclusion}
  \dot x (t) \in - \partial_P f (x(t)) .
\end{align}
Using the properties of the proximal subdifferential given by
Proposition~\ref{prop:proximal-convex}\textit{(ii)}, the existence of
a Caratheodory solution of~\eqref{eq:gradient-inclusion} starting from
every initial condition is guaranteed by
Proposition~\ref{prop:existence-solution}. Moreover, the uniqueness of
Caratheodory solutions can be established as follows.  Let $x,y \in
\real^d$, and take $\zeta_1 \in - \partial_P f(x)$ and $\zeta_2 \in -
\partial_P f(y)$. Using
Proposition~\ref{prop:proximal-convex}\textit{(i)}, we have
\begin{align*}
  f(y) \ge f(x) - \zeta_1 (y-x) , \qquad f(x) \ge f(y) - \zeta_2 (x-y)
  .
\end{align*}
From here, we deduce $- \zeta_1 (y-x) \le f(y) - f(x) \le - \zeta_2
(y-x)$, and therefore $(\zeta_2 - \zeta_1) (y-x) \le 0$, which, in
particular, implies that the set-valued map $x \mapsto -\partial_P
f(x)$ satisfies the one-sided Lipschitz
condition~\eqref{eq:Lipschitz-like}.
Proposition~\ref{prop:unique-sol-inclusion} then guarantees that there
exists a unique Caratheodory solution of~\eqref{eq:gradient-inclusion}
starting from every initial condition.

\section*{Nonsmooth Stability Analysis}

In this section, we study the stability of discontinuous dynamical
systems.  Unless explicitly mentioned otherwise, the stability notions
employed here correspond to the usual ones for differential equations,
see, for instance~\cite{HKK:02}.  The presentation focuses on the
time-invariant differential inclusion
\begin{align}\label{eq:auto-diff-inclusion}
  \dot x (t) \in \setvaluedmap (x(t)) ,
\end{align}
where $\map{\setvaluedmap}{\real^d}{\parts{\real^d}}$.  Throughout the
section, we assume that the set-valued map~$\setvaluedmap$ satisfies
the hypotheses of Proposition~\ref{prop:existence-solution}, so that
the existence of solutions of~\eqref{eq:auto-diff-inclusion} is
guaranteed.  The scenario~\eqref{eq:auto-diff-inclusion} has a direct
application to discontinuous differential equations and control
systems. For instance, the results presented here apply to
Caratheodory solutions by taking a singleton-valued map, to Filippov
solutions by taking $\setvaluedmap = F[X]$, and to control systems by
taking $\setvaluedmap = G[X]$.

Before proceeding with our exposition, we recall that solutions of a
discontinuous system are not necessarily unique.  Therefore, when
considering a property such as Lyapunov stability, we must specify
whether attention is being paid to a particular solution starting from
an initial condition (``weak'') or to all the solutions starting from
an initial condition (``strong'').  As an example, the set $M
\subseteq \real^d$ is \emph{weakly invariant}\index{set!weakly
  invariant} for~\eqref{eq:auto-diff-inclusion} if, for each $x_0 \in
M$, $M$ contains at least one maximal solution
of~\eqref{eq:auto-diff-inclusion} with initial condition $x_0$.
Similarly, $M \subseteq \real^d$ is \emph{strongly
  invariant}\index{set!strongly invariant}
for~\eqref{eq:auto-diff-inclusion} if, for each $x_0 \in M$, $M$
contains every maximal solution of~\eqref{eq:auto-diff-inclusion} with
initial condition $x_0$. Rather than reintroducing weak and strong
versions of standard stability notions, we rely on the guidance
provided above and the context for the specific meaning in each case.
Detailed definitions are given in~\cite{AFF:88,AB-LR:05}.

The following discussion requires the notion of a limit point of a
solution of a differential inclusion.  The point $x_* \in \real^d$ is
a \emph{limit point} of a solution $t \mapsto x(t)$
of~\eqref{eq:auto-diff-inclusion} if there exists a sequence
$\{t_n\}_{n \in \natural}$ such that $x (t_n) \rightarrow x_*$ as
$n\rightarrow \infty$. We denote by $\Omega (x)$ the set of limit
points of~$t \mapsto x(t)$. Under the hypotheses of
Proposition~\ref{prop:existence-solution}, $\Omega (x)$ is a weakly
invariant set. Moreover, if the solution $t \mapsto x(t)$ lies in a
bounded, open, and connected set, then $\Omega (x)$ is nonempty,
bounded, connected, and $x (t) \rightarrow \Omega (x)$ as $t
\rightarrow \infty$, see~\cite{AFF:88}.

\index{limit point}%
\sindex{LimitSet}{{$\Omega (x)$}}{Set of limit points of a curve $ t
  \mapsto x(t)$}%

\subsection*{Stability analysis by means of the generalized gradient
  of a nonsmooth Lyapunov
  function} 

In this section, we discuss nonsmooth stability analysis results that
use locally Lipschitz functions and generalized gradients.  We present
results taken from various sources in the literature.  This discussion
is not intended to be a comprehensive account of such a vast topic,
but rather serves as a motivation for further exploration.  The
books~\cite{AFF:88,AB-LR:05} and journal
articles~\cite{DS-BP:94,EPR:98,AB-FC:99,AB-FC:06} provide additional
information.

\bigskip
\noindent \textbf{\textit{{Lie derivative and
      monotonicity}}}
\medskip

A common theme in stability analysis is establishing the monotonic
evolution of a candidate Lyapunov function along the trajectories of
the system. Mathematically, the evolution of a function along
trajectories is captured by the notion of Lie derivative. Our first
task is to generalize this notion to the setting of differential
inclusions following~\cite{AB-FC:99}, see also~\cite{DS-BP:94,EPR:98}.

Given a locally Lipschitz function $\map{f}{\real^d}{\real}$ and a
set-valued map $\map{\setvaluedmap}{\real^d}{\parts{\real^d}}$, the
\emph{set-valued Lie derivative}
$\map{\setLiederG{\setvaluedmap}{f}}{\real^d}{\parts{\real}}$ of $f$
with respect to $\setvaluedmap$ at $x$ is defined as
\begin{align}\label{eq:LieDerG}
  \setLiederG{\setvaluedmap}{f} (x) = \setdef{a \in \real}{\text{there
      exists} \; v \in \setvaluedmap (x) \; \text{such that} \;
    \zeta^T v = a \; \text{for all} \; \zeta \in \partial f(x)} .
\end{align}
\index{Lie derivative!set-valued}
\sindex{SetLieDer}{$\setLiederG{\setvaluedmap}{f}$}{Set-valued Lie
  derivative of the locally Lipschitz function
  $\map{f}{\real^d}{\real}$ with respect to the set-valued map
  $\map{\setvaluedmap}{\real^d}{\parts{\real^d}}$}%
If $\setvaluedmap$ takes convex and compact values, then, for each $x
\in \real^d$, $\setLiederG{\setvaluedmap}{f}(x)$ is a closed and
bounded interval in $\real$, possibly empty.  If $f$ is continuously
differentiable at $x$, then
\begin{align*}
  \setLiederG{\setvaluedmap}{f}(x) = \setdef{(\nabla f (x))^T v}{v \in
    \setvaluedmap (x)}.
\end{align*}
The usefulness of the set-valued Lie derivative stems from the fact
that it allows us to study how the function $f$ evolves along the
solutions of a differential inclusion without having to explicitly
obtain the solutions.  Specifically, we have the following
result, see~\cite{AB-FC:99}.

\begin{proposition}\label{prop:fundamental}
  Let $x:[0,t_1] \rightarrow \real^d$ be a solution of the
  differential inclusion~\eqref{eq:auto-diff-inclusion}, and let
  $\map{f}{\real^d}{\real}$ be locally Lipschitz and regular. Then,
  the following statements hold:
  \begin{enumerate}
  \item The composition $ t \mapsto f (x(t))$ is differentiable at
    almost all $t \in [0,t_1]$.
  \item The derivative of $ t \mapsto f (x(t))$ satisfies
    \begin{align}\label{eq:evolution}
      \frac{d}{dt} \big( f (x(t)) \big) \in
      \setLiederG{\setvaluedmap}{f} (x(t)) \quad \text{for almost
        every $t \in [0,t_1]$.}
    \end{align}
  \end{enumerate}
\end{proposition}

A similar result can be established for a larger class of functions
when the set-valued map $\setvaluedmap$ is
singleton-valued~\cite{MV:89,AB-FC:06}.

Given the discontinuous vector field $\map{X}{\real^d}{\real^d}$,
consider the Filippov solutions of~\eqref{eq:ODE-auto}. In this case,
with a slight abuse of notation, we denote $\setLiederG{X}{f} =
\setLiederG{F[X]}{f}$. Note that if $X$ is continuous at $x$, then
$F[X] (x) = \{X(x)\}$, and therefore, $\setLiederG{X}{f}(x)$
corresponds to the singleton $\{ \Lc_{X} f (x) \}$, whose sole element
is the usual Lie derivative of $f$ in the direction of $X$ at $x$.
\sindex{SetLieDerII}{$\setLiederG{X}{f}$}{Set-valued Lie derivative of
  the locally Lipschitz function $\map{f}{\real^d}{\real}$ with
  respect to the Filippov set-valued map
  $\map{F[X]}{\real^d}{\parts{\real^d}}$}

\subsubsection*{Example~\ref{ex:oscillator} revisited: Monotonicity in
  the nonsmooth harmonic
  oscillator}

For the nonsmooth harmonic oscillator in Example~\ref{ex:oscillator},
consider the locally Lipschitz and regular map
$\map{f}{\real^2}{\real}$, $f(x_1,x_2) = |x_1| + \frac{x_2^2}{2}$.
Recall that Figure~\ref{fig:oscillator}(b) shows the contour plot of
$f$. Let us determine how $f$ evolves along the Filippov solutions of
the dynamical system by looking at the set-valued Lie derivative.
First, we compute the generalized gradient of $f$. To do so, we
rewrite $f$ as $f(x_1,x_2) = \max \{ x_1, -x_1 \} + \frac{x_2^2}{2}$,
and apply Proposition~\ref{prop:gradients-max-min}\textit{(ii)} and
the sum rule to find
\begin{align*}
  \partial f (x_1,x_2) =
  \begin{cases}
    \{ (\sign (x_1), x_2) \}, & x_1 \neq 0  , \\
    [-1,1] \times \{ x_2 \}, & x_1 = 0 .
  \end{cases}
\end{align*}
With this information, we can compute the set-valued Lie derivative
$\map{\setLiederG{X}{f}}{\real^2}{\parts{\real}}$ as
\begin{align}\label{eq:lie-der-harmonic}
  \setLiederG{X}{f} (x_1,x_2) = 
  \begin{cases}
    \{ 0 \}, & x_1 \neq 0  , \\
    \emptyset , & x_1 = 0 \; \text{and} \;
    x_2 \neq 0 ,\\
    \{ 0\} , & x_1 = 0 \; \text{and} \; x_2 = 0 .
  \end{cases}
\end{align}
From~\eqref{eq:evolution} and~\eqref{eq:lie-der-harmonic} we conclude
that $f$ is constant along the Filippov solutions of the discontinuous
dynamical system.  Indeed, the level sets of $f$ are exactly the
curves described by the solutions of the system in
Figure~\ref{fig:oscillator}(b).
\hfill
$\blacksquare$

\bigskip
\noindent \textbf{\textit{Stability results}}
\medskip

The above discussion on monotonicity is the stepping stone to
stability results using locally Lipschitz functions and generalized
gradient information.  Proposition~\ref{prop:fundamental} provides a
criterion for determining the monotonic behavior of a locally
Lipschitz function along the solutions of discontinuous dynamics. This
result, together with the appropriate positive definite assumptions on
the candidate Lyapunov function, allows us to synthesize checkable
stability tests.  We start by formulating~\cite{AB-FC:99} the natural
extension of Lyapunov's stability theorem for ordinary differential
equations.  Recall that, at each $x \in \real^d$, the Lie derivative
$\setLiederG{\setvaluedmap}{f}(x)$ is a set contained in $\real$. For
the empty set, we adopt the convention $\max \emptyset = -\infty$.

\begin{theorem}\label{th:Lyapunov}
  Let $\map{\setvaluedmap}{\real^d}{\parts{\real^d}}$ be a set-valued
  map satisfying the hypotheses of
  Proposition~\ref{prop:existence-solution}, let $x_e$ be an
  equilibrium of the differential
  inclusion~\eqref{eq:auto-diff-inclusion}, and let $\domain \subseteq
  \real^d$ be an open and connected set with $x_e \in \domain$.
  Furthermore, let $\map{f}{\real^d}{\real}$ be such that the
  following conditions hold:
  \begin{enumerate}
  \item $f$ is locally Lipschitz and regular on $\domain$.
  \item $f(x_e)=0$, and $f(x)>0$ for $x \in \domain \! \setminus \!
    \{x_e\}$.
  \item $\max \setLiederG{\setvaluedmap}{f}(x) \le 0$ for each $x \in
    \domain$.
  \end{enumerate}
  Then, $x_e$ is a strongly stable equilibrium
  of~\eqref{eq:auto-diff-inclusion}. In addition, if \textit{(iii)}
  above is replaced by
  \begin{enumerate}
  \item[(iii)'] $\max \setLiederG{\setvaluedmap}{f}(x) < 0$ for each
    $x \in \domain \! \setminus \! \{x_e \}$,
  \end{enumerate}
  then $x_e$ is a strongly asymptotically stable equilibrium
  of~\eqref{eq:auto-diff-inclusion}.
\end{theorem}

Let us apply this result to the nonsmooth harmonic oscillator.

\subsubsection*{Example~\ref{ex:oscillator} revisited: Stability
  analysis of the nonsmooth harmonic
  oscillator}

The function $(x_1,x_2) \rightarrow |x_1| + \frac{x_2^2}{2}$ satisfies
hypotheses \textit{(i)-(iii)} of Theorem~\ref{th:Lyapunov} on $\domain
= \real^d$.  Therefore, we conclude that $0$ is a strongly stable
equilibrium. The phase portrait in Figure~\ref{fig:oscillator}(a)
indicates that $0$ is not strongly asymptotically stable. Using
Theorem~\ref{th:Lyapunov}, it can be shown that the nonsmooth harmonic
oscillator under dissipation, with vector field $(x_1,x_2) \mapsto
(x_2,-\sign(x_1) - k \sign (x_2))$, where $k>0$, has $0$ as a strongly
asymptotically stable equilibrium.  The phase portrait of this system
is plotted in Figure~\ref{fig:oscillator-dissipation}(a), while
several Filippov solutions are plotted in
Figure~\ref{fig:oscillator-dissipation}(b).
\hfill $\blacksquare$

Another useful result in the theory of differential equations is the
invariance principle~\cite{HKK:02}. In many situations, this principle
allows us to determine the asymptotic convergence properties of the
solutions of a differential equation. Here, we build on the above
discussion to present a generalization to differential
inclusions~\eqref{eq:auto-diff-inclusion} and nonsmooth Lyapunov
functions.  This principle is thus applicable to discontinuous
differential equations.  The formulation is taken
from~\cite{AB-FC:99}, which slightly generalizes the presentation
in~\cite{DS-BP:94}.

\begin{theorem}\label{th:LaSalle}
  Let $\map{\setvaluedmap}{\real^d}{\parts{\real^d}}$ be a set-valued
  map satisfying the hypotheses of
  Proposition~\ref{prop:existence-solution}, and let
  $\map{f}{\real^d}{\real}$ be a locally Lipschitz and regular
  function.  Let $S \subset \real^d$ be compact and strongly invariant
  for~\eqref{eq:auto-diff-inclusion}, and assume that $\max
  \setLiederG{\setvaluedmap}{f}(y) \le 0$ for each $y \in S$.  Then,
  all solutions $\map{x}{[0,\infty)}{\real^d}$
  of~\eqref{eq:auto-diff-inclusion} starting at $S$ converge to the
  largest weakly invariant set $M$ contained in
  \begin{align*}
    S \cap \overline{\setdef{y \in \real^d}{0 \in
        \setLiederG{\setvaluedmap}{f}(y)}} .
  \end{align*}
  Moreover, if the set $M$ consists of a finite number of points, then
  the limit of each solution starting in $S$ exists and is an element
  of $M$.
\end{theorem}

We now apply Theorem~\ref{th:LaSalle} to nonsmooth gradient flows.

\bigskip
\noindent \textbf{\textit{Stability of nonsmooth gradient flows}}
\medskip

Consider the nonsmooth gradient flow~\eqref{eq:natural-gradient} of a
locally Lipschitz function $\map{f}{\real^d}{\real}$.  Assume further
that $f$ is regular.  Let us examine how $f$ evolves along the
solutions of the flow using the set-valued Lie derivative. Given $x
\in \real^d$, let $a \in \setLiederG{-\LN (\partial f)}{f}(x)$. By
definition, there exists $v \in F[-\LN (\partial f)](x)$ such that
\begin{align}\label{eq:aux-equality}
  a = \zeta^T v \quad \text{for all} \; \zeta \in \partial f (x) .
\end{align}
Recall from~\eqref{eq:magic-Filippov} that $ F[-\LN (\partial f)](x) =
-\partial f (x)$.  Since~\eqref{eq:aux-equality} holds for all
elements in the generalized gradient of $f$ at $x$, we can choose in
particular $\zeta = -v \in \partial f(x)$. Therefore,
\begin{align}\label{eq:lie-der-gradient}
  a = (-v)^T v = - \TwoNorm{v}^2 \le 0 .
\end{align}
From~\eqref{eq:lie-der-gradient}, we conclude that all of the elements
of $\setLiederG{-\LN (\partial f)}{f}$ belong to $\realnonpositive$,
and therefore, from~\eqref{eq:evolution}, $f$ monotonically decreases
along the solutions of its nonsmooth gradient
flow~\eqref{eq:natural-gradient}.  Moreover, we deduce that $0 \in
\setLiederG{-\LN (\partial f)}{f}(x)$ if and only if $0 \in
\partial f(x)$, that is, if $x$ is a critical point of $f$.  The
application of the Lyapunov stability theorem and the invariance
principle in theorems~\ref{th:Lyapunov} and~\ref{th:LaSalle},
respectively, gives rise to the following nonsmooth counterpart of the
classical smooth result for a gradient flow~\cite{WMH-SS:74}.

\begin{proposition}\label{prop:stability-nonsmooth-gradient-flows}
  Let $\map{f}{\real^d}{\real}$ be locally Lipschitz and regular.
  Then, the strict minimizers of $f$ are strongly stable equilibria of
  the nonsmooth gradient flow~\eqref{eq:natural-gradient} of $f$.
  Furthermore, if the level sets of $f$ are bounded, then the
  solutions of the nonsmooth gradient flow asymptotically converge to
  the set of critical points of~$f$.
\end{proposition}

The following example illustrates the above discussion.

\subsubsection*{Example~\ref{ex:sm} revisited: Stability analysis of
  the nonsmooth gradient flow of
  $-\sm_Q$}

\begin{figure}[htbp]
  \centering
  \includegraphics[width=.2\linewidth]{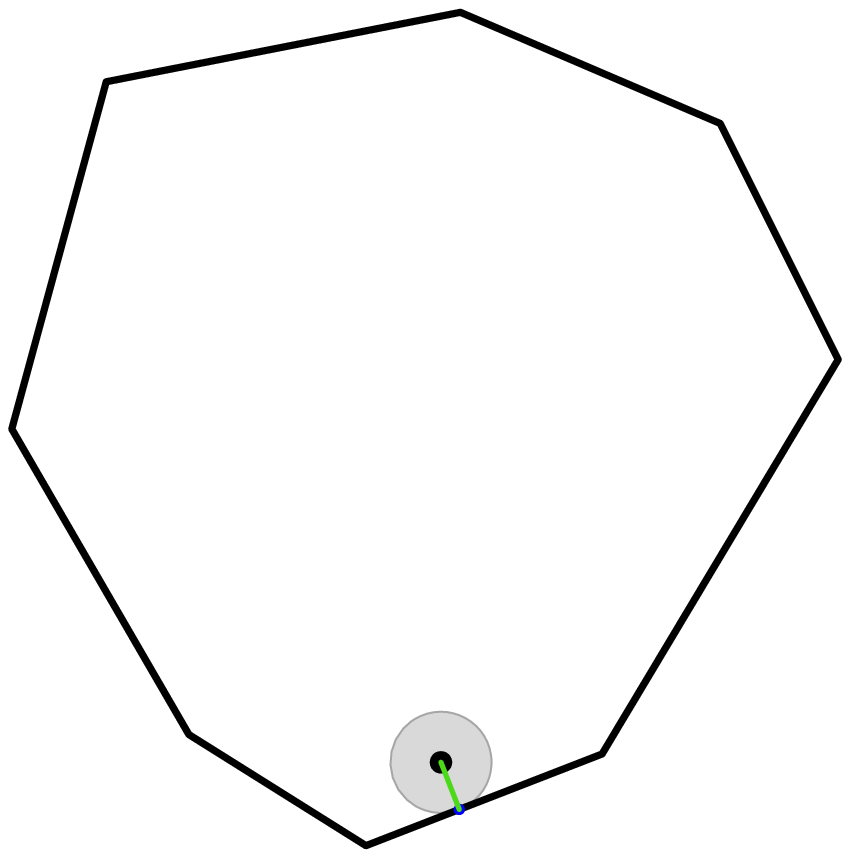}%
  \includegraphics[width=.2\linewidth]{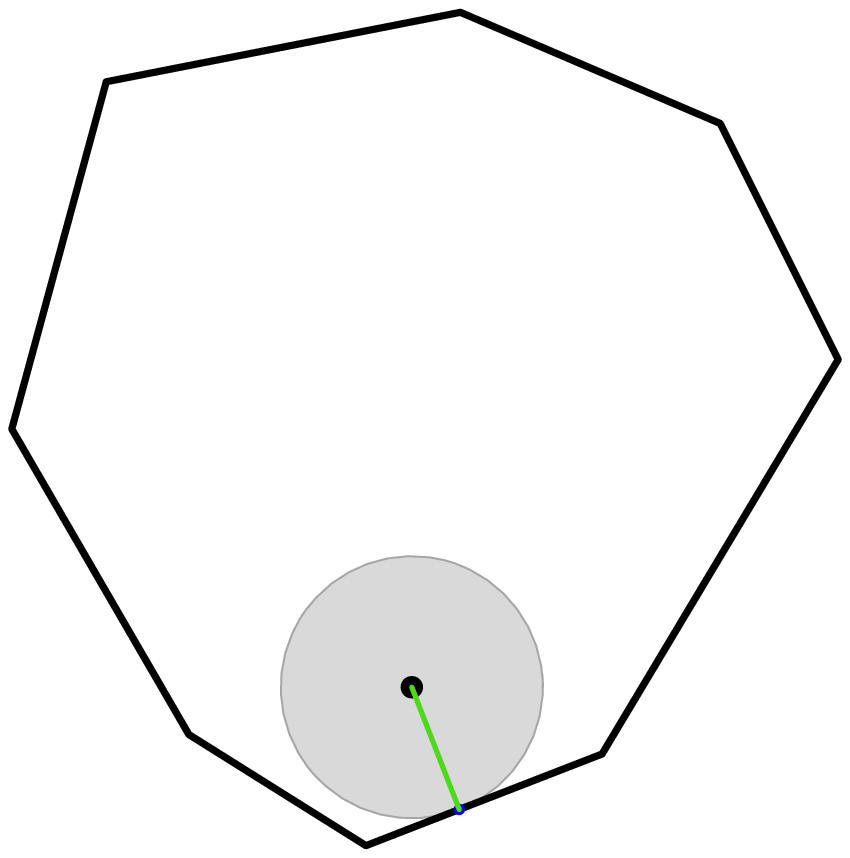}%
  \includegraphics[width=.2\linewidth]{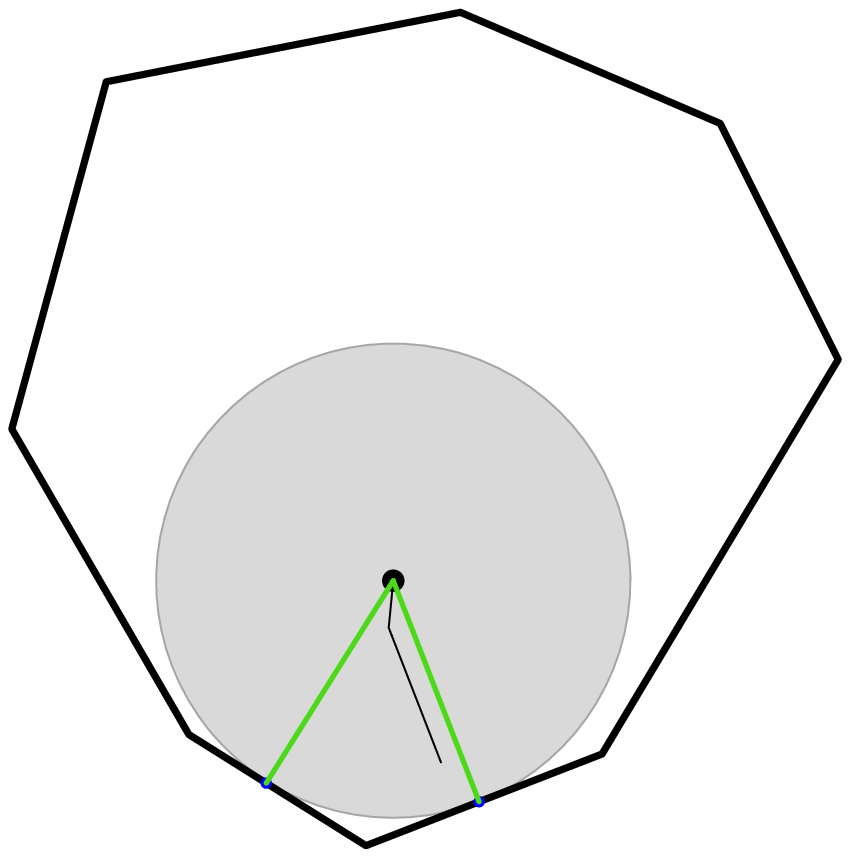}%
  \includegraphics[width=.2\linewidth]{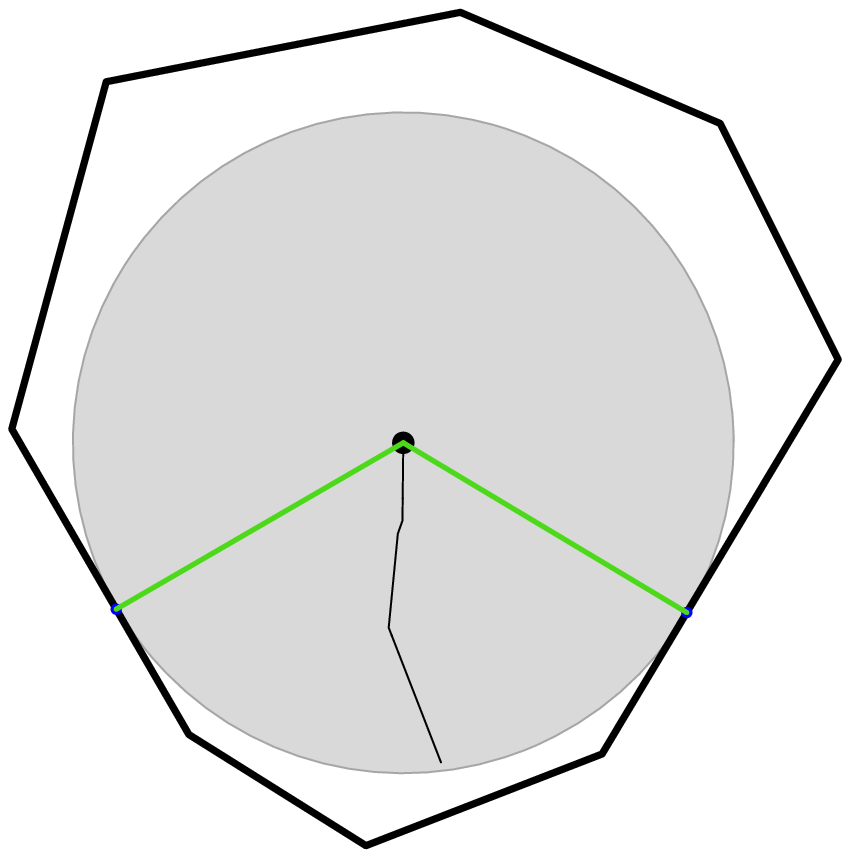}%
  \includegraphics[width=.2\linewidth]{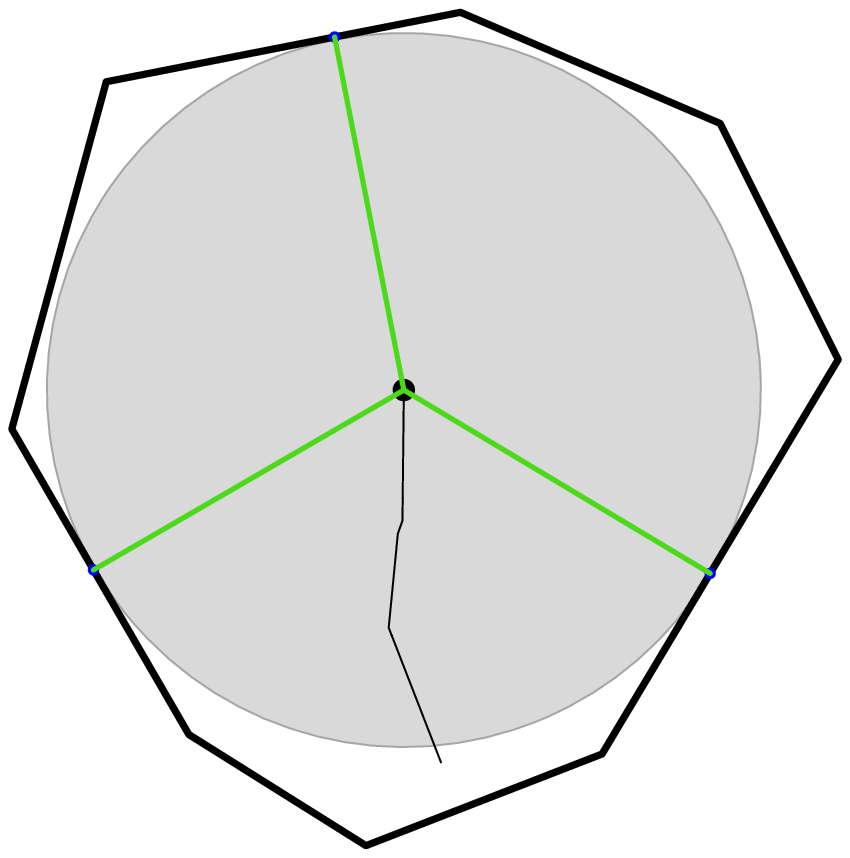} 
  \caption{From left to right, evolution of the nonsmooth gradient
    flow of the function $-\sm_Q$ in a convex polygon. At each
    snapshot, the value of $\sm_Q$ is the radius of the largest disk
    (plotted in gray) contained in the polygon with center at the
    current location. The flow converges in finite time to the
    incenter set, which, for this polygon, is a singleton whose only
    element is the center of the disk in the rightmost snapshot.}
  \label{fig:sm-flow}
\end{figure}

Consider the nonsmooth gradient flow of $-\sm_Q$, which is the
minimum-distance-to-polygonal-boundary function introduced in
Example~\ref{ex:sm}.  Proposition~\ref{prop:uniqueness-Filippov-II}
guarantees Uniqueness of solutions for this flow. Regarding
convergence, Proposition~\ref{prop:stability-nonsmooth-gradient-flows}
guarantees that solutions converge asymptotically to the incenter set.
Indeed, the incenter set is attained in finite time, and hence each
solution converges to a point of the incenter set~\cite{JC-FB:02m}.
The nonsmooth gradient flow can be interpreted as a sphere-packing
algorithm in the sense that, starting from an arbitrary initial
point, the flow monotonically maximizes the radius of the largest disk
contained in the polygon (that is, $\sm_Q$!) until it reaches an
incenter point.  This fact is illustrated in Figure~\ref{fig:sm-flow}.

\begin{figure}[htbp]
  \centering%
  \subfigure[]{
    \includegraphics[width=.3\linewidth]{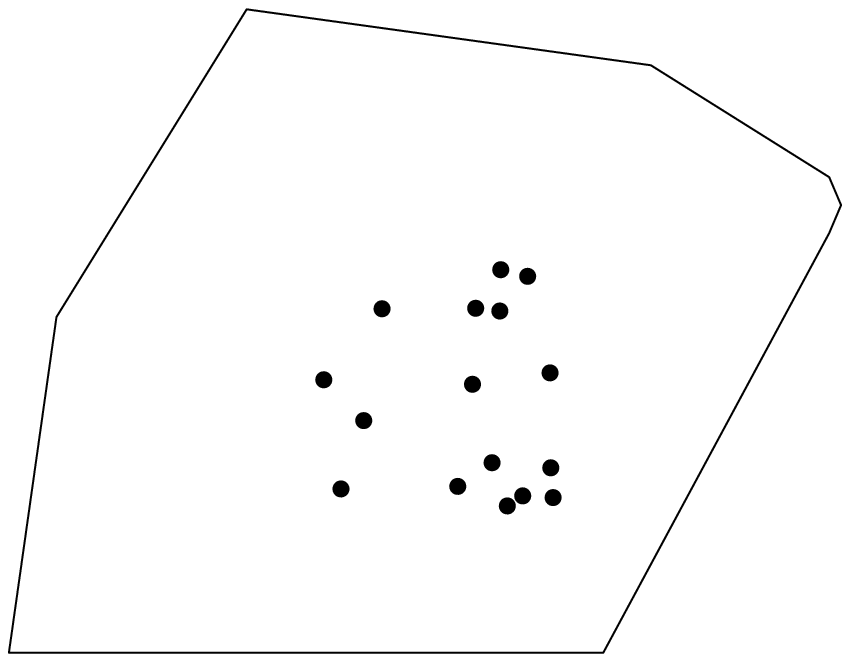}}%
  \subfigure[]{
    \includegraphics[width=.3\linewidth]{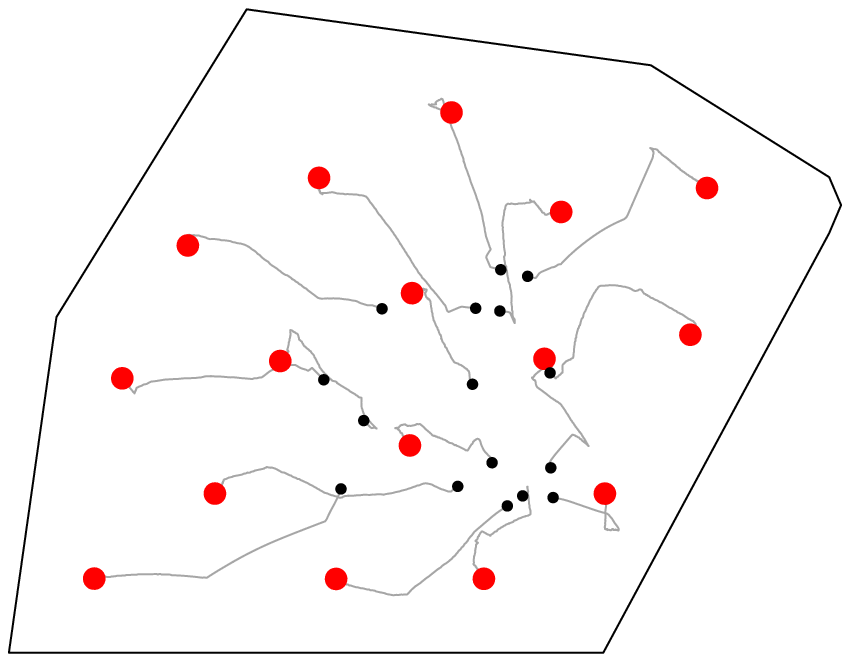}}%
  \subfigure[]{\includegraphics[width=.3\linewidth]{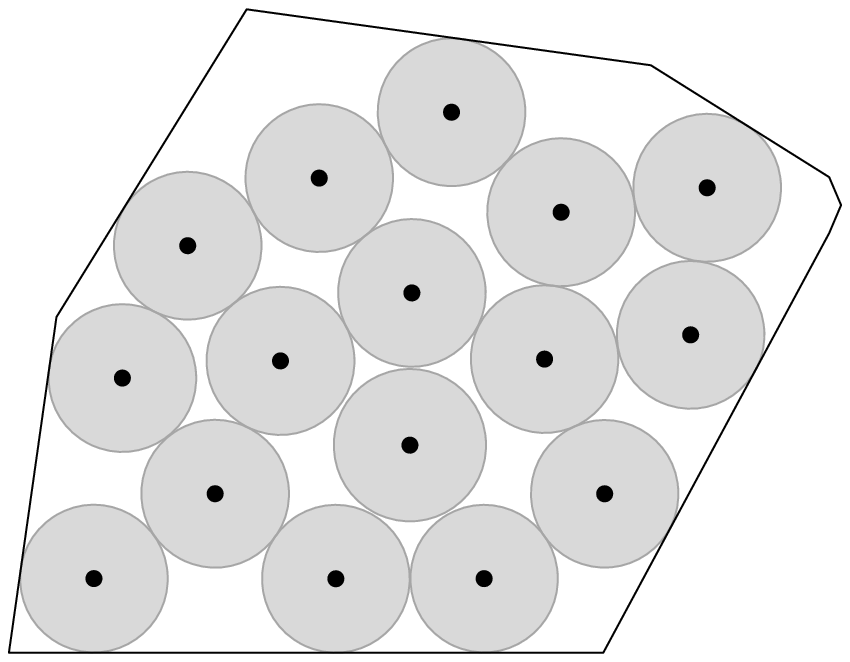}}
  \caption{The move-away-from-nearest-neighbor interaction law for
    solving a sphere-packing problem within the polygon $Q$.  (a) is
    the initial configuration, (b) is the evolution, and (c) is the
    final configuration of a Filippov solution. In (c), the minimum
    radius of the shaded spheres corresponds to the value of the
    locally Lipschitz function~$\HHSP$, defined in~\eqref{eq:HHSP}.
    The Filippov solutions of the move-away-from-nearest-neighbor
    dynamical system monotonically increase the value of $\HHSP$.  In
    (c), every node is at equilibrium since the infinitesimal motion
    of a node in any direction would place it closer to at least
    another node. This discontinuous dynamical system is an example of
    how simple local interactions can achieve a global objective. }
  \label{fig:multiple-sphere-packing}
\end{figure}

What if, instead, we want to pack an arbitrary number $n$ of spheres
within the polygon?  It turns out that the
move-away-from-nearest-neighbor interaction law is a discontinuous
dynamical system that solves this problem, where the solutions are
understood in the Filippov sense.
Figure~\ref{fig:multiple-sphere-packing} illustrates the evolution of
this dynamical system.  The stability properties of this law can be
determined using the Lyapunov function $\map{\HHSP}{Q^n}{\real}$
defined by
\begin{align}\label{eq:HHSP}
  \HHSP (p_1,\dots,p_n) = \min \setdef{\tfrac{1}{2} \TwoNorm{p_i-p_j}
    , \D(p_i,e)}{i \neq j\in\until{n} , \, e \; \text{edge of $Q$}}.
\end{align}
The value of $\HHSP$ corresponds to the largest radius such that
spheres with centers $p_1,\dots,p_n$ and radius $\HHSP
(p_1,\dots,p_n)$ fit within the environment and do not intersect each
other (except at most at the boundary).  The function $\HHSP$ is
locally Lipschitz, but not concave, on~$Q$, and its evolution is
monotonically nondecreasing along the move-away-from-nearest-neighbor
dynamical system.  The reader is referred to~\cite{JC-FB:02m} for
details on the relationship between $\HHSP$ and the minimum distance
function $\sm_Q$, as well as additional discontinuous dynamical
systems that solve this sphere-packing problem and other exciting
geometric optimization problems.
\hfill $\blacksquare$

\bigskip
\noindent 
\textbf{\textit{Finite-time-convergent nonsmooth gradient flows}}
\medskip

Convergence of the solutions of a dynamical system in finite time is a
desirable property in various settings. For instance, a motion
planning algorithm is effective if it can steer a robot from the
initial to the final configuration in finite time rather than
asymptotically.  Another example is given by a robotic network trying
to agree on the exact value of a quantity sensed by each agent.
Reaching agreement in finite time allows the network to use precise
information in the completion of other tasks.

Broadly applicable results on finite-time convergence for
discontinuous dynamical systems can be found
in~\cite{BP-SSS:87,JC:06-auto,SA-HA-AC:06,AC:07}.  Here, we briefly
discuss the finite-convergence properties of a class of nonsmooth
gradient flows.  Let $\map{f}{\real^d}{\real}$ be a continuously
differentiable function and assume that all of the level sets of $f$
are bounded.  It is well known~\cite{WMH-SS:74} that all solutions of
the gradient flow $\dot x (t) = - \nabla f (x(t))$ converge
asymptotically toward the set of critical points of $f$, although the
convergence does not occur in finite time.  Here, we slightly modify
the gradient flow into two different nonsmooth flows that achieve
finite-time convergence.  Consider the discontinuous differential
equations
\begin{align}
  \dot{x} (t) &= - \frac{\nabla f(x(t))}{\TwoNorm{\nabla f (x(t))}} ,
  \label{eq:gradient-flow-norm} \\
  \dot{x} (t) &= - \sign(\nabla f (x (t))) ,
  \label{eq:gradient-flow-sign}
\end{align}
where $\TwoNorm{\cdot}$ denotes the Euclidean distance, and $\sign (x)
= (\sign(x_1),\dots,\sign(x_d)) \in \real^d$.  The nonsmooth vector
field~\eqref{eq:gradient-flow-norm} always points in the direction of
the gradient with unit speed.  Alternatively, the nonsmooth vector
field~\eqref{eq:gradient-flow-sign} specifies the direction of motion
by means of a binary quantization of the direction of the gradient.
The following result~\cite{JC:06-auto} characterizes the properties of
the Filippov solutions of~\eqref{eq:gradient-flow-norm}
and~\eqref{eq:gradient-flow-sign}.

\begin{proposition}\label{prop:finite-convergent-gradient-flows}
  Let $\map{f}{\real^d}{\real}$ be a twice continuously differentiable
  function.  Let $S \subset \real^d$ be compact and strongly invariant
  for~\eqref{eq:gradient-flow-norm} (resp.,
  for~\eqref{eq:gradient-flow-sign}).  If the Hessian of~$f$ is
  positive definite at every critical point of $f$ in $S$, then every
  Filippov solution of~\eqref{eq:gradient-flow-norm}
  (resp.~\eqref{eq:gradient-flow-sign}) starting from $S$ converges in
  finite time to a minimizer of~$f$.
\end{proposition}

The proof of Proposition~\ref{prop:finite-convergent-gradient-flows}
builds on the stability tools presented in this section. Specifically,
the invariance principle can be used to establish convergence toward
the set of critical points of $f$.  Finite-time convergence can be
established by deriving bounds on the evolution of $f$ along the
solutions of the discontinuous dynamics using the set-valued Lie
derivative. This analysis also allows us to provide upper bounds on
the convergence time.  A more comprehensive exposition of results that
guarantee finite-time convergence of general discontinuous dynamics is
given in~\cite{JC:06-auto}.

\refstepcounter{example}%
\subsubsection*{Example~\theexample: Finite-time consensus}

The ability to reach consensus, or agreement, upon some (a priori
unknown) quantity is critical for multi-agent systems.  Network
coordination problems require that individual agents agree on the
identity of a leader, jointly synchronize their operation, decide
which specific pattern to form, balance the computational load, or
fuse consistently the information gathered on some spatial process.
Here, we briefly comment on two discontinuous algorithms that achieve
consensus in finite time, following~\cite{JC:06-auto}.
    
Consider a network of $n$ agents with states $p_1,\dots,p_n \in
\real$.  Let $G = (\until{n},E)$ be an undirected graph with $n$
vertices, describing the topology of the network.  Two agents $i$ and
$j$ \emph{agree} if and only if $p_i = p_j$.  The \emph{disagreement
  function} $\Phi_G : \real^n \rightarrow \realnonnegative$ quantifies
the group disagreement
\begin{align*}
  \Phi_G(p_1,\dots,p_n) = \frac{1}{2} \sum_{(i,j) \in E} (p_j -
  p_i)^2 .
\end{align*}
It is known~\cite{ROS-RMM:03c} that, if the graph is connected, the
gradient flow of $\Phi_G$ achieves consensus with an exponential rate
of convergence.  Actually, agents agree on the average value of their
initial states (this consensus is called \emph{average consensus}).
If $G$ is connected, the nonsmooth gradient
flow~\eqref{eq:gradient-flow-norm} of $\Phi_G$ achieves average
consensus in finite time, and the nonsmooth gradient
flow~\eqref{eq:gradient-flow-sign} of $\Phi_G$ achieves consensus on
the average of the maximum and the minimum of the initial states in
finite time, see~\cite{JC:06-auto}.
\hfill
$\blacksquare$

\subsection*{Stability analysis by means of the proximal
  subdifferential of a nonsmooth Lyapunov
  function}

This section presents stability tools for differential inclusions
using a lower semicontinuous function as a candidate Lyapunov
function.  We use the proximal subdifferential to study the monotonic
evolution of the candidate Lyapunov function along the solutions of
the differential inclusion.  As in the previous section, we present a
few representative and useful results. We refer the interested reader
to~\cite{FHC-YSL-RJS-PRW:98,EDS:99} for a more detailed exposition.

\bigskip
\noindent \textbf{\textit{Lie derivative and
  monotonicity}}
\medskip

Let $\domain \subseteq \real^d$ be an open and connected set.  A lower
semicontinuous function $\map{f}{\real^d}{\real}$ is \emph{weakly
  nonincreasing on~$\domain$}\index{function!weakly nonincreasing} for
a set-valued map $\map{\setvaluedmap}{\real^d}{\parts{\real^d}}$ if,
for all $y \in \domain$, there exists a solution
$\map{x}{[0,t_1]}{\real^d}$ of the differential
inclusion~\eqref{eq:auto-diff-inclusion} starting at $y$ and contained in
$\domain$ that satisfies
\begin{align*}
  f (x(t)) \le f (x (0)) = f (y) \quad \text{for all} \; t \in [0,t_1]
  .
\end{align*}
If, in addition, $f$ is continuous, then being weakly nonincreasing is
equivalent to the property of having a solution starting at $y$ such
that $t \mapsto f (x (t))$ is monotonically nonincreasing on
$[0,t_1]$.

Similarly, a lower semicontinuous function $\map{f}{\real^d}{\real}$
is \emph{strongly nonincreasing on $\domain$}\index{function!strongly
  nonincreasing} for a set-valued map
$\map{\setvaluedmap}{\real^d}{\parts{\real^d}}$ if, for all $y \in
\domain$, all solutions $\map{x}{[0,t_1]}{\real^d}$ of the
differential inclusion~\eqref{eq:auto-diff-inclusion} starting at $y$
and contained in $\domain$ satisfy
\begin{align*}
  f (x(t)) \le f (x (0)) = f (y) \quad \text{for all} \; t \in [0,t_1]
  .
\end{align*}
Note that $f$ is strongly nonincreasing if and only if $t \mapsto f (x
(t))$ is monotonically nonincreasing on $[0,t_1]$ for all solutions $t
\mapsto x(t)$ of the differential inclusion.

Given the set-valued map
$\map{\setvaluedmap}{\real^d}{\parts{\real^d}}$, which takes nonempty,
compact values, and the lower semicontinuous function
$\map{f}{\real^d}{\real}$, the \emph{lower and upper set-valued Lie
  derivatives}
\index{Lie derivative!lower set-valued}%
\sindex{SetLieDerLower}{$\setLiederlow{\setvaluedmap}{f}$}{Lower
  set-valued Lie derivative of the lower semicontinuous function
  $\map{f}{\real^d}{\real}$ with respect to the set-valued map
  $\map{\setvaluedmap}{\real^d}{\parts{\real^d}}$}%
\index{Lie derivative!upper set-valued}
\sindex{SetLieDerUpper}{$\setLiederup{\setvaluedmap}{f}$}{Upper
  set-valued Lie derivative of the lower semicontinuous function
  $\map{f}{\real^d}{\real}$ with respect to the set-valued map
  $\map{\setvaluedmap}{\real^d}{\parts{\real^d}}$}%
$\map{\setLiederlow{\setvaluedmap}{f},
  \setLiederup{\setvaluedmap}{f}}{\real^d}{\parts{\real}}$ of $f$ with
respect to $\setvaluedmap$ at $y$ are defined by, respectively,
\begin{align*}
  \setLiederlow{\setvaluedmap}{f}(y) & \triangleq \setdef{a \in
    \real}{\text{there exists} \; \zeta \in \partial_P f(y) \;
    \text{such that} \; a = \min
    \setdef{\zeta^T v}{v \in \setvaluedmap(y)}} , \\
  \setLiederup{\setvaluedmap}{f}(y) & \triangleq \setdef{a \in
    \real}{\text{there exists} \; \zeta \in \partial_P f(y) \;
    \text{such that} \; a = \max \setdef{\zeta^T v}{v \in
      \setvaluedmap(y)}} .
\end{align*}
If, in addition, $\setvaluedmap$ takes convex values, then for each
$\zeta \in \partial_P f(y)$, the set $ \setdef{\zeta^T v}{v \in
  \setvaluedmap(y)}$ is a closed interval of the form $[\min
\setdef{\zeta^T v}{v \in \setvaluedmap(y)},\max \setdef{\zeta^T v}{v
  \in \setvaluedmap(y)}]$. Note that the lower and upper set-valued
Lie derivatives at a point $y$ might be empty.

The lower and upper set-valued Lie derivatives play a role for a lower
semicontinuous function that is similar to the role played by the
set-valued Lie derivative $\setLiederG{\setvaluedmap}{f}$ for a
locally Lipschitz function.  These objects allow us to study how $f$
evolves along the solutions of a differential inclusion without having
to obtain the solutions in closed form.  Specifically, we have the
following result~\cite{FHC-YSL-RJS-PRW:98}. In this and in forthcoming
statements, it is convenient to adopt the convention $\sup \emptyset =
-\infty$.

\begin{proposition}\label{prop:weakly-strongly-decreasing}
  Let $\map{\setvaluedmap}{\real^d}{\parts{\real^d}}$ be a set-valued
  map satisfying the hypotheses of
  Proposition~\ref{prop:existence-solution}, and consider the
  associated differential inclusion~\eqref{eq:auto-diff-inclusion}.
  Let $\map{f}{\real^d}{\real}$ be a lower semicontinuous function,
  and let $\domain \subseteq \real^d$ be open. Then, the following
  statements hold:
  \begin{enumerate}
  \item The function $f$ is weakly nonincreasing on $\domain$ if and
    only if
    \begin{align*}
      \sup \setLiederlow{\setvaluedmap}{f} (y) \le 0, \quad \text{for
        all} \; y \in \domain .
    \end{align*}
  \item If, in addition, either $\setvaluedmap$ is locally Lipschitz
    on $\domain$, or $\setvaluedmap$ is continuous on $\domain$ and
    $f$ is locally Lipschitz on $\domain$, then $f$ is strongly
    nonincreasing on $\domain$ if and only if
    \begin{align*}
      \sup \setLiederup{\setvaluedmap}{f} (y) \le 0 , \quad \text{for
        all} \; y \in \domain .
    \end{align*}
  \end{enumerate}
\end{proposition}

Let us illustrate this result in a particular example.

\refstepcounter{example}%
\subsubsection*{Example~\theexample: Cart on a circle}\label{ex:cart}

Consider, following~\cite{ZA:83,EDS:99}, the control system on
$\real^2$
\begin{align}
  \dot x_1 & = (x_1^2 - x_2^2) u , \label{eq:cart-a}\\
  \dot x_2 & = 2 x_1 x_2 u ,  \label{eq:cart-b}
\end{align}
with $u \in \real$.  Equations~\eqref{eq:cart-a}-\eqref{eq:cart-b} are
named after the fact that the integral curves of the vector field
$(x_1,x_2) \mapsto (x_1^2 - x_2^2 ,2 x_1 x_2)$ are circles with center
on the $x_2$-axis and tangent to the $x_1$-axis, see the phase
portrait in Figure~\ref{fig:cart}(a).

Let $\map{X}{\real^2 \times \real}{\real^2}$ be defined by
$X((x_1,x_2), u) = (x_1^2 - x_2^2 ,2 x_1 x_2) \, u$.
Following~\eqref{eq:control-set-valued-map}, consider the associated
set-valued map $\map{G[X]}{\real^2}{\parts{\real^2}}$ defined by
$G[X](x_1,x_2) = \setdef{ X((x_1,x_2), u)}{u \in \real}$.  Since
$G[X]$ does not take compact values, we instead take a nondecreasing
map $\map{\sigma}{\realnonnegative}{\realnonnegative}$, and consider
the set-valued map
$\map{\setvaluedmap_\sigma}{\real^2}{\parts{\real^2}}$ given by
\begin{align}\label{eq:set-valued-map-cart}
  \setvaluedmap_\sigma(x_1,x_2) = \setdef{ X((x_1,x_2), u) \in
    \real^2}{|u| \le \sigma (\TwoNorm{(x_1,x_2)})}.
\end{align}

Consider the locally Lipschitz function $\map{f}{\real^2}{\real}$
defined by
\begin{align*}
  f(x_1,x_2) =
  \begin{cases}
    \frac{x_1^2 + x_2^2}{\sqrt{x_1^2 + x_2^2} + |x_1|} , & x \neq
    0 ,\\
    0 , & x = 0 .
  \end{cases}
\end{align*}
The level set curves of $f$ are depicted in Figure~\ref{fig:cart}(c).
Let us determine how $f$ evolves along the solutions of the control
system by using the lower and upper set-valued Lie derivatives. First,
we compute the proximal subdifferential of $f$. Using the fact
that $f$ is twice continuously differentiable on the open right and
left half-planes, together with the geometric interpretation of
proximal subgradients, we obtain
\begin{align*}
  \partial_P f (x_1,x_2) =
  \begin{cases}
    \Big \{ \Big( - \frac{x_1^2+x_2^2 -2
      x_1\sqrt{x_1^2+x_2^2}}{x_1^2+x_2^2 +x_1\sqrt{x_1^2+x_2^2}} ,
    \frac{ x_2 (2 x_1 + \sqrt{x_1^2 + x_2^2})}{(x_1 + \sqrt{x_1^2
        + x_2^2})^2} \Big) \Big \} , & x_1>0 ,\\
    \emptyset, & x_1 = 0, \\
    \Big \{ \Big( \frac{x_1^2+x_2^2 +2
      x_1\sqrt{x_1^2+x_2^2}}{x_1^2+x_2^2 -x_1\sqrt{x_1^2+x_2^2}} ,
    \frac{ x_2 (-2 x_1 + \sqrt{x_1^2 + x_2^2})}{(x_1 - \sqrt{x_1^2
        + x_2^2})^2} \Big) \Big \} , & x_1<0 .
  \end{cases}
\end{align*}
With this information, we compute the set
\begin{multline*}
   \setdef{\zeta^T v}{\zeta \in \partial_P f(x_1,x_2) , \; v \in
    \setvaluedmap_\sigma(x_1,x_2)} \\
  =
  \begin{cases}
    \big \{ u \frac{(x_1^2+x_2^2)^2}{x_1^2+x_2^2 +x_1
      \sqrt{x_1^2+x_2^2}} \; : \; |u| \le \sigma (\TwoNorm{(x_1,x_2)})
    \big \},  & x_1 > 0 , \\
    \emptyset ,  & x_1 = 0 ,\\
    \big \{ - u \frac{(x_1^2+x_2^2)^2}{x_1^2+x_2^2 -x_1
      \sqrt{x_1^2+x_2^2}} \; : \; |u| \le \sigma (\TwoNorm{(x_1,x_2)})
    \big \} , & x_1 < 0 .
  \end{cases}
\end{multline*}
We are now ready to compute the lower and upper set-valued Lie
derivatives as
\begin{align}
  \setLiederlow{\setvaluedmap_\sigma}{f} (x_1,x_2) &=
  \begin{cases}
    -\sigma (\TwoNorm{(x_1,x_2)})
    \frac{(x_1^2+x_2^2)^{3/2}}{\sqrt{x_1^2+x_2^2}
      + |x_1|} , & x_1 \neq 0  ,\\
    -\infty, & x_1 = 0 ,
  \end{cases}
  \label{eq:setLiederlow}
  \\
  \setLiederup{\setvaluedmap_\sigma}{f} (x_1,x_2) & =
  \begin{cases}
    \sigma (\TwoNorm{(x_1,x_2)})
    \frac{(x_1^2+x_2^2)^{3/2}}{\sqrt{x_1^2+x_2^2}
      + |x_1|} , & x_1 \neq 0  ,\\
    -\infty, & x_1 = 0 .
  \end{cases}
  \label{eq:setLiederup}
\end{align}
Therefore, $\sup \setLiederlow{\setvaluedmap_\sigma}{f}(x_1,x_2) \le
0$ for all $(x_1,x_2) \in \real^2$.  Now using
Proposition~\ref{prop:weakly-strongly-decreasing}\textit{(i)}, we
deduce that $f$ is weakly nonincreasing on $\real^2$. Since $f$ is
continuous, this fact is equivalent to saying that there exists a
control input $u$ such that the solution $t \mapsto x(t)$ of the
dynamical system resulting from substituting $u$
in~\eqref{eq:cart-a}-\eqref{eq:cart-b} has the property that $t
\mapsto f (x (t))$ is monotonically nonincreasing.
\hfill
$\blacksquare$

\bigskip
\noindent \textbf{\textit{Stability
    results}}
\medskip

The results presented in the previous section establishing the
monotonic behavior of a lower semicontinuous function allow us to
provide tools for stability analysis. We present here an exposition
parallel to the discussion for the locally Lipschitz function and
generalized gradient case. We start by presenting a result on Lyapunov
stability that follows from the exposition
in~\cite{FHC-YSL-RJS-PRW:98}.

\begin{theorem}\label{th:Lyapunov-prox}
  Let $\map{\setvaluedmap}{\real^d}{\parts{\real^d}}$ be a set-valued
  map satisfying the hypotheses of
  Proposition~\ref{prop:existence-solution}.  Let $x_e$ be an
  equilibrium of the differential
  inclusion~\eqref{eq:auto-diff-inclusion}, and let $\domain \subseteq
  \real^d$ be a domain with $x_e \in \domain$.  Let
  $\map{f}{\real^d}{\real}$ and assume that the following conditions
  hold:
  \begin{enumerate}
  \item $\setvaluedmap$ is continuous on $\domain$ and $f$ is locally
    Lipschitz on $\domain$, or $\setvaluedmap$ is locally Lipschitz on
    $\domain$ and $f$ is lower semicontinuous on $\domain$, and $f$ is
    continuous at $x_e$.
  \item $f(x_e)=0$, and $f(x)>0$ for $x \in \domain \! \setminus \!
    \{x_e\}$.
  \item $\sup \setLiederup{\setvaluedmap}{f}(x) \le 0$ for all $x \in
    \domain$.
  \end{enumerate}
  Then, $x_e$ is a strongly stable equilibrium
  of~\eqref{eq:auto-diff-inclusion}. In addition, if \textit{(iii)}
  above is replaced by
  \begin{enumerate}
  \item[(iii)'] $\sup \setLiederup{\setvaluedmap}{f}(x) < 0$ for all
    $x \in \domain \! \setminus \! \{x_e \}$,
  \end{enumerate}
  then $x_e$ is a strongly asymptotically stable equilibrium
  of~\eqref{eq:auto-diff-inclusion}.
\end{theorem}

A similar result can be stated for weakly stable equilibria
substituting \textit{(i)} by ``\textit{(i)'} $f$ is continuous on
$\domain$,'' and the upper set-valued Lie derivative by the lower
set-valued Lie derivative in \textit{(iii)} and \textit{(iii)'}.  Note
that, if the differential inclusion~\eqref{eq:auto-diff-inclusion} has
a unique solution starting from every initial conditions, then the
notions of strong and weak stability coincide, and it is sufficient to
satisfy the simpler requirements \textit{(i)'} and \textit{(iii)'} for
weak stability.

Similarly to the case of a continuous differential equation, global
asymptotic stability can be established by requiring the Lyapunov
function $f$ to be continuous and radially unbounded.  The equivalence
between global strong asymptotic stability and the existence of
infinitely differentiable Lyapunov functions is discussed
in~\cite{FHC-YSL-RJS:98}.  This type of global result is commonly
invoked when dealing with the stabilization of control systems by
referring to control Lyapunov functions~\cite{EDS:99} or Lyapunov
pairs~\cite{FHC-YSL-RJS-PRW:98}.  Two lower semicontinuous functions
$\map{f,g}{\real^d}{\real}$ are a \emph{Lyapunov pair}\index{Lyapunov
  pair} for an equilibrium $x_e \in \real^d$ if they satisfy $f(x),
g(x) \ge 0$ for all $x \in \real^d$, and $g(x) = 0$ if and only if $x
= x_e$; $f$ is radially unbounded, and, moreover,
\begin{align*}
  \sup \setLiederlow{\setvaluedmap}{f}(x) \le - g(x) \quad \text{for
    all} \; x \in \real^d .
\end{align*}
If an equilibrium $x_e$ of~\eqref{eq:auto-diff-inclusion} admits a
Lyapunov pair, then there exists at least one solution starting from
every initial condition that asymptotically converges to the
equilibrium, see~\cite{FHC-YSL-RJS-PRW:98}.

\subsubsection*{Example~\ref{ex:cart} revisited: Asymptotic stability
  of the origin in the cart example}

As an application of the above discussion and the version of
Theorem~\ref{th:Lyapunov-prox} for weak stability, consider the cart
on a circle introduced in Example~\ref{ex:cart}. Setting $x_e = (0,0)$
and $\domain = \real^2$, and taking into account the
computation~\eqref{eq:setLiederlow} of the lower set-valued Lie
derivative, we conclude that $(0,0)$ is a globally weakly
asymptotically stable equilibrium.
\hfill
$\blacksquare$

We now turn our attention to the extension of the invariance principle
for differential inclusions using the proximal subdifferential of a
lower semicontinuous function. The following result can be derived
from the exposition in~\cite{FHC-YSL-RJS-PRW:98}.

\begin{theorem}\label{th:LaSalle-prox}
  Let $\map{\setvaluedmap}{\real^d}{\parts{\real^d}}$ be a set-valued
  map satisfying the hypotheses of
  Proposition~\ref{prop:existence-solution}, and let
  $\map{f}{\real^d}{\real}$.  Assume that either $\setvaluedmap$ is
  continuous and $f$ is locally Lipschitz, or $\setvaluedmap$ is
  locally Lipschitz and $f$ is continuous. Let $S \subset \real^d$ be
  compact and strongly invariant for~\eqref{eq:auto-diff-inclusion},
  and assume that $\sup \setLiederup{\setvaluedmap}{f}(y) \le 0$ for
  all $y \in S$.  Then, every solution $\map{x}{[0,\infty)}{\real^d}$
  of~\eqref{eq:auto-diff-inclusion} starting at $S$ converges to the
  largest weakly invariant set $M$ contained in
  \begin{align*}
    S \cap \overline{\setdef{y \in \real^d}{0 \in
        \setLiederup{\setvaluedmap}{f}(y)}} .
  \end{align*}
  Moreover, if the set $M$ consists of a finite number of points, then
  the limit of each solution starting in $S$ exists and is an element
  of $M$.
\end{theorem}

Next, we apply Theorem~\ref{th:LaSalle-prox} to gradient differential
inclusions.

\bigskip
\noindent
\textbf{\textit{Stability of gradient differential inclusions for
    convex functions}} \medskip

Consider the gradient differential
inclusion~\eqref{eq:gradient-inclusion} associated with the convex
function $\map{f}{\real^d}{\real}$. We study here the asymptotic
behavior of the solutions of~\eqref{eq:gradient-inclusion}.  From the
section ``Gradient differential inclusion of a convex function,'' we
know that solutions exist and are unique. Consequently, the notions of
weakly nonincreasing and strongly nonincreasing function coincide.
Therefore, it suffices to show that $f$ is weakly nonincreasing on
$\real^d$ for the gradient differential
inclusion~\eqref{eq:gradient-inclusion}.

For all $\zeta \in \partial_P f(x)$, there exists $v = -\zeta \in -
\partial_P f(x)$ such that $ \zeta^T v = - \TwoNorm{\zeta}^2 \le 0$.
Hence,
\begin{align*}
  \setLiederlow{-\partial_P f}{f}(x) \le 0 \quad \text{for all} \; x
  \in \real^d.
\end{align*}
Proposition~\ref{prop:weakly-strongly-decreasing}\textit{(i)} now
guarantees that $f$ is weakly nonincreasing on $\real^d$. Since the
solutions of the gradient differential
inclusion~\eqref{eq:gradient-inclusion} are unique, $t \mapsto f (x
(t))$ is monotonically nonincreasing for all solutions $t \mapsto
x(t)$.

The application of the Lyapunov stability theorem and the invariance
principle in theorems~\ref{th:Lyapunov-prox}
and~\ref{th:LaSalle-prox}, respectively, gives rise to the following
nonsmooth counterpart of the classical smooth result for a gradient
flow~\cite{WMH-SS:74}.

\begin{proposition}
  Let $\map{f}{\real^d}{\real}$ be convex. Then, the strict minimizers
  of~$f$ are strongly stable equilibria of the gradient differential
  inclusion~\eqref{eq:gradient-inclusion} associated with $f$.
  Furthermore, if the level sets of $f$ are bounded, then the
  solutions of the gradient differential inclusion asymptotically
  converge to the set of minimizers of $f$.
\end{proposition}

Convergence rate estimates of~\eqref{eq:gradient-inclusion} can be
found in~\cite{OG:05}.

\bigskip
\noindent \textbf{\textit{Stabilization of control systems}}
\medskip

Consider the control system on $\real^d$ given by
\begin{align}
  \label{eq:control-system}
  \dot x = X(x,u) ,
\end{align}
where $\map{X}{\real^d \times \real^m}{\real^d}$.  The
system~\eqref{eq:control-system} is locally (respectively, globally)
\emph{continuously stabilizable} if there exists a continuous map
$\map{k}{\real^d}{\real^m}$ such that the closed-loop system
\begin{align*}
  \dot x = X(x,k(x))  
\end{align*}
is locally (respectively globally) asymptotically stable at the
origin. The following result by Brockett~\cite{RWB:83a}, see
also~\cite{EDS:98,AB-LR:05}, states that a large class of control
systems are not stabilizable by a continuous feedback controller.
\begin{theorem}\label{th:Brockett}
  Let $\map{X}{\real^d \times \real^m}{\real^d}$ be continuous and
  $X(0,0)=0$. If there exists a continuous stabilizer of the control
  system~\eqref{eq:control-system}, then there exists a neighborhood
  of the origin in $\real^d \times \real^m$ whose image by $X$ is a
  neighborhood of the origin in $\real^d$.
\end{theorem}

In particular, Theorem~\ref{th:Brockett} implies that control systems
of the form
\begin{align}\label{eq:driftless}
  \dot x = u_1 X_1(x) + \dots + u_m X_m(x) ,
\end{align}
with $m<n$ and $\map{X_i}{\real^d}{\real^d}$, $i \in \until{m}$,
continuous with $\text{rank}(X_1(0),\dots,X_m(0)) = m$, cannot be
stabilized by a continuous feedback controller.

\refstepcounter{example}%
\subsubsection*{Example~\theexample: Nonholonomic
  integrator}\label{ex:Brockett}

Consider the nonholonomic integrator~\cite{RWB:81}
\begin{align}
  \dot x_1 &= u_1 , \label{eq:Brockett-1} \\
  \dot x_2 &= u_2 , \label{eq:Brockett-2} \\
  \dot x_3 &= x_1 u_2 - x_2 u_1 , \label{eq:Brockett-3}
\end{align}
where $(u_1,u_2) \in \real^2$. The nonholonomic integrator is a
control system of the form~\eqref{eq:driftless}, with $m=2$,
$X_1(x_1,x_2,x_3) = (1,0,-x_2)$, and $X_2(x_1,x_2,x_3) = (0,1,x_1)$.
The nonholonomic integrator is controllable, that is, for each pair of
states $x,y \in \real^3$, there exists an input $t \mapsto u(t)$ such
that the corresponding solution
of~\eqref{eq:Brockett-1}-\eqref{eq:Brockett-3} starting at $x$ reaches
$y$. However, the nonholonomic integrator does not satisfy the
condition in Theorem~\ref{th:Brockett}, and thus is not continuously
stabilizable.  Indeed, no point of the form $(0,0,\eps)$, where $\eps
\neq 0$, belongs to the image of the map $\map{X}{\real^3 \times
  \real^2}{\real^3}$, $X(x_1,x_2,x_3,u_1,u_2) = (u_1,u_2, x_1 u_2 -
x_2 u_1)$.
\hfill
$\blacksquare$

The condition in Theorem~\ref{th:Brockett} is necessary but not
sufficient. There exist control systems that satisfy the condition,
and still cannot be stabilized by means of a continuous stabilizer.
The cart on a circle in Example~\ref{ex:cart} is one of these systems.
The map $((x_1,x_2),u) \to X(x_1,x_2,u) $ satisfies the necessary
condition in Theorem~\ref{th:Brockett} but cannot be stabilized with a
continuous $\map{k}{\real^2}{\real}$, see~\cite{EDS:99}.

Another obstruction to the existence of continuous stabilizing
controllers is given by Milnor's theorem~\cite{EDS:98}, which states
that the domain of attraction of an asymptotically stable equilibrium
of a locally Lipschitz vector field must be diffeomorphic to Euclidean
space.  Since environments with obstacles are not diffeomorphic to
Euclidean space, one can use Milnor's theorem to justify~\cite{EDS:99}
the nonexistence of continuous globally stabilizing controllers in
environments with obstacles.

The obstruction to the existence of continuous stabilizers has
motivated the search for time-varying and discontinuous feedback
stabilizers.  Regarding the latter, given a discontinuous
$\map{k}{\real^d}{\real^m}$, the immediate question that arises is how
to understand the solutions of the resulting discontinuous dynamical
system $\dot x = X(x,k(x))$.  From the previous discussion, we know
that Caratheodory solutions are not a good candidate, since in many
situations they fail to exist. The following
result~\cite{EPR:94bis,JMC-LR:94} shows that Filippov solutions are
also not a good candidate.

\begin{theorem}
  Let $\map{X}{\real^d \times \real^m}{\real^d}$ be continuous and
  $X(0,0)=0$. Assume that, for each $U \subseteq \real^m$ and each $x
  \in \real^d$, $X(x,\co U) = \co X(x,U)$ holds. With solutions
  understood in the Filippov sense, if there exists a measurable,
  locally bounded stabilizer of the control
  system~\eqref{eq:control-system}, then there exists a neighborhood
  of the origin in $\real^d \times \real^m$ whose image by $X$ is a
  neighborhood of the origin in $\real^d$.
\end{theorem}

In particular, control systems of the form~\eqref{eq:driftless} cannot
be stabilized by means of a discontinuous feedback if solutions are
understood in the Filippov sense. This impossibility result, however,
can be overcome if solutions are understood in the sample-and-hold
sense, as shown in~\cite{FHC-YSL-EDS-AIS:97}.  Let us briefly discuss
this result in the light of the above exposition.  For more details,
see~\cite{FHC-YSL-RJS-PRW:98,FHC:04,EDS:99}.  Consider the
differential inclusion~\eqref{eq:control-inclusion} associated with
the control system~\eqref{eq:control-system}.  The
system~\eqref{eq:control-system} is \emph{(open loop) globally
  asymptotically controllable} (to the origin) if $0$ is a Lyapunov
weakly stable equilibrium of the differential
inclusion~\eqref{eq:control-inclusion}, and every point $y \in
\real^d$ has the property that there exists a solution
of~\eqref{eq:control-inclusion} satisfying $x(0) = y$ and $\lim_{t \to
  \infty} x(t) = 0$.  On the other hand, a map
$\map{k}{\real^d}{\real^m}$ \emph{stabilizes} the
system~\eqref{eq:control-system} \emph{in the sample-and-hold sense}
if, for all $x_0 \in \real^d$ and all $\eps \in \realpositive$, there
exist $\delta, T \in \realpositive$ such that, for all partitions
$\pi$ of $[0,t_1]$ with $\diam (\pi) < \delta$, the corresponding
$\pi$-solution $t \mapsto x(t)$ of~\eqref{eq:control-system} starting
at $x_0$ satisfies $\TwoNorm{x(t)} \le \eps$ for all $t \ge T$.

\index{globally asymptotically controllable}

The following result~\cite{FHC-YSL-EDS-AIS:97} states that both
notions, global asymptotic controllability and the existence of a
feedback stabilizer in the sample-and-hold sense, are equivalent.
\begin{theorem}\label{th:hard-equivalence}
  Let $\map{X}{\real^d \times \real^m}{\real^d}$ be continuous and
  $X(0,0)=0$. Then, the control system~\eqref{eq:control-system} is
  globally asymptotically controllable if and only if it admits a
  measurable, locally bounded stabilizer in the sample-and-hold sense.
\end{theorem}

The ``only if'' implication is clear. The converse implication is
proved by explicit construction of the stabilizer, and is based on the
fact that the control system~\eqref{eq:control-system} is globally
asymptotically controllable if and only if it admits a continuous
Lyapunov pair, see~\cite{EDS:83}. Using the continuous Lyapunov
function provided by this characterization, the discontinuous feedback
for the control system~\eqref{eq:control-system} can be constructed
explicitly~\cite{FHC-YSL-EDS-AIS:97,EDS:99}.

The existence of a Lyapunov pair in the sense of generalized gradients
(that is, when the set-valued Lie derivative involving the generalized
gradient is used instead of the lower set-valued Lie derivative
involving the proximal subdifferential), however, turns out to be
equivalent to the existence of a stabilizing feedback in the sense of
Filippov, see~\cite{LR:01}.

\subsubsection*{Example~\ref{ex:cart} revisited: Cart stabilization by
  discontinuous feedback}

As an illustration, consider Example~\ref{ex:cart}. We have already
shown that $(0,0)$ is a globally weakly asymptotically stable
equilibrium of the differential
inclusion~\eqref{eq:set-valued-map-cart} associated with the control
system. Therefore, the control system is globally asymptotically
controllable, and can be stabilized in the sample-and-hold sense by
means of a discontinuous feedback. The stabilizing feedback that
results from the proof of Theorem~\ref{th:hard-equivalence} can be
described as follows, see~\cite{GAL-EDS:93,EDS:99}. If placed to the
left of the $x_2$ axis, move in the direction of the vector field~$g$.
If placed to the right of the $x_2$ axis, move in the opposite
direction of the vector field $g$. Finally, on the $x_2$-axis, choose
an arbitrary direction.  The stabilizing nature of this feedback can
be graphically checked in Figure~\ref{fig:cart}(a) and~(b).
\hfill $\blacksquare$

Remarkably, for a system that is affine in the controls, it is
possible to show~\cite{LR:02} that there exists a stabilizing feedback
controller whose discontinuities form a set of measure zero, and,
moreover, the discontinuity set is repulsive for the solutions of the
closed-loop system. In particular, this fact means that, for the
closed-loop system, the solutions can be understood in the
Caratheodory sense.  This situation exactly corresponds to the
situation in Example~\ref{ex:cart}.

\section*{Conclusions}

This article has presented an introductory tutorial on discontinuous
dynamical systems.  Various examples illustrate the pertinence of the
continuity and Lipschitzness properties that guarantee the existence
and uniqueness of classical solutions to ordinary differential
equations. The lack of these properties in examples drawn from various
disciplines motivates the need for more general notions than the
classical one.  This observation is the starting point into the three
main themes of our discussion. First, we introduced notions of
solution for discontinuous systems.  Second, we reviewed the available
tools from nonsmooth analysis to study the gradient information of
candidate Lyapunov functions. And, third, we presented nonsmooth
stability tools to characterize the asymptotic behavior of solutions.

The physical significance of the solution notions discussed in this
article depend on the specific setting.  Caratheodory solutions are
employed for time-dependent vector fields that depend discontinuously
on time, such as dynamical systems involving impulses as well as
control systems with discontinuous open-loop inputs.  Filippov
solutions are used in problems involving electrical circuits with
switches, relay control, friction, and sliding. This observation is
also valid for solution notions similar to Filippov's, such as
Krasovskii and Sentis solutions, see ``Additional Solution Notions for
Discontinuous Systems.''  The notion of a $\pi$-solution for control
systems has the physical interpretation of iteratively evaluating the
input at the current state and holding it steady for some time while
the closed-loop dynamical system evolves. As illustrated above, this
solution notion plays a pivotal role in the stabilization of
asymptotically controllable systems.

There are numerous important issues that are not treated here;
``Additional Topics on Discontinuous Systems and Differential
Inclusions'' lists some of them.  The topic of discontinuous dynamical
systems is vast, and our focus on the above-mentioned themes is aimed
at providing a coherent exposition. We hope that this tutorial serves
as a guided motivation for the reader to further explore the exciting
topic of discontinuous systems. The list of references of this
manuscript provides a good starting point to undertake this endeavor.

\section*{Acknowledgments}
This research was supported by NSF CAREER Award ECS-0546871. The
author wishes to thank Francesco Bullo and Anurag Ganguli for
countless hours of fun with Filippov solutions, and Bernard Brogliato,
Rafal Goebel, Rishi Graham, and three anonymous reviewers for numerous
comments that improved the presentation.

\clearpage

{
  \psfrag{a1}[cc][cc]{{\tiny $x$}}%
  \psfrag{a2}[cc][cc]{{\tiny $u$}}%
  \psfrag{-3}[cc][cc]{{\tiny $-3$}}%
  \psfrag{-2}[cc][cc]{{\tiny $-2$}}%
  \psfrag{-1}[cc][cc]{{\tiny $-1$}}%
  \psfrag{1}[cc][cc]{{\tiny $1$}}%
  \psfrag{2}[cc][cc]{{\tiny $2$}}%
  \psfrag{3}[cc][cc]{{\tiny $3$}}%
  \begin{figure}
    \centering
    \includegraphics[width=.475\linewidth]{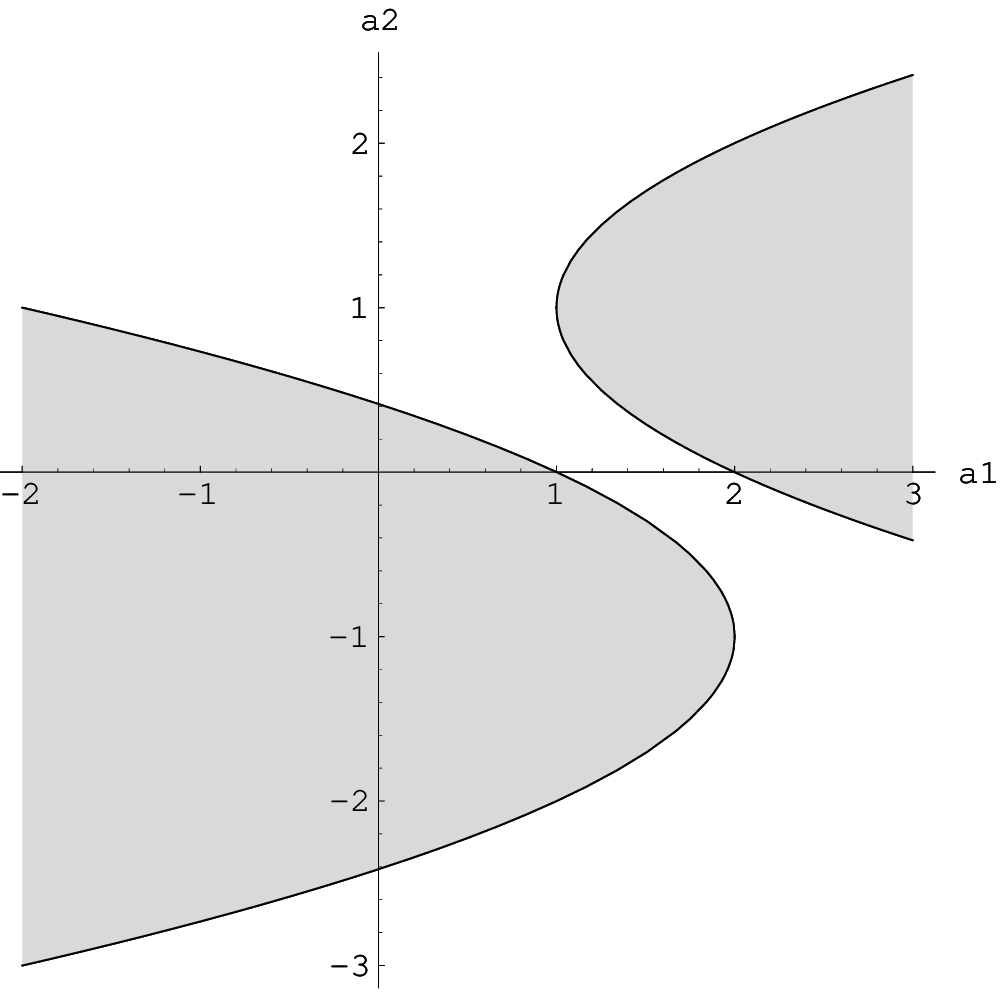}
    \caption{Obstruction to the existence of continuous state-dependent
      feedback~\cite{EDS:98}. The shaded areas represent the values of
      the control $u$ that are needed to stabilize the equilibrium point
      $0$ of system~\eqref{eq:sontag}. There does not exist a function
      $x \mapsto u(x)$ defined on $\real$ that at the same time is
      continuous and belongs to the shaded areas.}
    \label{fig:sontag}
  \end{figure}
}

\clearpage

{
  \psfrag{4}[cc][cc]{{\tiny $4$}}%
  \psfrag{3.5}[cc][cc]{{\tiny $3.5$}}%
  \psfrag{3}[cc][cc]{{\tiny $3$}}%
  \psfrag{2.5}[cc][cc]{{\tiny $2.5$}}%
  \psfrag{2}[cc][cc]{{\tiny $2$}}%
  \psfrag{1.5}[cc][cc]{{\tiny $1.5$}}%
  \psfrag{1}[cc][cc]{{\tiny $1$}}%
  \psfrag{0.5}[cc][cc]{{\tiny $0.5$}}%
  \psfrag{0}[cc][cc]{{\tiny $0$}}%
  \psfrag{-1}[cc][cc]{{\tiny $-1$}}%
  \psfrag{-2}[cc][cc]{{\tiny $-2$}}%
  \psfrag{a1}[cc][cc]{{\footnotesize $v$}}%
  \psfrag{a2}[cc][cc]{{\footnotesize $\nu$}}%
  \begin{figure}[htbp]
    \centering 
    \subfigure[]{
      \resizebox{.425\linewidth}{!}{\input{figure2a.tex}}
    }\qquad
    \subfigure[]{
      \includegraphics[width=.475\linewidth]{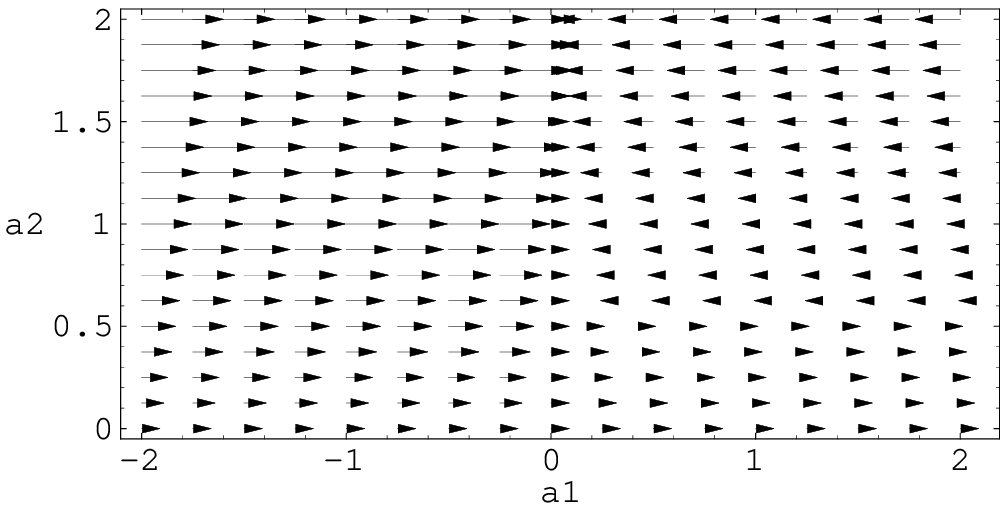}}
    \caption{Brick sliding on a frictional ramp.  (a) shows the
      physical quantities used to describe the example.  (b) shows the
      one-dimensional phase portraits of~\eqref{eq:brick}
      corresponding to values of the friction coefficient $\nu$
      between $0$ and $2$, with a ramp inclination of $30$ degrees.
      The phase portrait shows that, for sufficiently small values of
      $\nu$, every trajectory that starts with positive initial
      velocity (``moving to the right'') never stops.  However, for
      sufficiently large values of $\nu$, every trajectory that starts
      with positive initial velocity eventually stops and remains
      stopped.  However, there is no continuously differentiable
      solution of~\eqref{eq:brick} that exhibits this type of
      behavior.  We thus need to expand our notion of solution beyond
      continuously differentiable solutions by understanding the
      effect of the discontinuity
      in~\eqref{eq:brick}.}\label{fig:brick}
  \end{figure}
}

\clearpage

{
\psfrag{1}[cc][cc]{{\tiny $1$}}%
\psfrag{0.5}[cc][cc]{{\tiny $.5$}}%
\psfrag{0}[cc][cc]{{\tiny $0$}}%
\psfrag{-0.5}[cc][cc]{{\tiny $-.5$}}%
\psfrag{-1}[cc][cc]{{\tiny $-1$}}%
\psfrag{a1}[cc][cc]{{\footnotesize $x_1$}}%
\psfrag{a2}[cc][cc]{{\footnotesize $x_2$}}%
\begin{figure}[htbp]
  \centering%
  \subfigure[]{
    \includegraphics[width=.355\linewidth]{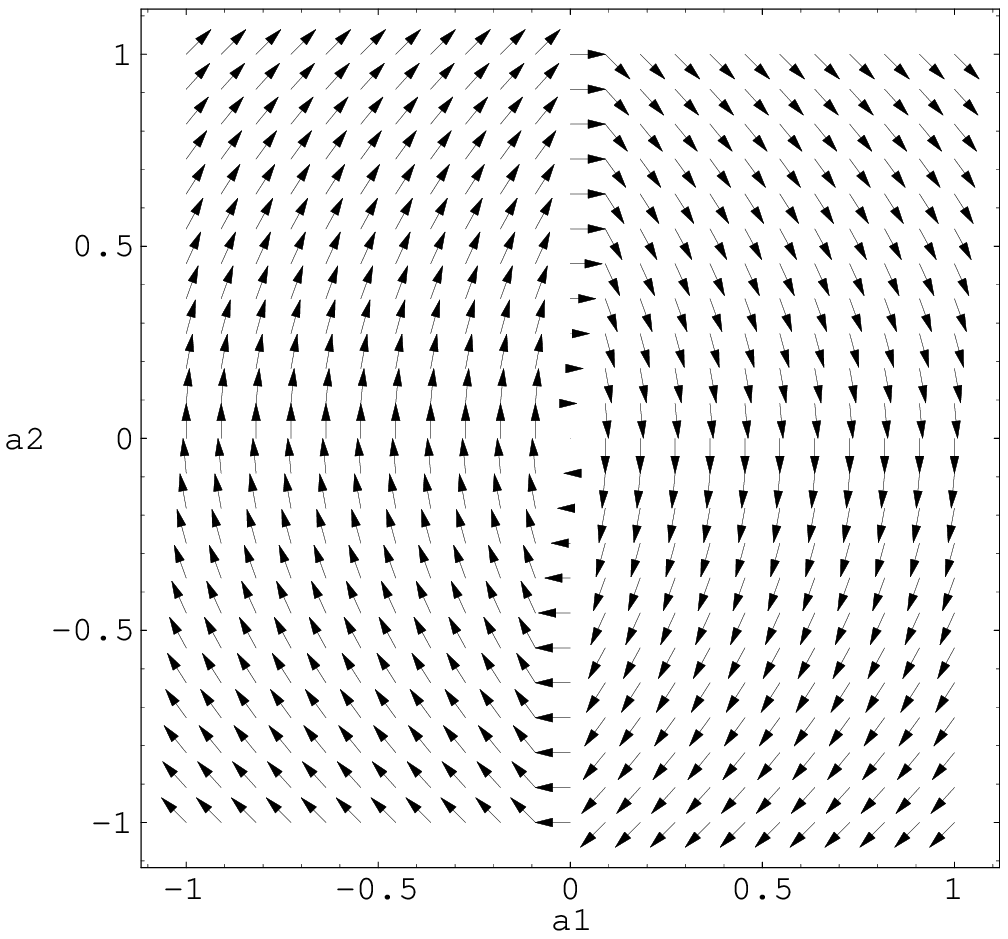}}\qquad
  \subfigure[]{
    \includegraphics[width=.35\linewidth]{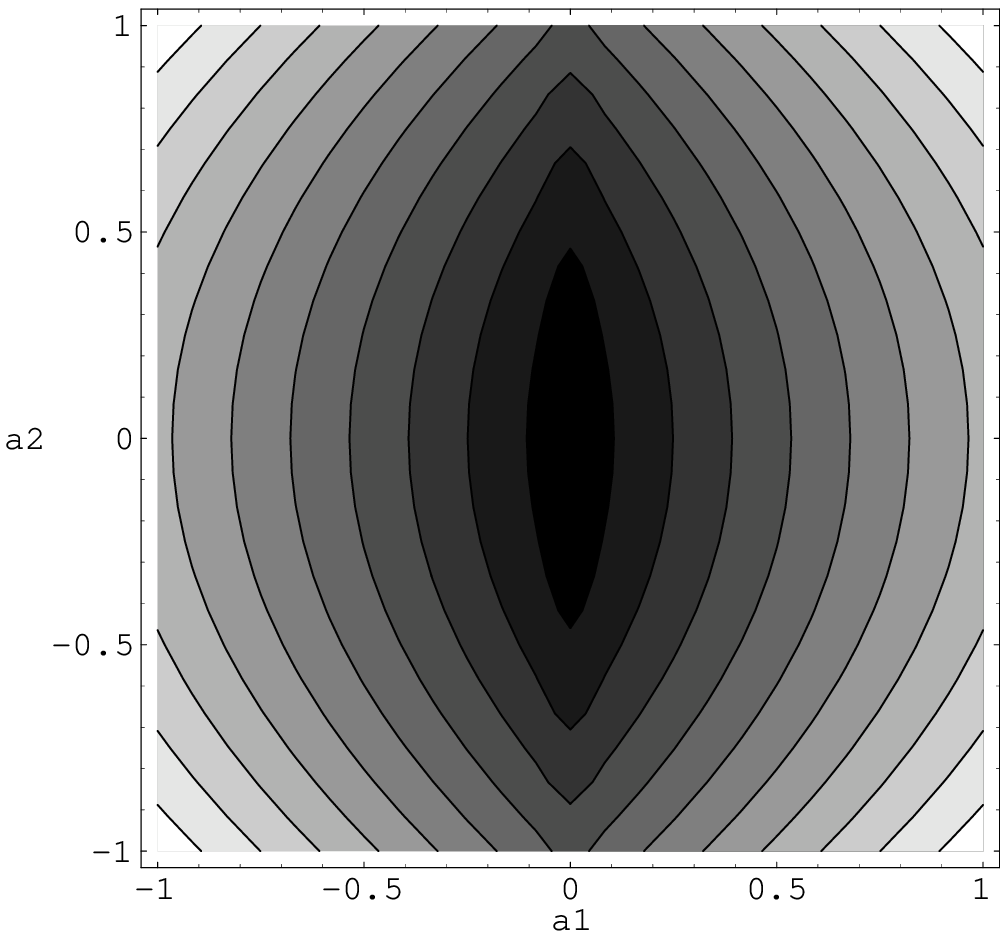}}
  \caption{Nonsmooth harmonic oscillator.  (a) shows the phase
    portrait on $[-1,1]^2$ of the vector field $(x_1,x_2) \mapsto
    (x_2,-\sign (x_1))$, while (b) shows the contour plot on
    $[-1,1]^2$ of the function $(x_1,x_2) \mapsto |x_1| +
    \frac{x_2^2}{2}$. The discontinuity of the vector field along the
    $x_2$-coordinate axis makes it impossible to find continuously
    differentiable solutions of~\eqref{eq:ODE-auto}. On the other
    hand, the level sets in (b) match the phase portrait in (a)
    everywhere except for the $x_2$-coordinate axis, which suggests
    that trajectories along the level sets are candidates for
    solutions of~\eqref{eq:ODE-auto} in a sense different from the
    classical one.}
  \label{fig:oscillator}
\end{figure}
}

\clearpage

\clearpage

{
\psfrag{1}[cc][cc]{{\tiny $1$}}%
\psfrag{0.5}[cc][cc]{{\tiny $.5$}}%
\psfrag{-0.5}[cc][cc]{{\tiny $-.5$}}%
\psfrag{-1}[cc][cc]{{\tiny $-1$}}%
\psfrag{-0.1}[cc][cc]{{\tiny $-.1$}}%
\psfrag{0.1}[cc][cc]{{\tiny $.1$}}%
\psfrag{-0.2}[cc][cc]{{\tiny $-.2$}}%
\psfrag{0.2}[cc][cc]{{\tiny $.2$}}%
\psfrag{0.3}[cc][cc]{{\tiny $.3$}}%
\psfrag{0.4}[cc][cc]{{\tiny $.4$}}%
\psfrag{0.6}[cc][cc]{{\tiny $.6$}}%
\psfrag{0.8}[cc][cc]{{\tiny $.8$}}%
\begin{figure}[htbp]
  \centering
  \subfigure[]{\includegraphics[width=.25\linewidth]{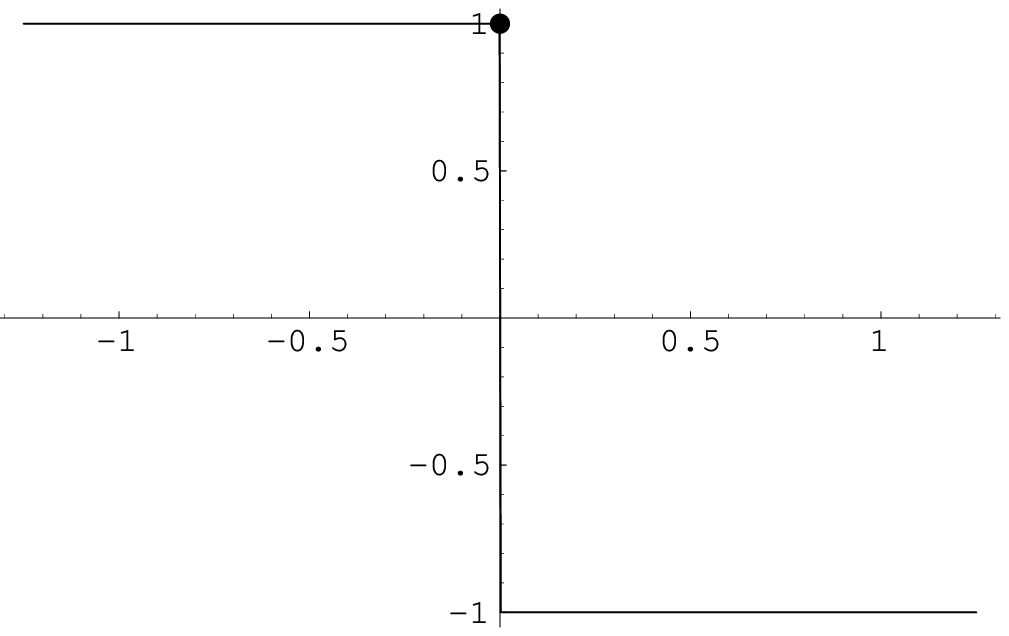}}%
  \subfigure[]{\includegraphics[width=.25\linewidth]{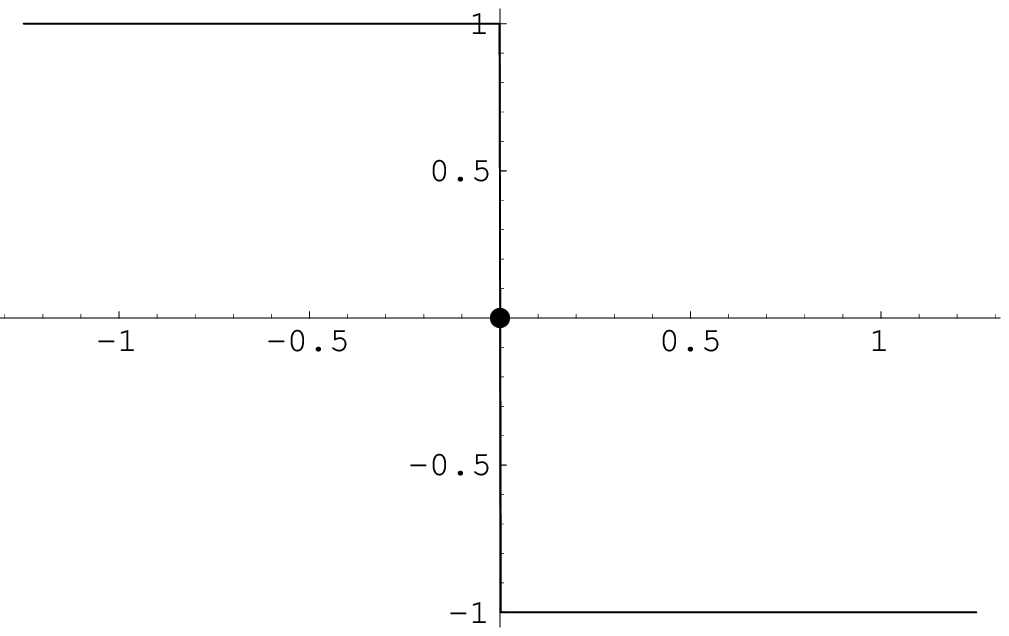}}%
  \subfigure[]{\includegraphics[width=.25\linewidth]{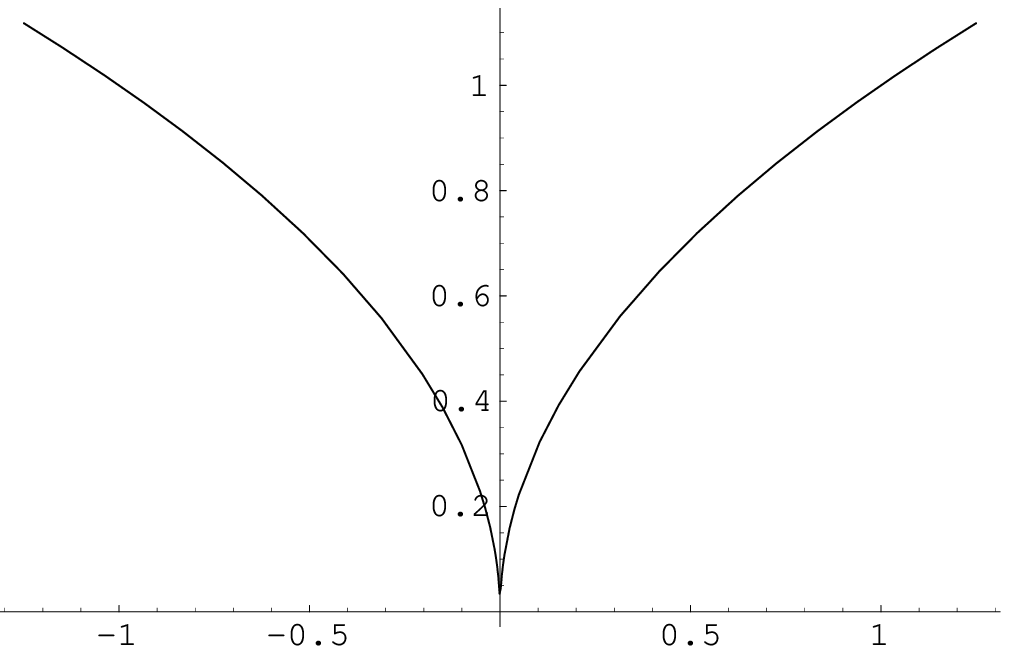}}%
  \subfigure[]{\includegraphics[width=.25\linewidth]{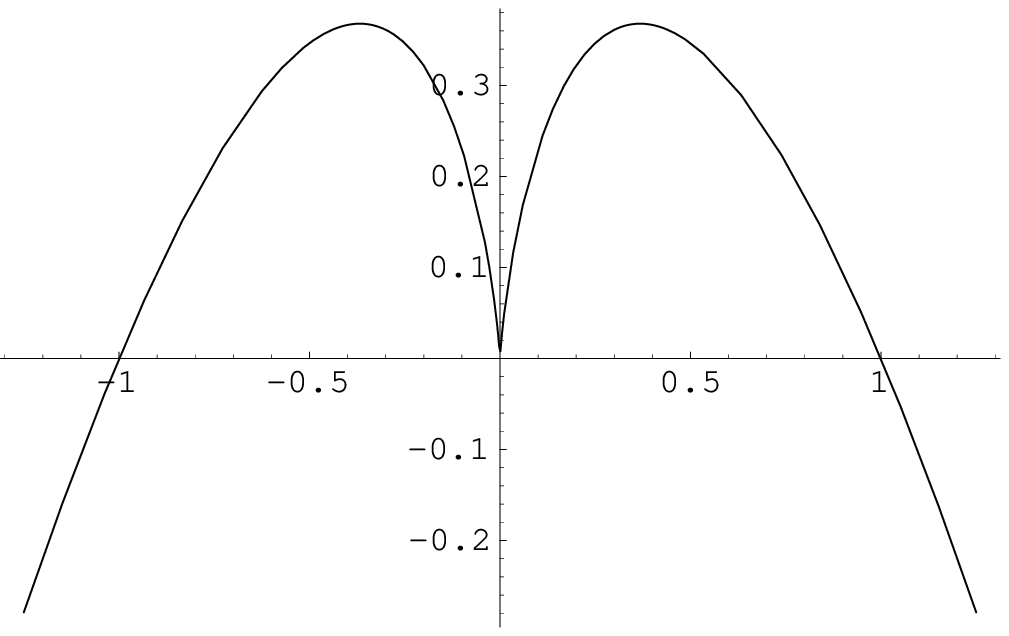}}
  \caption{Discontinuous and not-one-sided Lipschitz vector fields.
    The vector fields in (a) and (b), which differ only in their
    values at $0$, are discontinuous, and thus do not satisfy the
    hypotheses of Proposition~\ref{prop:ODE-existence}.  Therefore,
    the existence of solutions is not guaranteed. In fact, the vector
    field in (a) has no solution starting from $0$, whereas the vector
    field in (b) has a solution starting from all initial conditions.
    The vector fields in (c) and (d) are neither locally Lipschitz nor
    one-sided Lipschitz, and thus do not satisfy the hypotheses of
    Proposition~\ref{prop:ODE-uniqueness}. Therefore, uniqueness of
    solutions is not guaranteed. The vector field in (c) has two
    solutions starting from $0$. However, the vector field in (d) has
    a unique solution starting from all initial conditions.  }
  \label{fig:counter-examples}
\end{figure}
}

\clearpage

{
  \psfrag{1}[cc][cc]{{\tiny $1$}}%
  \psfrag{0.5}[cc][cc]{{\tiny $.5$}}%
  \psfrag{0}[cc][cc]{{\tiny $0$}}%
  \psfrag{-0.5}[cc][cc]{{\tiny $-.5$}}%
  \psfrag{-1}[cc][cc]{{\tiny $-1$}}%
  \psfrag{a1}[cc][cc]{{\footnotesize $x_1$}}%
  \psfrag{a2}[cc][cc]{{\footnotesize $x_2$}}%
  \begin{figure}[htbp]
    \centering
    \subfigure[]{
      \includegraphics[width=.4\linewidth]{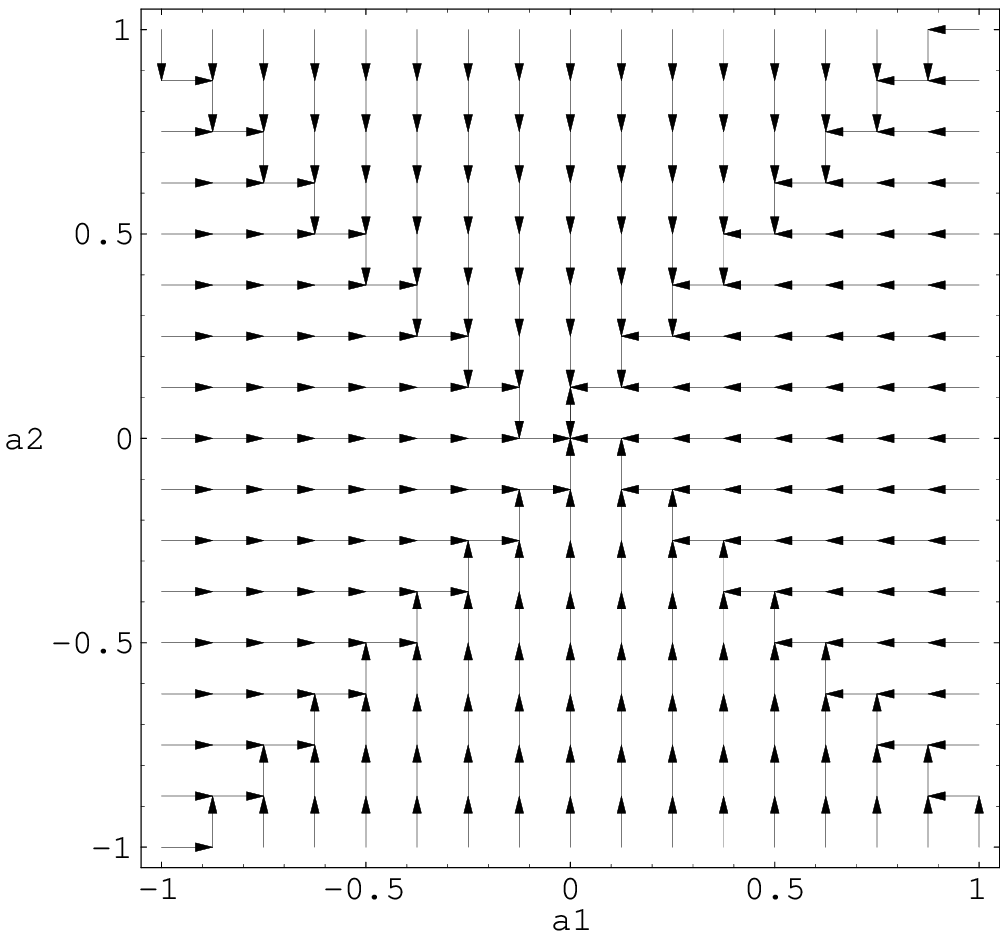}}
    \qquad
    \subfigure[]{
      \includegraphics[width=.4\linewidth]{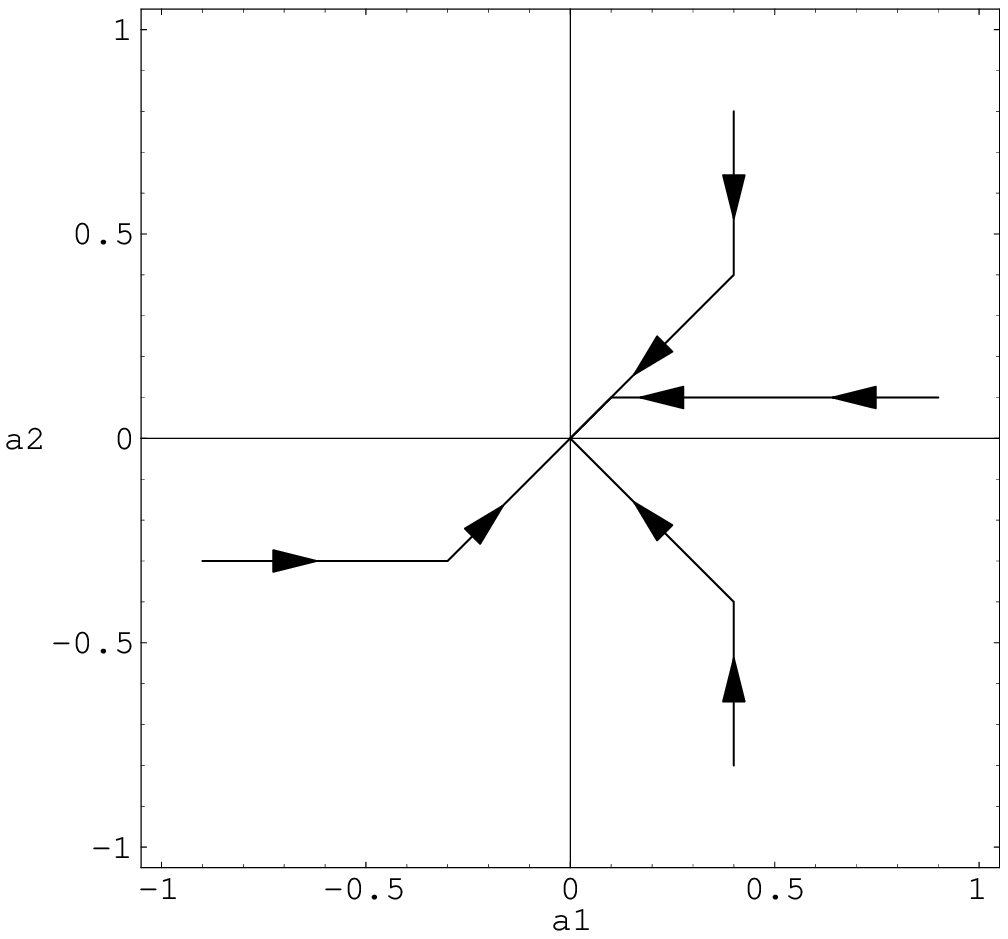}}
    \caption{Move-away-from-nearest-neighbor interaction law for one
      agent moving in the square $[-1,1]^2 \subset \real^2$. (a) shows
    the phase portrait, while (b) shows several Filippov solutions. On
    the diagonals, the vector field pushes trajectories out, whereas,
    outside the diagonals, the vector field pushes trajectories in.
    This fact makes it impossible to construct a Caratheodory solution
    starting from any initial condition on the diagonals of the
    square.  }
    \label{fig:ex1}
  \end{figure}
}

\clearpage

{
  \psfrag{0.2}[cc][cc]{{\tiny $.2$}}%
  \psfrag{0.4}[cc][cc]{{\tiny $.4$}}%
  \psfrag{0.6}[cc][cc]{{\tiny $.6$}}%
  \psfrag{0.8}[cc][cc]{{\tiny $.8$}}%
  \psfrag{1}[cc][cc]{{\tiny $1$}}%
  \psfrag{1.25}[cc][cc]{{\tiny $1.25$}}%
  \psfrag{1.5}[cc][cc]{{\tiny $1.5$}}%
  \psfrag{1.75}[cc][cc]{{\tiny $1.75$}}%
  \psfrag{2}[cc][cc]{{\tiny $2$}}%
  \psfrag{2.25}[cc][cc]{{\tiny $2.25$}}%
  \psfrag{2.5}[cc][cc]{{\tiny $2.5$}}%
  \psfrag{2.75}[cc][cc]{{\tiny $2.75$}}%
  \psfrag{a1}[cc][cc]{{\footnotesize $t$}}%
  \psfrag{a2}[cc][cc]{{\footnotesize $x$}}%
  \begin{figure}[htbp]
    \centering
    \includegraphics[width=.4\linewidth]{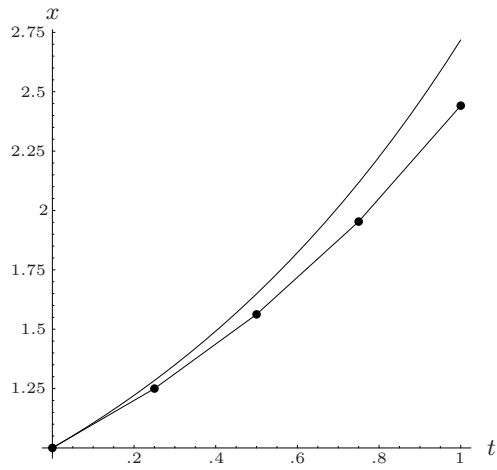}
    \caption{Illustration of the notion of sample-and-hold solution.
      For the control system $\dot x = u$, we choose the control input
      $u(x) = x$. The upper curve is the classical solution starting
      from $x_0 = 1$, while the lower curve is the sample-and-hold
      solution starting from $x_0=1$ corresponding to the
      $\pi$-partition $\{0,\frac{1}{4},\frac{1}{2}, \frac{3}{4}, 1\}$
      of $[0,1]$.}
    \label{fig:sample-and-hold}
  \end{figure}
}

\clearpage

\clearpage

{
  \psfrag{1}[cc][cc]{{\tiny $1$}}%
  \psfrag{0.5}[cc][cc]{{\tiny $.5$}}%
  \psfrag{0}[cc][cc]{{\tiny $0$}}%
  \psfrag{-0.5}[cc][cc]{{\tiny $-.5$}}%
  \psfrag{-1}[cc][cc]{{\tiny $-1$}}%
  \psfrag{a1}[cc][cc]{{\footnotesize $x_1$}}%
  \psfrag{a2}[cc][cc]{{\footnotesize $x_2$}}%
  \begin{figure}[htbp]
    \centering
    \includegraphics[width=.4\linewidth]{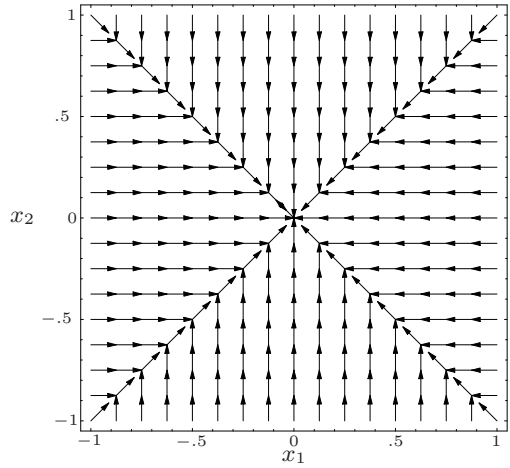}
    \caption{Generalized gradient vector field. This plot shows the
      generalized gradient vector field of the
      minimum-distance-to-polygonal-boundary function
      $\map{\sm_Q}{Q}{\real}$ on the square $[-1,1]^2$.  The vector
      field is discontinuous on the diagonals of the square.  Note the
      similarity with the phase portrait of the
      move-away-from-nearest-neighbor interaction law for one agent
      moving in the square $[-1,1]^2 \subset \real^2$ plotted in
      Figure~\ref{fig:ex1}.}
    \label{fig:sm}
  \end{figure}
}

\clearpage

\clearpage

{
  \psfrag{1}[cc][cc]{{\tiny $1$}}%
  \psfrag{0.5}[cc][cc]{{\tiny $.5$}}%
  \psfrag{0}[cc][cc]{{\tiny $0$}}%
  \psfrag{-0.5}[cc][cc]{{\tiny $-.5$}}%
  \psfrag{-1}[cc][cc]{{\tiny $-1$}}%
  \psfrag{a1}[cc][cc]{{\footnotesize $x_1$}}%
  \psfrag{a2}[cc][cc]{{\footnotesize $x_2$}}%
  \begin{figure}[htbp]
    \centering%
    \subfigure[]{
      \includegraphics[width=.355\linewidth]{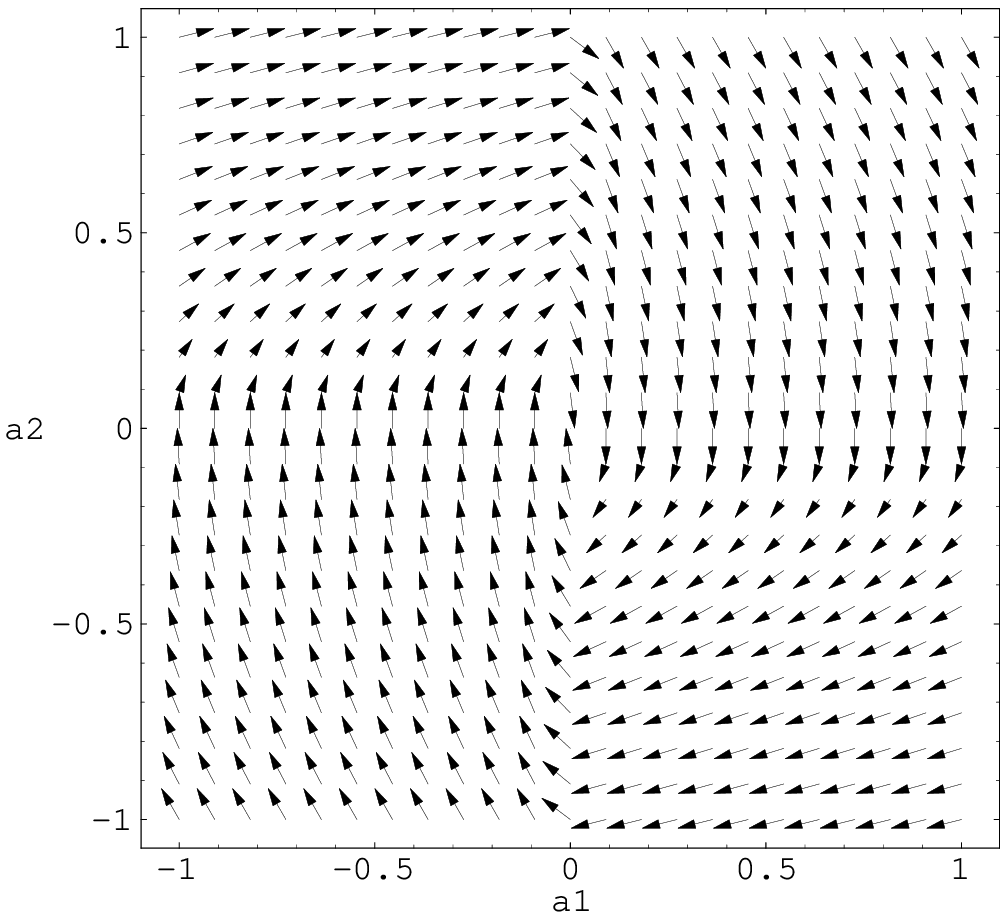}}%
    \qquad
    \subfigure[]{\includegraphics[width=.35\linewidth]{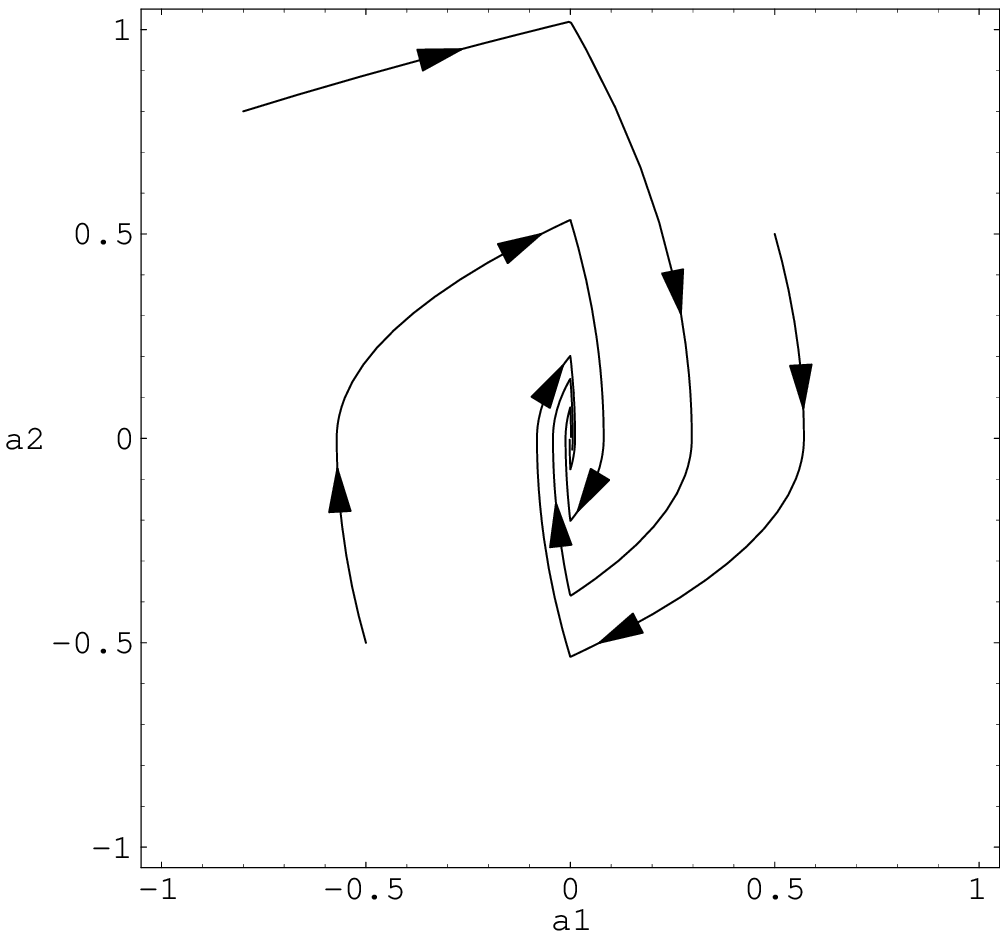}}%
    \caption{Nonsmooth harmonic oscillator with dissipation.  (a)
      shows the phase portrait on $[-1,1]^2$ of the vector field
      $(x_1,x_2) \mapsto (x_2,-\sign (x_1)- k \sign(x_2))$, where $k
      =0.75$, while (b) shows some Filippov solutions of the
      associated dynamical system~\eqref{eq:ODE-auto}. In the
      nonsmooth harmonic oscillator in Example~\ref{ex:oscillator},
      the equilibrium at the origin is strongly stable, but not
      strongly asymptotically stable, see Figure~\ref{fig:oscillator}.
      The addition of dissipation renders the equilibrium at the
      origin strongly asymptotically stable.}
    \label{fig:oscillator-dissipation}
  \end{figure}
}

\clearpage

\clearpage

\clearpage

  {
  \psfrag{3}[cc][cc]{{\tiny $3$}}%
  \psfrag{2}[cc][cc]{{\tiny $2$}}%
  \psfrag{1.5}[cc][cc]{{\tiny $1.5$}}%
  \psfrag{1}[cc][cc]{{\tiny $1$}}%
  \psfrag{0.5}[cc][cc]{{\tiny $.5$}}%
  \psfrag{0}[cc][cc]{{\tiny $0$}}%
  \psfrag{-0.5}[cc][cc]{{\tiny $-.5$}}%
  \psfrag{-1}[cc][cc]{{\tiny $-1$}}%
  \psfrag{-1.5}[cc][cc]{{\tiny $-1.5$}}%
  \psfrag{-2}[cc][cc]{{\tiny $-2$}}%
  \psfrag{-3}[cc][cc]{{\tiny $-3$}}%
  \psfrag{a1}[cc][cc]{{\footnotesize $x_1$}}%
  \psfrag{a2}[cc][cc]{{\footnotesize $x_2$}}%
  \begin{figure}[htbp]
    \centering%
    \subfigure[]{
      \includegraphics[width=.35\linewidth]{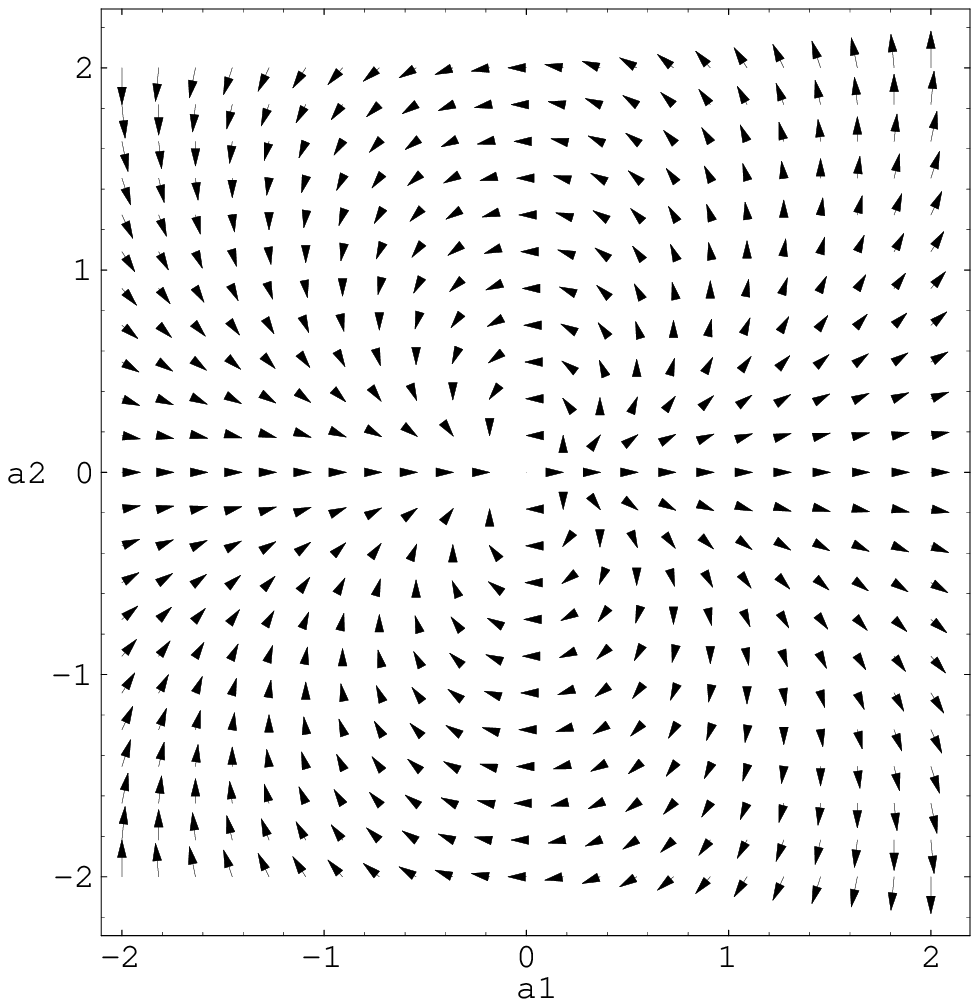}}\quad
    \subfigure[]{
      \includegraphics[width=.19\linewidth]{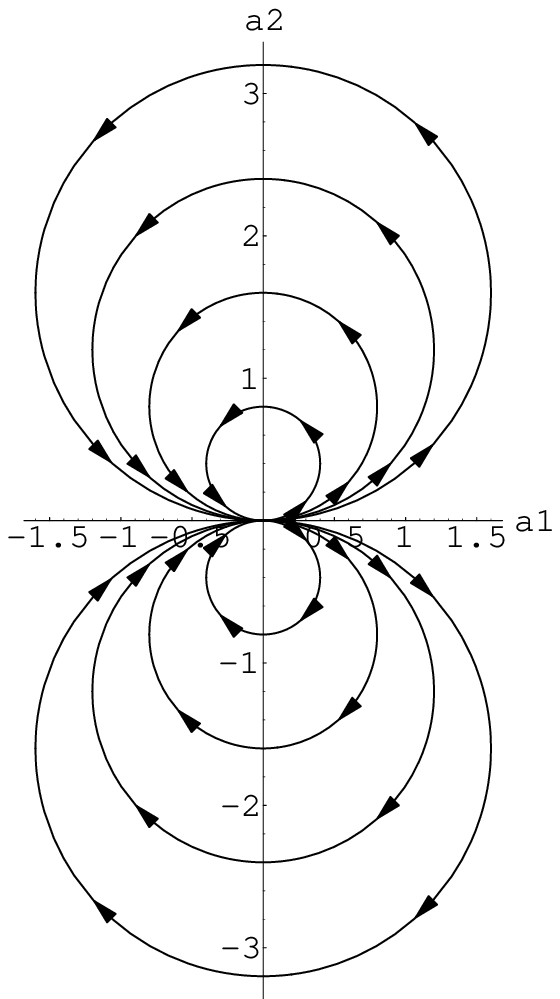}}\quad
    \subfigure[]{
      \includegraphics[width=.37\linewidth]{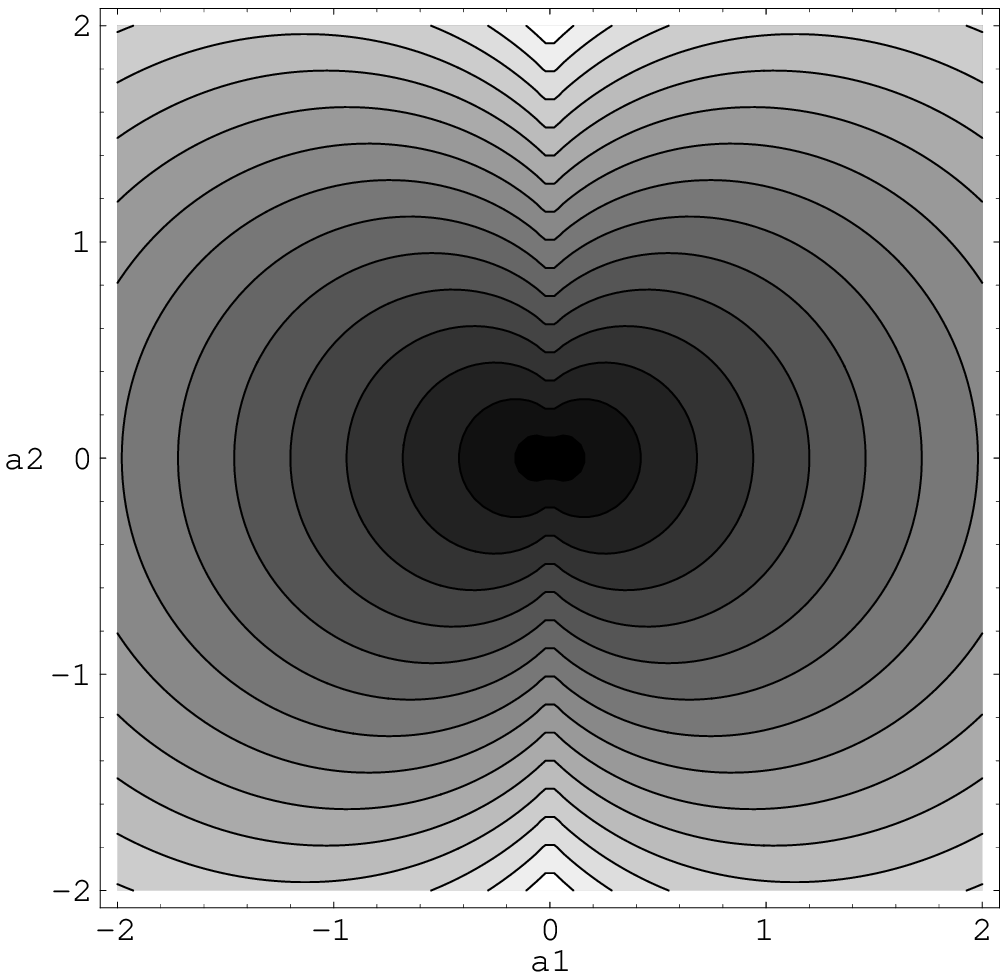}}
    \caption{Cart on a circle.  (a) shows the phase portrait of the
      input vector field $(x_1,x_2) \mapsto (x_1^2 - x_2^2 ,2 x_1
      x_2)$, (b) shows its integral curves, and (c) shows the contour
      plot of the function $0 \neq(x_1,x_2) \mapsto \frac{x_1^2 +
        x_2^2}{\sqrt{x_1^2 + x_2^2} + |x_1|}$, $(0,0) \mapsto 0$. The
      origin cannot be asymptotically stabilized by means of
      continuous feedback in the cart
      dynamics~\eqref{eq:cart-a}-\eqref{eq:cart-b}. However, the
      origin can be asymptotically stabilized by means of
      discontinuous feedback when solutions are understood in the
      sample-and-hold sense.}
    \label{fig:cart}
  \end{figure}
}%

\clearpage

\renewcommand{\thefigure}{S\arabic{figure}} 
\setcounter{figure}{0} 

\bibcite{JJM:88}{{S1}{}{{}}{{}}}
\bibcite{MDPMM:03}{{S2}{}{{}}{{}}}
\bibcite{RWC-JSP-RES:92}{{S3}{}{{}}{{}}}
\bibcite{GI:92}{{S4}{}{{}}{{}}}
\bibcite{MCF-JSP:97}{{S5}{}{{}}{{}}}
\bibcite{BB:03}{{S6}{}{{}}{{}}}
\bibcite{PJA-AN:98}{{S7}{}{{}}{{}}}
\bibcite{ASM-CCP-SSS-JMS:99}{{S8}{}{{}}{{}}}
\bibcite{AJS-HS:00}{{S9}{}{{}}{{}}}
\bibcite{RG-ART:06}{{S10}{}{{}}{{}}}

\noindent \textbf{Sidebar~\thesidebar: Solutions with Jumps}\\
\refstepcounter{sidebar}%
\renewcommand{\theequation}{S\arabic{sidebarequation}}%
\renewcommand{\theproposition}{S\arabic{sidebarproposition}}%
\noindent
In this article, we focus entirely on absolutely continuous solutions
of ordinary differential equations. However, nonsmooth continuous-time
systems can possess discontinuous solutions that admit jumps in the
state. Such notions are appropriate for dealing with, for instance,
mechanical systems subject to unilateral constraints~\cite{DES:00}.
As an example, a bouncing ball hitting the ground experiences an
instantaneous change of velocity.  This change corresponds to a
discontinuous jump in the trajectory describing the evolution of the
velocity of the ball.  Both measure differential
inclusions~\citejump{JJM:88,MDPMM:03} and linear complementarity
systems~\citejump{RWC-JSP-RES:92,GI:92,MCF-JSP:97,BB:03} are
approaches that specifically allow for discontinuous solutions.
Within hybrid systems
theory~\citejump{PJA-AN:98,ASM-CCP-SSS-JMS:99,AJS-HS:00,RG-ART:06},
discontinuous solutions arise in systems that involve both continuous-
and discrete-time evolutions.

\clearpage

\bibcite{HH:67}{{S11}{}{{}}{{}}}
\bibcite{OH:79a}{{S12}{}{{}}{{}}}
\bibcite{LA:88}{{S13}{}{{}}{{}}}
\bibcite{VAY-GAL-AKG:04}{{S14}{}{{}}{{}}}
\bibcite{JSS-DCB:96}{{S15}{}{{}}{{}}}
\bibcite{SH:91}{{S16}{}{{}}{{}}}
\bibcite{AB:04}{{S17}{}{{}}{{}}}

\noindent \textbf{Sidebar~\thesidebar: Additional Solution Notions
  for Discontinuous Systems}\\
\refstepcounter{sidebar}%
\renewcommand{\theequation}{S\arabic{sidebarequation}}%
\renewcommand{\theproposition}{S\arabic{sidebarproposition}}%
\noindent
Several solution notions are available in addition to Caratheodory,
Filippov, and sample-and-hold solutions.  These notions include the
ones considered by Krasovskii~\cite{NNK:63},
Hermes~\citesol{HH:67,OH:79a}, Ambrosio~\citesol{LA:88},
Sentis~\cite{RS:78}, and
Yakubovich-Leonov-Gelig~\citesol{VAY-GAL-AKG:04}, see
Table~\ref{tab:solutions}.  As in the case in which the vector field
is continuous, Euler solutions~\cite{AFF:88,FHC-YSL-RJS-PRW:98} are
useful for establishing existence and in characterizing basic
mathematical properties of the dynamical system.  Additional notions
of solutions for discontinuous systems are provided
in~\citesol[Section~1.1.3]{VAY-GAL-AKG:04}.  With so many notions of
solution available, various works explore the relationships among
them. For example, Caratheodory and Filippov solutions are compared
in~\citesol{JSS-DCB:96}; Caratheodory and Krasovskii solutions are
compared in~\citesol{SH:91}; Caratheodory, Euler, sample-and-hold,
Filippov, and Krasovskii solutions are compared in~\cite{FC:99};
Hermes, Filippov, and Krasovskii solutions are compared
in~\citesol{OH:79a}; and Caratheodory, Euler, and Sentis solutions are
compared in~\citesol{AB:04}.

\renewcommand{\thetable}{S\arabic{table}} 
\setcounter{table}{0} 

\begin{table}[htbp]
  \centering
  \begin{tabular}{|l|c|}
    \hline
    Notion of solution & References\\
    \hline
    Caratheodory & \cite{AFF:88}\\
    Filippov & \cite{AFF:88}\\
    Krasovskii & \cite{NNK:63}\\
    Euler & \cite{AFF:88,FHC-YSL-RJS-PRW:98} \\
    Sample-and-hold & \cite{NNK-AIS:88}\\
    Hermes & \citesol{HH:67,OH:79a}\\
    Sentis & \cite{RS:78}, \citesol{AB:04}\\
    Ambrosio & \citesol{LA:88}\\
    Yakubovich-Leonov-Gelig & \citesol{VAY-GAL-AKG:04}\\
    \hline
  \end{tabular}
  \caption{Several notions of solution for discontinuous
    dynamics. Depending on the specific problem, some notions give
    more physically meaningful solution trajectories than others.}
  \label{tab:solutions}
\end{table}

\clearpage

\bibcite{YSL-EDS:99}{{S18}{}{{}}{{}}}
\bibcite{FHC-YSL-LR-RJS:00}{{S19}{}{{}}{{}}}
\bibcite{ART-LP:00}{{S20}{}{{}}{{}}}
\bibcite{JPA:91}{{S21}{}{{}}{{}}}
\bibcite{AN-DZ:96}{{S22}{}{{}}{{}}}
\bibcite{AD-FL:92}{{S23}{}{{}}{{}}}
\bibcite{FL-VV:98}{{S24}{}{{}}{{}}}
\bibcite{VA-BB:08}{{S25}{}{{}}{{}}}
\bibcite{BB-AD-CL-VA:06}{{S26}{}{{}}{{}}}

\noindent \textbf{Sidebar~\thesidebar: Additional Topics on
  Discontinuous Systems and Differential Inclusions}\\
\refstepcounter{sidebar}%
\renewcommand{\theequation}{S\arabic{sidebarequation}}%
\renewcommand{\theproposition}{S\arabic{sidebarproposition}}%
\noindent
Beyond the topics discussed in this article, we briefly mention
several that are relevant to systems and control. These topics include
continuous dependence of solutions with respect to initial conditions
and parameters~\cite{AFF:88,JPA-AC:94}, robustness properties against
external disturbances and state measurement
errors~\cite{FHC-YSL-RJS-PRW:98} and~\citeadd{YSL-EDS:99}, conditions
for the existence of periodic solutions~\cite{AFF:88,GVS:01},
bifurcations~\cite{MDB-CJB-ARC-PK:07}, converse Lyapunov
theorems~\citeadd{FHC-YSL-LR-RJS:00,ART-LP:00}, controllability of
differential inclusions~\cite{FHC-YSL-RJS-PRW:98,GVS:01}, output
tracking in differential inclusions~\cite{JPA-HF:90}, systems subject
to gradient nonlinearities~\cite{VAY-GAL-AKG:04}, viability
theory~\citeadd{JPA:91}, maximal monotone inclusions~\cite{JPA-AC:94},
projected dynamical systems and variational
inequalities~\citeadd{AN-DZ:96}, and numerical methods for
discontinuous systems and differential
inclusions~\citeadd{AD-FL:92,FL-VV:98,VA-BB:08}.  Several works report
equivalence results among different approaches to nonsmooth systems,
including~\cite{JPA-AC:94} on the equivalence between differential
inclusions and projected dynamical systems,
and~\citeadd{BB-AD-CL-VA:06} on the equivalence of these formalisms
with complementarity systems.

\clearpage

\noindent \textbf{Sidebar~\thesidebar: Index of Symbols}\\
\refstepcounter{sidebar}%
\noindent 
The following is a list of the symbols used throughout the article.
\begin{symbollist}
\item[\textbf{Symbol}] \textbf{Description and page(s) when
    applicable} \vspace{3\lineskip}
\item[{$G[X]$}] Set-valued map associated with a control system
  $\map{X}{\real^d \times \Uc}{\real^d}$
\item[{$\overline{\co}(S)$}] Convex closure of a set $S \subseteq
  \real^d$
\item[{$\co(S)$}] Convex hull of a set $S \subseteq   \real^d$
\item[{$\diam (\pi)$}] Diameter of the partition $\pi$
\item[{$f^o(x;v)$}] Generalized directional derivative of the function
  $\map{f}{\real^d}{\real}$ at $x \in \real^d$ in the direction of $v
  \in \real^d$
\item[{$f'(x;v)$}] Right directional derivative of the function
  $\map{f}{\mathbb{R}^d}{\real}$ at $x \in \real^d$ in the direction
  of $v \in \real^d$
\item[{$S_X$}] Set of points where the vector field
  $\map{X}{\real^d}{\real^d}$ is discontinuous
\item[$\D(p,S)$] Euclidean distance from the point $p \in \real^d$ to
  the set $S \subseteq \real^d$
\item[{$F[X]$}] Filippov set-valued map associated with a vector field
  $X:\real^d \rightarrow \real^d$
\item[$\partial f$] Generalized gradient of the locally Lipschitz
  function $\map{f}{\real^d}{\real}$
\item[$\nabla f$] Gradient of the differentiable function
  $\map{f}{\real^d}{\real}$
\item[{$\LN (S)$}] Least-norm elements in the closure of the set $S
  \subseteq \real^d$
\item[{$\Omega (x)$}] Set of limit points of a curve $ t \mapsto x(t)$
\item[{$\nearestneighbor$}] Nearest-neighbor map
\item[$\normal_e$] Unit normal to the edge $e$ of a polygon $Q$
  pointing toward the interior of $Q$
\item[{$\Omega_f$}] Set of points where the locally Lipschitz function
  $\map{f}{\real^d}{\real}$ fails to be differentiable
\item[{$\pi$}] Partition of a closed interval
\item[$\parts{S}$] Set whose elements are all the possible subsets of
  $S \subseteq \real^d$
\item[$\partial_P f$] Proximal subdifferential of the lower
  semicontinuous function $\map{f}{\real^d}{\real}$
\item[$\setLiederG{\setvaluedmap}{f}$] Set-valued Lie derivative of
  the locally Lipschitz function $\map{f}{\real^d}{\real}$ with
  respect to the set-valued map
  $\map{\setvaluedmap}{\real^d}{\parts{\real^d}}$
\item[$\setLiederG{X}{f}$] Set-valued Lie derivative of the locally
  Lipschitz function $\map{f}{\real^d}{\real}$ with respect to the
  Filippov set-valued map $\map{F[X]}{\real^d}{\parts{\real^d}}$
\item[$\setLiederlow{\setvaluedmap}{f}$] Lower set-valued Lie
  derivative of the lower semicontinuous function
  $\map{f}{\real^d}{\real}$ with respect to the set-valued map
  $\map{\setvaluedmap}{\real^d}{\parts{\real^d}}$
\item[$\setLiederup{\setvaluedmap}{f}$] Upper set-valued Lie
  derivative of the lower semicontinuous function
  $\map{f}{\real^d}{\real}$ with respect to the set-valued map
  $\map{\setvaluedmap}{\real^d}{\parts{\real^d}}$
\item[{$\setvaluedmap$}] Set-valued map
\item[$\sm_Q$] Minimum distance function from a point in a convex
  polygon $Q \subset \real^d$ to the boundary of $Q$
\end{symbollist}

\clearpage

\bibcite{DMWSU:72}{{S27}{}{{}}{{}}}

\noindent \textbf{Sidebar~\thesidebar: Locally Lipschitz
  Functions}\\
\refstepcounter{sidebar}%
\noindent A function $\map{f}{\real^d}{\real^m}$ is \emph{locally
  Lipschitz at $x \in \real^d$} if there exist $L_x, \eps \in
\realpositive$ such that
\begin{align*}
  \TwoNorm{f(y) - f(y')} \le L_x \TwoNorm{y-y'} ,
\end{align*}
for all $y,y' \in \oball{\eps}{x}$.
\index{function!locally Lipschitz}
A function that is locally Lipschitz at $x$ is continuous at $x$, but
the converse is not true. For example, $\map{f}{\real}{\real}$, $f(x)
= \sqrt{|x|}$, is continuous at $0$, but not locally Lipschitz at $0$,
see Figure~\ref{fig:lispchitz-not-lipschitz}(a).
{
  \psfrag{1}[cc][cc]{{\tiny $1$}}%
  \psfrag{0.8}[cc][cc]{{\tiny $.8$}}%
  \psfrag{0.6}[cc][cc]{{\tiny $.6$}}%
  \psfrag{0.5}[cc][cc]{{\tiny $.5$}}%
  \psfrag{0.4}[cc][cc]{{\tiny $.4$}}%
  \psfrag{0.2}[cc][cc]{{\tiny $.2$}}%
  \psfrag{0}[cc][cc]{{\tiny $0$}}%
  \psfrag{-0.5}[cc][cc]{{\tiny $-.5$}}%
  \psfrag{-1}[cc][cc]{{\tiny $-1$}}%
  \begin{figure}[htbp]
    \centering%
    \subfigure[]{\includegraphics[width=.4\linewidth]{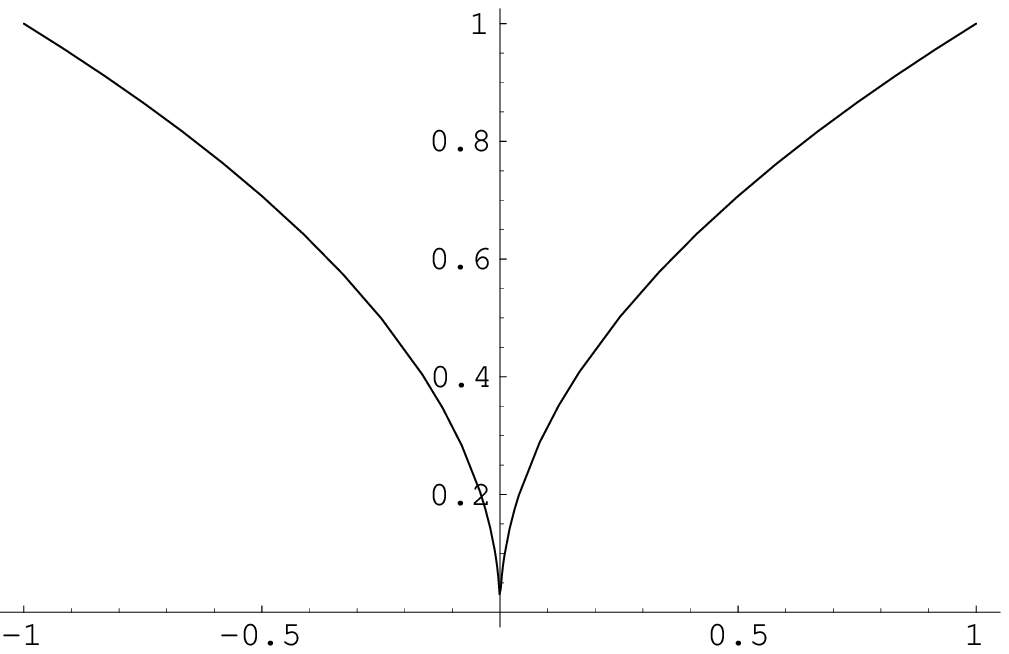}}\quad
    \subfigure[]{\includegraphics[width=.4\linewidth]{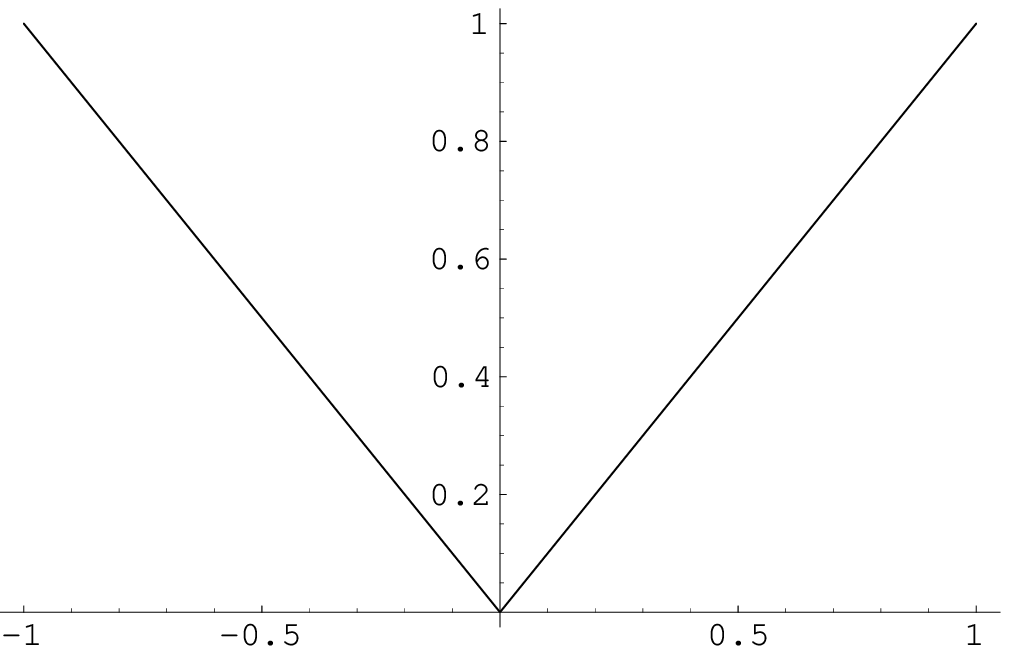}}
    \caption{Illustration of the difference among continuous, locally
      Lipschitz, and differentiable functions.  (a) shows the graph of
      $\map{f}{\real}{\real}$, $f(x) = \sqrt{|x|}$, which is continuous
      at $0$, but not locally Lipschitz at $0$. (b) shows the graph of
      $\map{f}{\real}{\real}$, $f(x) = |x|$, which is is locally
      Lipschitz at $0$, but not differentiable at $0$.}
    \label{fig:lispchitz-not-lipschitz}
\end{figure}
}%
A function is \emph{locally Lipschitz on $S \subseteq \real^d$} if it
is locally Lipschitz at $x$ for all $x \in S$. If $f$ is locally
Lipschitz on $\real^d$, we simply say $f$ is locally Lipschitz.
Convex functions are locally Lipschitz~\citelip{DMWSU:72}, and hence
concave functions are also locally Lipschitz. Note that a function
that is continuously differentiable at $x$ is locally Lipschitz at
$x$, but the converse is not true. For example,
$\map{f}{\real}{\real}$, $f(x) = |x|$, is locally Lipschitz at $0$,
but not differentiable at $0$, see
Figure~\ref{fig:lispchitz-not-lipschitz}(b).  A function
$\map{f}{\real \times \real^d}{\real^m}$ that depends explicitly on
time is \emph{locally Lipschitz at $x \in \real^d$} if there exists
$\eps \in \realpositive$ and $\map{L_x}{\real}{\realpositive}$ such
that $\TwoNorm{f(t, y) - f(t, y')} \le L_x(t) \TwoNorm{y-y'}$ for all
$t \in \real$ and $y,y' \in \oball{\eps}{x}$.

\clearpage

\bibcite{DCB-PAB:97}{{S28}{}{{}}{{}}}

\noindent \textbf{Sidebar~\thesidebar: Caratheodory
  Conditions for Time-varying Vector Fields}\\
\refstepcounter{sidebar}%
\renewcommand{\theequation}{S\arabic{sidebarequation}}%
\renewcommand{\theproposition}{S\arabic{sidebarproposition}}%
\noindent

Consider the differential equation
\begin{align}
  \label{eq:ODE-nonauto}
  \dot x (t) = X(t,x(t)) ,
\end{align}
where $\map{X}{\realnonnegative \times \real^d}{\real^d}$ is a
time-varying vector field.  The following result is taken
from~\cite{AFF:88}: a weaker version of the Caratheodory conditions is
given, for instance, in~\citetime{DCB-PAB:97}.

\begin{proposition}\label{prop:Caratheodory-existence}
  Let $\map{X}{\realnonnegative \times \real^d}{\real^d}$. Assume that
  \textit{(i)} for almost all $t \in \realnonnegative$, the map $x
  \mapsto X(t,x)$ is continuous, \textit{(ii)} for each $x \in
  \real^d$, the map $t \mapsto X(t,x)$ is measurable, and
  \textit{(iii)} $X$ is locally essentially bounded, that is,
  \index{vector field!locally essentially bounded}%
  for all $(t,x) \in \realnonnegative \times \real^d$, there exist
  $\eps \in \realpositive$ and an integrable function
  $\map{m}{[t,t+\delta]}{\realpositive}$ such that $\TwoNorm{X(s,y)}
  \le m(s)$ for almost all $s \in [t,t+\delta]$ and all $y \in
  \oball{\eps}{x}$.  Then, for all $(t_0,x_0) \in \realnonnegative
  \times \real^d$, there exists a Caratheodory solution
  of~\eqref{eq:ODE-nonauto} with initial condition $x(t_0) = x_0$.
\end{proposition}

The specialization of Proposition~\ref{prop:Caratheodory-existence} to
a time-invariant vector field requires that the vector field be
continuous, which in turn guarantees the existence of a classical
solution.

\clearpage

\sindex{Parts}{$\parts{S}$}{Set whose elements are all the possible
  subsets of $S \subseteq \real^d$}%
\sindex{SetValuedMap}{{$\setvaluedmap$}}{Set-valued map}%
\noindent \textbf{Sidebar~\thesidebar: Set-valued Maps}\\
\refstepcounter{sidebar}%
\noindent
A set-valued map, as its name suggests, is a map that assigns sets to
points. We consider time-varying set-valued maps of the form
$\map{\setvaluedmap}{\realnonnegative \times
  \real^d}{\parts{\real^d}}$.  Recall that $\parts{\real^d}$ denotes
the collection of all subsets of $\real^d$.  The map $\setvaluedmap$
assigns to each point $(t,x) \in \realnonnegative \times \real^d$ the
set $\setvaluedmap (t,x) \subseteq \real^d$.  Note that a standard map
$\map{f}{\realnonnegative \times \real^d}{\real^d}$ can be interpreted
as a singleton-valued map.  A complete analysis for set-valued maps
can be developed, as in the case of standard maps~\cite{JPA-HF:90}.
Here, we are mainly interested in concepts related to boundedness and
continuity, which we define next.
  
The set-valued map $\map{\setvaluedmap}{\realnonnegative \times
  \real^d}{\parts{\real^d}}$ is \emph{locally bounded} (respectively,
\emph{locally essentially bounded})
\index{set-valued map!locally bounded}%
\index{set-valued map!locally essentially bounded}%
at $(t,x) \in \realnonnegative \times \real^d$ if there exist $\eps,
\delta \in \realpositive$ and an integrable function
$\map{m}{[t,t+\delta]}{\realpositive}$ such that $\TwoNorm{z} \le
m(s)$ for all $z \in \setvaluedmap (s,y)$, all $s \in [t,t+\delta]$,
and all $y \in \oball{\eps}{x}$ (respectively, almost all $y \in
\oball{\eps}{x}$ in the sense of Lebesgue measure).

The time-invariant set-valued map
$\map{\setvaluedmap}{\real^d}{\parts{\real^d}}$ is \emph{upper
  semicontinuous}\index{set-valued map!upper semicontinuous}
(respectively, \emph{lower semicontinuous}\index{set-valued map!lower
  semicontinuous}) at $x \in \real^d$ if, for all $\eps \in
\realpositive$, there exists $\delta \in \realpositive$ such that
$\setvaluedmap (y) \subseteq \setvaluedmap (x) + \oball{\eps}{0}$
(respectively, $\setvaluedmap (x) \subseteq \setvaluedmap (y) +
\oball{\eps}{0}$) for all $y \in B(x,\delta)$.  The set-valued map
$\map{\setvaluedmap}{\real^d}{\parts{\real^d}}$ is
\emph{continuous}\index{set-valued map!continuous} at $x \in \real^d$
if it is both upper and lower semicontinuous at $x \in \real^d$.
Finally, the set-valued map
$\map{\setvaluedmap}{\real^d}{\parts{\real^d}}$ is \emph{locally
  Lipschitz}\index{set-valued map!locally Lipschitz} at $x \in
\real^d$ if there exist $L_x, \eps \in \realpositive$ such that
\begin{align*}
  \setvaluedmap (y') \subseteq \setvaluedmap (y) + L_x \TwoNorm{y -
    y'} \cball{1}{0} ,
\end{align*}
for all $y,y' \in \oball{\eps}{x}$. A locally Lipschitz set-valued map
at $x$ is upper semicontinuous at $x$, but the converse is not true.

The notion of upper semicontinuity of a map $\map{f}{\real^d}{\real}$,
which is defined in the section ``The proximal subdifferential of a
lower semicontinuous function,'' is weaker than the notion of upper
semicontinuity of $f$ when viewed as a (singleton-valued) set-valued
map from $\real^d$ to $\parts{\real}$.  Indeed, the latter is
equivalent to the condition that $\map{f}{\real^d}{\real}$ is
continuous.

\clearpage

\noindent \textbf{Sidebar~\thesidebar: Caratheodory
  Solutions of Differential Inclusions}\\
\refstepcounter{sidebar}%
\renewcommand{\theequation}{S\arabic{sidebarequation}}%
\renewcommand{\theproposition}{S\arabic{sidebarproposition}}%
\noindent
A differential inclusion~\cite{GVS:01,JPA-AC:94} is a generalization
of a differential equation. At each state, a differential inclusion
specifies a set of possible evolutions, rather than a single one.
This object is defined by means of a set-valued map, see ``Set-valued
Maps.''  The \emph{differential inclusion}\index{differential
  inclusion} associated with a time-varying set-valued map
$\map{\setvaluedmap}{\realnonnegative \times
  \real^d}{\parts{\real^d}}$ is an equation of the form
\begin{align}\label{eq:diff-inclusion}
  \dot{x}(t) \in \setvaluedmap (t,x(t)) .
\end{align}
\addtocounter{sidebarequation}{1}%
\addtocounter{equation}{-1}%
The point $x_e \in \real^d$ is an \emph{equilibrium} of the
differential inclusion if $0 \in \setvaluedmap(t,x_e)$ for all $t \in
\realnonnegative$.  We now define the notion of solution of a
differential inclusion in the sense of Caratheodory.

A \emph{Caratheodory solution
  of~\eqref{eq:diff-inclusion}}\index{differential inclusion!solution}
defined on $[t_0,t_1] \subset \realnonnegative$ is an absolutely
continuous map $\map{x}{[t_0,t_1]}{\real^d}$ such that $\dot{x}(t) \in
\setvaluedmap (t,x(t))$ for almost every $t \in [t_0,t_1]$.  The
existence of at least one solution starting from each initial
condition is guaranteed by the following result (see, for
instance,~\cite{JPA-AC:94,AB-LR:05}).

\begin{proposition}\label{prop:existence-solution}
  Let $\map{\setvaluedmap}{\realnonnegative \times
    \real^d}{\parts{\real^d}}$ be locally bounded and take nonempty,
  compact, and convex values.  Assume that, for each $t \in \real$,
  the set-valued map $x \mapsto \setvaluedmap (t,x)$ is upper
  semicontinuous, and, for each $x \in \real^d$, the set-valued map $t
  \mapsto \setvaluedmap (t,x)$ is measurable.  Then, for all
  $(t_0,x_0) \in \realnonnegative \times \real^d$, there exists a
  Caratheodory solution of~\eqref{eq:diff-inclusion} with initial
  condition $x(t_0) = x_0$.
\end{proposition}
\addtocounter{sidebarproposition}{1}%
\addtocounter{proposition}{-1}%

This result is sufficient for our purposes.  Additional existence
results based on alternative assumptions are given
in~\cite{AB-LR:05,FHC-YSL-RJS-PRW:98}.  The uniqueness of Caratheodory
solutions is guaranteed by the following result.
\begin{proposition}\label{prop:unique-sol-inclusion}
  In addition to the hypotheses of
  Proposition~\ref{prop:existence-solution}, assume that, for all $x
  \in \real^d$, there exist $\eps \in \realpositive$ and an integrable
  function $\map{L_x}{\real}{\realpositive}$ such that
  \begin{align}\label{eq:Lipschitz-like}
    (v - w)^T (y-y') \le L_x(t) \, \TwoNorm{y-y'}^2 ,
  \end{align}
  \addtocounter{sidebarequation}{1}%
  \addtocounter{equation}{-1}%
  for almost every $y,y' \in \oball{\eps}{x}$, every $t \in
  \realnonnegative$, every $v \in \setvaluedmap (t,y)$, and every $w
  \in \setvaluedmap (t,y')$.
  Then, for all $(t_0,x_0) \in \realnonnegative \times \real^d$, there
  exists a unique Caratheodory solution of~\eqref{eq:diff-inclusion}
  with initial condition $x(t_0) = x_0$.
\end{proposition}
\addtocounter{sidebarproposition}{1}%
\addtocounter{proposition}{-1}%

The following example illustrates
propositions~\ref{prop:existence-solution}
and~\ref{prop:unique-sol-inclusion}.  Following~\cite{JPA-HF:90},
consider the set-valued map
$\map{\setvaluedmap}{\real}{\parts{\real}}$ defined by
\begin{align*}
  \setvaluedmap(x) =
  \begin{cases}
    0, & x \neq 0 ,\\
    [-1,1], & x = 0 .
  \end{cases}
\end{align*}
Note that $\setvaluedmap$ is upper semicontinuous, but not lower
semicontinuous, and thus it is not continuous. This set-valued map
satisfies all of the hypotheses in
Proposition~\ref{prop:existence-solution}, and therefore Caratheodory
solutions exist starting from all initial conditions. In addition,
$\setvaluedmap$ satisfies~\eqref{eq:Lipschitz-like} as long as $y$ and
$y'$ are nonzero. Therefore,
Proposition~\ref{prop:unique-sol-inclusion} guarantees the uniqueness
of Caratheodory solutions.  In fact, for every initial condition, the
Caratheodory solution of $ \dot{x}(t) \in \setvaluedmap (x(t))$ is
just the equilibrium solution.

\clearpage

\begin{figure}[htbp]
  \centering
  \resizebox{0.4\linewidth}{!}{\input{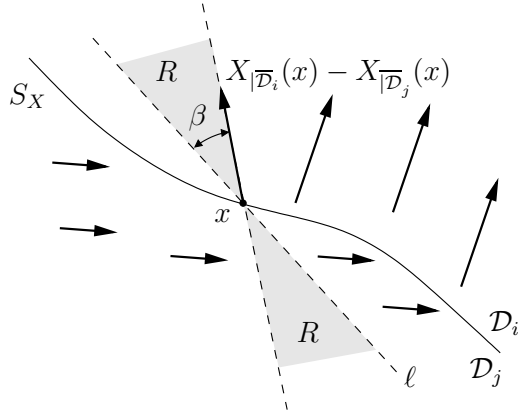}}
  \caption{Piecewise continuous vector field.  The vector field has a
    unique Filippov solution starting from all initial conditions.
    Note that solutions that reach $S_X$ coming from $\domain_j$ cross
    it, and then continue in $\domain_i$.  However, the vector field
    is not essentially one-sided Lipschitz, and hence
    Proposition~\ref{prop:uniqueness-Filippov-I} cannot be invoked to
    conclude uniqueness.}
  \label{fig:Lipschitz-type-counter-example}
\end{figure}

\noindent \textbf{Sidebar~\thesidebar: Uniqueness of Filippov
  Solutions of Piecewise Continuous Vector Fields}\\
\refstepcounter{sidebar}%
\noindent
Here we justify why, in order to guarantee the uniqueness of Filippov
solutions for a piecewise continuous vector field, we cannot resort to
Proposition~\ref{prop:uniqueness-Filippov-I} and instead must use
Proposition~\ref{prop:uniqueness-Filippov-II}.  To see this, consider
a piecewise continuous vector field $\map{X}{\real^d}{\real^d}$, $d
\ge 2$, and let $x \in S_X$ be a point of discontinuity.  Let us show
that $X$ is not essentially one-sided Lipschitz on every neighborhood
of $x$. For simplicity, assume $x$ belongs to the boundaries of just
two sets, that is, $x \in \boundary{\domain_i} \intersection
\boundary{\domain_j}$ (the argument proceeds similarly for the general
case). For $\eps>0$, we show that~\eqref{eq:unique-Filipov-sol-I} is
violated on a set of nonzero measure contained in $\oball{\eps}{x}$.
Notice that
\begin{align*}
  (X(y) - X(y'))^T (y-y') = \TwoNorm{X(y) - X(y')} \, \TwoNorm{y-y'}
  \, \cos \alpha (y,y') ,
\end{align*}
where $ \alpha (y,y') = \angle (X(y) - X(y'),y - y')$ is the angle
between the vectors $X(y) - X(y')$ and $y-y'$.
Therefore,~\eqref{eq:unique-Filipov-sol-I} is equivalent to
\begin{align}\label{eq:unique-Filipov-sol-I-equiv}
  \TwoNorm{X(y) - X(y')} \, \cos \alpha (y,y') \le L \TwoNorm{y - y '}
  .
\end{align}
\addtocounter{sidebarequation}{1}%
\addtocounter{equation}{-1}%
Consider the vectors $X_{|\ov{\domain_i}}(x)$ and $
X_{|\ov{\domain_j}}(x)$.  Since $X$ is discontinuous at $x$, we have
$X_{|\ov{\domain_i}}(x) \neq X_{|\ov{\domain_j}}(x)$.  Take $y \in
\domain_i \intersection \oball{\eps}{x}$ and $y' \in \domain_j
\intersection \oball{\eps}{x}$.  Note that, as $y$ and $y'$ tend to
$x$, the vector $X(y) - X(y')$ tends to $X_{|\ov{\domain_i}}(x) -
X_{|\ov{\domain_j}}(x)$.  Consider then the straight line $\ell$ that
crosses $S_X$, passes through $x$, and forms a small angle $\beta>0$
with $X_{|\ov{\domain_i}}(x) - X_{|\ov{\domain_j}}(x)$, see
Figure~\ref{fig:Lipschitz-type-counter-example}.

Let $R$ be the set enclosed by the line $\ell$ and the line in the
direction of the vector $X_{|\ov{\domain_i}}(x) -
X_{|\ov{\domain_j}}(x)$.  If $y \in \domain_i \intersection R$ and $y'
\in \domain_j \intersection R$ tend to $x$, we deduce that $\TwoNorm{y
  - y'} \rightarrow 0$ while, at the same time,
\begin{multline*}
  \TwoNorm{X(y) - X(y')} \, |\cos \alpha (y,y')| \\
  \ge \TwoNorm{X(y) - X(y')} \, \cos \beta \longrightarrow
  \TwoNorm{X_{|\ov{\domain_i}}(x) - X_{|\ov{\domain_j}}(x)} \, \cos
  \beta > 0 .
\end{multline*}
Therefore, there does not exist $L \in \realpositive$ such
that~\eqref{eq:unique-Filipov-sol-I-equiv} is satisfied for $y \in R
\intersection \domain_i \cap \oball{\eps}{x}$ and $y' \in R
\intersection \domain_j \cap \oball{\eps}{x}$. Thus, $X$ is not
essentially one-sided Lipschitz on any neighborhood of $x$, and the
hypotheses of Proposition~\ref{prop:uniqueness-Filippov-I} do not
hold.

\clearpage

\sindex{DirectionalDerivativeRight}{{$f'(x;v)$}}{Right directional
  derivative of the function $\map{f}{\mathbb{R}^d}{\real}$ at $x
  \in \real^d$ in the direction of $v \in \real^d$}
\sindex{DirectionalDerivativeGeneralized}{{$f^o(x;v)$}}{Generalized
  directional derivative of the function $\map{f}{\real^d}{\real}$ at
  $x \in \real^d$ in the direction of $v \in \real^d$}
\noindent \textbf{Sidebar~\thesidebar: Regular Functions}\\
\refstepcounter{sidebar}%
\noindent
To introduce the notion of regular function, we need to first define
the right directional derivative and the generalized right directional
derivative. Given $\map{f}{\real^d}{\real}$, the \emph{right
  directional derivative} of $f$ at $x$ in the direction of $v \in
\real^d$ is defined as
\begin{align*}
  f'(x;v) = \lim_{h \rightarrow 0^+} \frac{f(x+hv) - f(x)}{h} ,
\end{align*}
when this limits exists. On the other hand, the \emph{generalized
  directional derivative} of $f$ at $x$ in the direction of $v \in
\real^d$ is defined as
\begin{align*}
  f^o(x;v) = \limsup_{
    \begin{subarray}{l}
      y \rightarrow x\\
      h \rightarrow 0^+
    \end{subarray}
  } \frac{f(y+hv)-f(y)}{h} = \lim_{
    \begin{subarray}{l}
      \delta \rightarrow 0^+\\
      \eps \rightarrow 0^+
    \end{subarray}
  } \sup_{
    \begin{subarray}{l}
      y\in \oball{\delta}{x}\\
      h \in [0,\eps)
    \end{subarray}
  } \frac{f(y+hv)-f(y)}{h} .
\end{align*}
The advantage of the generalized directional derivative compared to
the right directional derivative is that the limit always exists. When
the right directional derivative exists, these quantities may be
different. When they are equal, we call the function regular. More
formally, a function $\map{f}{\real^d}{\real}$ is \emph{regular at
  $x\in \real^d$} if, for all $v \in \real^d$, the right directional
derivative of $f$ at $x$ in the direction of $v$ exists, and $f'(x;v)
= f^o(x;v)$.
\index{function!regular}
A function that is continuously differentiable at $x$ is regular at
$x$.  Also, a convex function is regular
(cf.~\cite[Proposition~2.3.6]{FHC:83}).

The function $\map{g}{\real}{\real}$, $g(x) = -|x|$, is not regular.
Since $g$ is continuously differentiable everywhere except for zero,
it is regular on $\real \! \setminus \! \{0\}$. However, its directional
derivatives
\begin{align*}
  g'(0;v)  = 
  \begin{cases}
    -v , & v>0 , \\
    \phantom{-}v , & v<0,
  \end{cases}
  \qquad
  g^o(0;v) =
  \begin{cases}
    \phantom{-}v , & v>0 , \\
    -v , & v<0,
  \end{cases}
\end{align*}
do not coincide. Hence, $g$ is not regular at $0$.

\clearpage

\textbf{Jorge Cort\'es (Department of Mechanical and Aerospace
  Engineering, University of California at San Diego, 9500 Gilman Dr,
  La Jolla, CA 92093, phone 1-858-822-7930, fax 1-858-822-3107,
  cortes@ucsd.edu)} received the Licenciatura degree in mathematics
from the Universidad de Zaragoza, Spain, in 1997, and his Ph.D. degree
in engineering mathematics from the Universidad Carlos III de Madrid,
Spain, in 2001.  He held postdoctoral positions at the Systems,
Signals, and Control Department of the University of Twente, and at
the Coordinated Science Laboratory of the University of Illinois at
Urbana-Champaign.  From 2004 to 2007, he was an assistant professor
with the Department of Applied Mathematics and Statistics, University
of California, Santa Cruz.  He is currently an assistant professor in
the Department of Mechanical and Aerospace Engineering, University of
California, San Diego.  His research interests focus on mathematical
control theory, distributed motion coordination for groups of
autonomous agents, and geometric mechanics and geometric integration.
He is the author of \emph{Geometric, Control and Numerical Aspects of
  Nonholonomic Systems} (Springer Verlag, 2002), and the recipient of
the 2006 Spanish Society of Applied Mathematics Young Researcher
Prize. He is currently an associate editor for the \textit{European
  Journal of Control.}

\end{document}